%% file: aap-sample.tex
\documentclass[aap]{imsart}

\input{math}

\RequirePackage{amsthm,amsmath,amsfonts,amssymb}
\RequirePackage[numbers]{natbib}
\RequirePackage[colorlinks,citecolor=blue,urlcolor=blue]{hyperref}
\RequirePackage{graphicx}

\startlocaldefs
\theoremstyle{plain}
\newtheorem{theorem}{Theorem}
\newtheorem{lemma}{Lemma}

\newtheorem{remark}{Remark}
\newtheorem{definition}{Definition}
\newtheorem{assumption}{Assumption}
\newtheorem{corollary}{Corollary}
\usepackage{todonotes}
\usepackage{listings}
\usepackage{algorithm}
\usepackage{algpseudocode}

\makeatletter
\def\tagform@#1{\maketag@@@{(\ignorespaces#1\unskip\@@italiccorr)}}
\makeatother

\let\originaleqref=\ref

\renewcommand*{\eqref}[1]{Eq.~(\originaleqref{#1})}

\DeclareMathOperator{\EE}{\mathbb{E}}
\DeclareMathOperator{\supp}{\mathcal{S}}
\DeclareMathOperator{\dist}{dist}
\DeclareMathOperator{\tr}{tr}

\DeclareMathOperator{\mat}{Mat}
\DeclareMathOperator{\asto}{\,\xrightarrow{\text{a.s.}}\,}

\usepackage{comment}
\definecolor{cblue}{rgb}{0,0,0}
\endlocaldefs

\begin{document}

\begin{frontmatter}
\title{When random tensors meet random matrices}
\runtitle{When random tensors meet random matrices}

\begin{aug}
\author[A]{\fnms{Mohamed El Amine} \snm{Seddik}\ead[label=e1]{mohamed.seddik@huawei.com}},
\author[A]{\fnms{Maxime} \snm{Guillaud}\ead[label=e2]{maxime.guillaud@huawei.com}}
\and
\author[B]{\fnms{Romain} \snm{Couillet}\ead[label=e3]{romain.couillet@univ-grenoble-alpes.fr}}
\address[A]{Mathematical and Algorithmic Sciences Laboratory, Huawei Technologies France,\\
\printead{e1,e2}}

\address[B]{Université Grenoble Alpes, CNRS, GIPSA-lab, Grenoble INP,\\
\printead{e3}}
\end{aug}

\begin{abstract}
Relying on random matrix theory (RMT), this paper studies asymmetric order-$d$ spiked tensor models with Gaussian noise. Using the variational definition of the singular vectors and values of (Lim, 2005) \cite{lim2005singular}, we show that the analysis of the considered model boils down to the analysis of an equivalent spiked symmetric \textit{block-wise} random matrix, that is constructed from \textit{contractions} of the studied tensor with the singular vectors associated to its best rank-1 approximation. Our approach allows the exact characterization of the almost sure asymptotic singular value and alignments of the corresponding singular vectors with the true spike components, when $\frac{n_i}{\sum_{j=1}^d n_j}\to c_i\in (0, 1)$ with $n_i$'s the tensor dimensions. In contrast to other works that rely mostly on tools from statistical physics to study random tensors, our results rely solely on classical RMT tools such as Stein's lemma. Finally, classical RMT results concerning spiked random matrices are recovered as a particular case. 
\end{abstract}




\end{frontmatter}
\tableofcontents

\paragraph*{Notations:} $[n]$ denotes the set $\{1, \ldots, n\}$. The set of rectangular matrices of size $m\times n$ is denoted $\sM_{m,n}$. The set of square matrices of size $n$ is denoted $\sM_n$. The set of $d$-order tensors of size $n_1\times \cdots \times n_d$ is denoted $\sT_{n_1, \ldots, n_d}$. The set of hyper-cubic tensors of size $n$ and order $d$ is denoted $\sT_n^d$. The notation $\tX\sim \sT_{n_1, \ldots, n_d}(\mathcal{N}(0, 1))$ means that $\tX$ is a random tensor with i.i.d.\@ Gaussian $\mathcal{N}(0, 1)$ entries. Scalars are denoted by lowercase letters as $a,b,c$. Vectors are denoted by bold lowercase letters as $\va,\vb,\vc$. $\ve_i^d$ denotes the canonical vector in $\sR^d$ with $[\ve_i^d]_j = \delta_{ij}$. Matrices are denoted by bold uppercase letters as $\mA,\mB,\mC$. Tensors are denoted as $\tA, \tB, \tC$. $T_{i_1,\ldots, i_d}$ denotes the entry $(i_1,\ldots, i_d)$ of the tensor $\tT$. $\langle \vu, \vv\rangle = \sum_{i} u_i v_i$ denotes the scalar product between $\vu$ and $\vv$, the $\ell_2$-norm of a vector $\vu$ is denoted as $\Vert \vu \Vert^2 = \langle \vu, \vu\rangle$. $\Vert \cdot \Vert$ denotes the spectral norm for tensors. $\tT(\vu^{(1)},\ldots,\vu^{(d)}) \equiv \sum_{i_1, \ldots, i_d} u_{i_1}^{(1)}\ldots u_{i_d}^{(d)} T_{i_1,\ldots, i_d}$ denotes the contraction of tensor $\tT$ on the vectors given as arguments. Given some vectors $\vu^{(1)},\ldots,\vu^{(k)}$ with $k<d$, the contraction $\tT(\vu^{(1)},\ldots,\vu^{(k)},:,\ldots,:)$ denotes the resulting $(d-k)$-th order tensor. $\sS^{N-1}$ denotes the $N$-dimensional unit sphere.

\section{Introduction}\label{sec:intro}
The extraction of latent and low-dimensional structures from raw data is a key step in various machine learning and signal processing applications. Our present interest is in those modern techniques which rely on the extraction of such structures from a low-rank random tensor model \cite{anandkumar2014tensor} and which extends the ideas from matrix-type data to tensor structured data. We refer the reader to \cite{sun2014tensors, rabanser2017introduction, zare2018extension} and the references therein which introduce an extensive set of applications of tensor decomposition methods to machine learning, including dimensionality reduction, supervised and unsupervised learning, learning subspaces for feature extraction, low-rank tensor recovery etc. Although random matrix models have been extensively studied and well understood in the literature, the understanding of random tensor models is still in its infancy and the ideas from random matrix analysis do not easily extend to higher-order tensors. Indeed, the \textit{resolvent} notion (see Definition \ref{def:resolvent}) which is at the heart of random matrices does not generalize to tensors.
In our present investigation, we consider the spiked tensor model, which consists in an observed $d$-order tensor $\tT \in \sT_{n_1,\ldots, n_d}$ of the form
\begin{align} 
\label{eq_asymmetric_spike_model}
\tT =  \beta\, \vx^{(1)} \otimes \cdots \otimes \vx^{(d)} + \frac{1}{\sqrt N} \tX,
\end{align}
where $(\vx^{(1)} , \ldots , \vx^{(d)})\in \sS^{n_1 - 1} \times \cdots \times \sS^{n_d - 1}$, $\tX \sim \sT_{n_1,\ldots, n_d}(\mathcal{N}(0,1))$ and $N = \sum_{i=1}^d n_i$. {\color{cblue} Note that the tensor noise is normalized by $\sqrt{\sum_{i=1}^d n_i}$ with $n_i$ the tensor dimensions, since the spectral norm of $\tX$ is of order $\sqrt{\sum_{i=1}^d n_i}$ from Lemma \ref{lemma_spectral_norm_of_X}.}
One aims at retrieving the rank-1 component (or \textit{spike}) $\beta\, \vx^{(1)} \otimes \cdots \otimes \vx^{(d)}$ from the noisy tensor $\tT$, where $\beta$ can be seen as controlling the signal-to-noise ratio (SNR).
The identification of the dominant rank-1 component is an important special case of the low-rank tensor approximation problem, a noisy version of the classical canonical polyadic decomposition (CPD) \cite{hitchcock1927expression, landsberg2012tensors}.
Extensive efforts have been made to study the performance of low rank tensor approximation methods in the large dimensional regime -- when the tensor dimensions $n_i\to \infty$, however considering symmetric tensor models where  $n_1=\ldots=n_d$, $\vx^{(1)} = \ldots = \vx^{(d)}$, and assuming that the noise is symmetric \cite{montanari2014statistical, perry2020statistical, lesieur2017statistical, handschy2019phase, jagannath2020statistical, de2021random}.

In particular, in the matrix case (i.e., $d=2$), the above spiked tensor model becomes a so-called spiked matrix model. It is well-known that in the large dimensional regime, there exists an order one critical value $\beta_c$ of the SNR below which it is information-theoretically impossible to detect/recover the spike, while above $\beta_c$, it is possible to detect the spike and recover the corresponding components in (at least) polynomial time using singular value decomposition (SVD). This phenomenon is sometimes known as the BBP (Baik, Ben Arous, and Péché) phase transition \cite{baik2005phase, benaych2011eigenvalues, capitaine2009largest, peche2006largest}. 

In the (symmetric) spiked tensor model for $d\geq 3$, there also exists an order one critical value\footnote{Depending on the tensor order $d$. We will sometimes omit the dependence on $d$ if there is no ambiguity.} $\beta_c(d)$ (in the high-dimensional asymptotic) below which it is information-theoretically impossible to detect/recover the spike, while above $\beta_c(d)$ recovery is theoretically possible with the maximum likelihood (ML) estimator. Computing the maximum likelihood in the matrix case corresponds to the computation of the largest singular vectors of the considered matrix which has a polynomial time complexity, while for $d\geq 3$, computing ML is NP-hard \cite{montanari2014statistical,biroli2020iron}. As such, a more practical phase transition for tensors is to characterize the algorithmic critical value $\beta_a(d,n)$ (which might depend on the tensor dimension $n$) above which the recovery of the spike is possible in polynomial time. Richard and Montanari \cite{montanari2014statistical} first introduced the symmetric spiked tensor model (of the form $ \tY =\mu \vx^{\otimes d} + \tW\in \sT_n^d$ with symmetric $\tW$) and also considered the related algorithmic aspects. In particular, they used heuristics to highlight that spike recovery is possible, with Approximate Message Passing (AMP) or the tensor power iteration method, in polynomial time\footnote{Using tensor power iteration or AMP with random initialization.} provided $\mu \gtrsim n^{\frac{d-1}{2}}$. This phase transition was later proven rigorously for AMP by \cite{lesieur2017statistical,jagannath2020statistical} and recently for tensor power iteration by \cite{huang2020power}.

Richard and Montanari \cite{montanari2014statistical} further introduced a method for tensor decomposition based on tensor unfolding, which consists in unfolding $\tY$ to an $n^{q}\times n^{d-q}$ matrix $\mat(\tY) = \mu \vx \vy^\top + \mZ$ for $q\in[d-1]$, to which a SVD is then performed. Setting $q = 1$, they predicted that their proposed method recovers successively the spike if $\mu \gtrsim n^{\frac{d-2}{4}}$. In a very recent paper by Ben Arous et al. \cite{arous2021long}, a study of spiked long rectangular random matrices\footnote{Number of rows $m$ are allowed to grow polynomially in the number of columns $n$, i.e., $\frac{m}{n} = n^{\alpha}$.} has been proposed under fairly general (bounded fourth-order moment) noise distribution assumptions. They particularly proved the existence of a critical SNR for which the extreme singular value and singular vectors exhibit a BBP-type phase transition. They applied their result for the asymmetric rank-one spiked model in \eqref{eq_asymmetric_spike_model} (with equal dimensions) using the tensor unfolding method, and found the exact threshold obtained by \cite{montanari2014statistical}, i.e., $\beta \gtrsim n^{\frac{d-2}{4}}$ for tensor unfolding to succeed in signal recovery.

For the asymmetric spiked tensor model in \eqref{eq_asymmetric_spike_model}, few results are available in the literature (to the best of our knowledge only \cite{arous2021long} considered this setting by applying the tensor unfolding method proposed by \cite{montanari2014statistical}). This is precisely the model we consider in the present work and our more general result is derived as follows. Given the asymmetric model from \eqref{eq_asymmetric_spike_model}, the maximum likelihood (ML) estimator of the best rank-one approximation of $\tT$ is given by
\begin{align} \label{eq_variational_definition}
    (\lambda_*, \vu_*^{(i)}) = \argmin_{\lambda\,\in\, \sR^+,\, (\vu^{(1)},\ldots,\vu^{(d)})\,\in\, \sS^{n_1 - 1} \times \cdots \times \sS^{n_d - 1}} \Vert \tT - \lambda \, \vu^{(1)} \otimes  \cdots \otimes \vu^{(d)} \Vert_{\text{F}}^2.
\end{align}
In the above equation, $\lambda_*$ and the $\vu_*^{(i)}$ can be respectively interpreted as the generalization to the tensor case of the concepts of dominant singular value and associated singular vectors \cite{lim2005singular}.
Following variational arguments therein, \eqref{eq_variational_definition} can be reformulated using contractions of $\tT$ as
\begin{align}\label{eq:objective}
    \max_{ \prod_{i=1}^d \Vert \vu^{(i)}\Vert = 1} \vert \tT(\vu^{(1)},\ldots,\vu^{(d)})\vert,
\end{align}
the Lagrangian of which writes as $\mathcal L \equiv \tT(\vu^{(1)},\ldots,\vu^{(d)}) -\lambda \left(\prod_{i=1}^d \Vert \vu^{(i)}\Vert - 1\right)$ with $\lambda >0$. Hence, the stationary points $(\lambda, \vu^{(1)},\ldots,\vu^{(d)})$, with the $\vu^{(i)}$'s being unitary vectors, must satisfy the Karush-Kuhn-Tucker conditions, for $i\in [d]$
\begin{align}\label{eq:order_d_singular}
\begin{cases}
 \tT(\vu^{(1)},\ldots,\vu^{(i-1)}, :, \vu^{(i+1)},\ldots ,\vu^{(d)}) = \lambda \vu^{(i)}, \\
\lambda = \tT(\vu^{(1)},\ldots,\vu^{(d)}).
\end{cases}
\end{align}

\begin{figure}[t!]
    \centering
    \includegraphics[width=10cm]{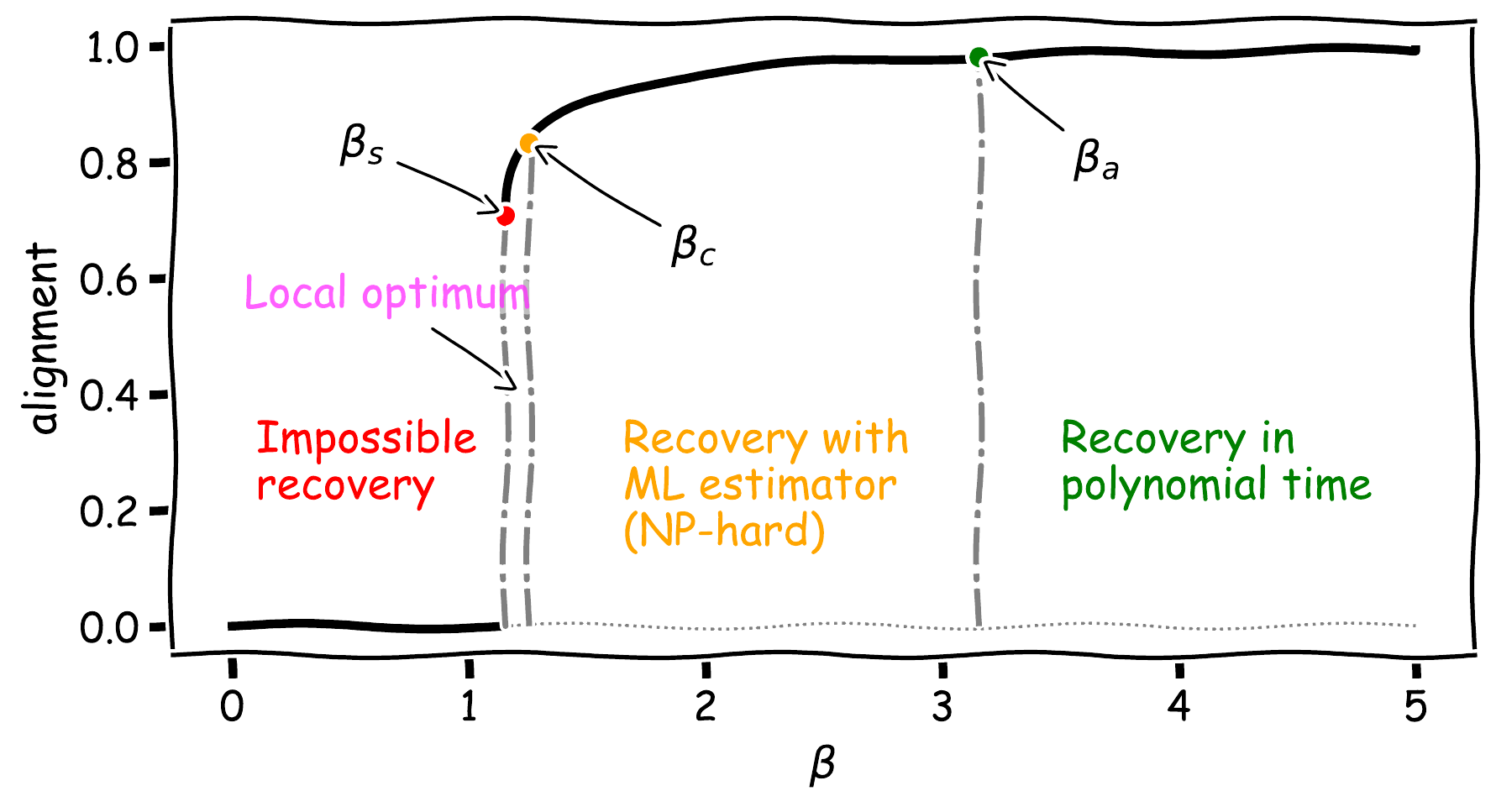}
    \caption{Illustration of the different thresholds for the SNR $\beta$. Our approach exhibits $\beta_s$ such that for $\beta>\beta_s$ there exists a local optimum that correlates with the spike, the threshold $\beta_c$ is unknown for asymmetric tensors and corresponds to the ML phase transition (above $\beta_c$ the global maximum correlates with the spike), while $\beta_a$ corresponds to the algorithmic phase transition (recovery of the spike in polynomial time).}
    \label{fig:illustration}
\end{figure}

An interesting question concerns the computation of the expected number of stationary points (local optima or saddle points) satisfying the identities in \eqref{eq:order_d_singular}. \cite{arous2019landscape} studied the landscape of a symmetric spiked tensor model and found that for $\beta < \beta_c$ the values of the objective function of all local maxima (including the global one) tend to concentrate on a small interval, while for $\beta > \beta_c$ the value achieved by the global maximum exits that interval and increases with $\beta$. In contrast, very recently Goulart et al. \cite{de2021random} have studied an order $3$ symmetric spiked random tensor $\tY$ using a RMT approach, where it was stated that there exists a threshold $0<\beta_s < \beta_c$ such that for $\beta\in [\beta_s, \beta_c]$ there exists a local optimum of the ML problem that correlates with the spike and such local optimum coincides with the global one for $\beta>\beta_c$. We conjecture that such observations extend to asymmetric spiked tensors, namely that there exists an order one critical value $\beta_c$ above which the ML problem in \eqref{eq:objective} admits a global maximum. As for \cite{de2021random}, our present findings do not allow to express such $\beta_c$ and its exact characterization is left for future investigation. However, for asymmetric spiked random tensors, we also exhibit a threshold $\beta_s$ such that for $\beta>\beta_s$ there exists a local optimum of the ML objective that correlates with the true spike. Figure \ref{fig:illustration} provides an illustration of the different thresholds of $\beta$ and see the last part of Subsection \ref{sec:equivalent_order_3} for a more detailed discussion.

\paragraph*{Main Contributions}
Starting form the conditions in \eqref{eq:order_d_singular}, we provide an exact expression of the asymptotic singular value and alignments $\langle \vx^{(i)}, \vu_*^{(i)} \rangle$, when the tensor dimensions $n_i\to \infty$ with $\frac{n_i}{\sum_{j=1}^d n_j}\to c_i\in(0, 1)$, where the tuple $(\lambda_*, \vu_*^{(1)}, \ldots, \vu_*^{(d)})$ is associated to a local optimum of the ML problem verifying some technical conditions (detailed in Assumption \ref{ass:lambda_outside_order_d}). 
We conjecture that when the SNR $\beta$ is large enough, there is a unique local optimum verifying Assumption \ref{ass:lambda_outside_order_d} and for which our results characterize the corresponding alignments.
We further conjecture that $(\lambda_*, \vu_*^{(1)}, \ldots, \vu_*^{(d)})$ coincides with the global maximum above some\footnote{Such a critical value has been characterized by \cite{jagannath2020statistical} for symmetric tensors. See \cite{de2021random} for a detailed discussion about this aspect in the case of symmetric tensors.} $\beta_c$ -- that needs to be characterized. 

Technically, we first show that the considered random tensor $\tT$ can be mapped to an equivalent \textit{symmetric} random matrix $\mT\in\sM_{N}$, constructed through contractions of $\tT$ with $d-2$ directions among $\vu_*^{(1)},\ldots,\vu_*^{(d)}$. 
Then, leveraging on random matrix theory, we first characterize the limiting spectral measure of $\mT$ and then provide estimates of the asymptotic alignments $\langle \vx^{(i)}, \vu_*^{(i)} \rangle$. We precisely show (see Theorem \ref{thm:asymptotics_general}) that under Assumption \ref{ass:lambda_outside_order_d}, for $d\geq 3$, there exists $\beta_s>0$ such that for $\beta > \beta_s$
\begin{align}\label{eq:asymptotic_d}
\begin{cases}
\lambda_* \asto \lambda^\infty(\beta),\\
\left\vert \langle \vx^{(i)}, \vu_*^{(i)} \rangle \right\vert \asto q_i(\lambda^\infty(\beta)),
\end{cases}
\end{align}
where $\lambda^\infty(\beta)$ satisfies $f(\lambda^\infty(\beta), \beta) = 0$ with $f(z, \beta) = z + g(z) - \beta \prod_{i=1}^d q_i(z) $ and
\begin{align*}
\color{cblue}
 q_i(z)  = \sqrt{1 - \frac{g_i^2(z)}{c_i} }, \quad g_i(z) = \frac{g(z) + z}{2} - \frac{\sqrt{4c_i + (g(z)+z)^2}}{2},
\end{align*}
$g(z)$ being the solution to the fixed point equation $g(z)=\sum_{i=1}^d g_i(z)$ (for $z$ large enough); see Section \ref{proof_existence_of_g} for the existence of $g(z)$. Besides, for $\beta \in [0, \beta_s]$ with\footnote{For arbitrary values of $c_i$'s, the upper bound of $\lambda^\infty$ can be computed numerically as the minimum non-negative real number $z$ for which Algorithm \ref{alg:stieltjes_transform} converges.} $c_i = \frac1d$ for all $i\in [d]$
\begin{align}
\lambda_* \asto \lambda^\infty \leq  2 \sqrt{\frac{d-1}{d}} ,\quad
\left\vert \langle \vx^{(i)}, \vu_*^{(i)} \rangle \right\vert \asto 0,
\end{align}
that is $\vu_*^{(i)}$ (asymptotically) ceases to correlate with $\vx^{(i)}$.

{\color{cblue}
\begin{remark}\label{remark_expression}
    Note that $q_i(z)$ can be equivalently expressed as $q_i(z) =  \left( \frac{\alpha_i(z)^{d-3}}{ \prod_{j\neq i} \alpha_j(z) } \right)^{\frac{1}{2d-4}}$ with $\alpha_i(z) = \frac{\beta}{z + g(z) - g_i(z)}$. Note that such expression is defined for $c_i\in[0, 1]$ with $d\geq 3$. See details in Section \ref{proof:asymptotics_general}.
\end{remark}}

We highlight that in our formulation the threshold $\beta_s$ corresponds to the minimal value of the SNR $\beta$ above which the derived asymptotic formulas are algebraically defined, which may differ from the true phase transition $\beta_c$ of the ML problem: In the case of symmetric tensors, the results form  \cite{de2021random} seem to indicate that $\beta_s$ is slightly below $\beta_c$ obtained by \cite{jagannath2020statistical} where the ML problem was studied.

\begin{figure}[t!]
    \centering
    \includegraphics[width=14cm]{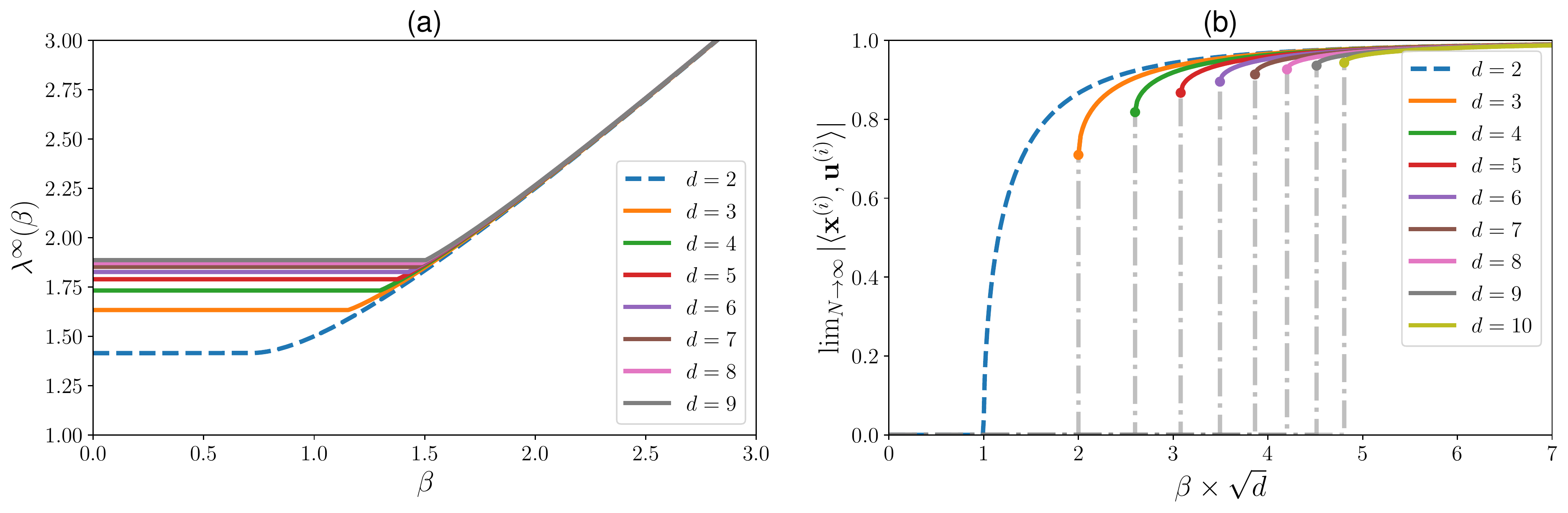}
    \caption{Asymptotic singular value and alignments as per \eqref{eq:asymptotic_d} when all the dimensions $n_i$ are equal ($c_i=\frac1d$ for all $i\in [d]$) for different values of the tensor order $d$. The matrix case corresponds to setting $d=3$ and $c_3 = 0$ (see Section \ref{sec:tensors2matrices} for more details).}
    \label{fig:align_order}
\end{figure}

Figure \ref{fig:align_order} notably depicts the asymptotic alignments as in \eqref{eq:asymptotic_d} when all the tensor dimensions $n_i$ are equal. Since our result characterizes the alignments for $d\geq 3$, it is not possible to recover the matrix case by simply setting $d=2$ since the equations are not defined for $d=2$ (see Remark \ref{remark_expression}). However, the matrix case is recovered by considering an order $3$ tensor ($d=3$) and then taking the limit $c_3\to 0$ (equivalent to the degenerate case $n_3=1$ which results in a spiked random matrix model), see Section \ref{sec:tensors2matrices} for more details. From Figure \ref{fig:align_order}{-(b)}, unlike the matrix case, i.e., $d=2$, the predicted asymptotic alignments are not continuous for orders $d\geq 3$, this phenomenon has already been observed in the case of symmetric random tensors \cite{jagannath2020statistical}. In particular, the predicted theoretical threshold $\beta_s$ in the matrix case $d=2$ coincides with the classical so-called BBP (Baik, Ben Arous, and Péché) phase transition $\beta_c(2)$ \cite{baik2005phase, benaych2011eigenvalues, capitaine2009largest, peche2006largest}. Moreover, our result for the matrix case characterizes the asymptotic alignments for the long rectangular matrices studied by \cite{arous2021long} and we also recover the threshold $\beta \gtrsim n^{\frac{d-2}{4}}$ for the case of tensor unfolding method (See remark \ref{remark:unfolding}). From a methodological view point, our results are derived based solely on Stein's lemma without the use of complex contour integrals as classically performed in RMT. In essence, we follow the approach in \cite{de2021random} and further introduce a new object (the mapping $\Phi_d$ defined subsequently in \eqref{eq:phi_d}) that simplifies drastically the use of RMT tools.\\

The remainder of the paper is organized as follows. Section \ref{sec:rmt_tools} recalls some fundamental random matrix theory tools. In Section \ref{sec:order_3}, we study asymmetric spiked tensors of order $3$ where we provide the main steps of our approach. Section \ref{sec:tensors2matrices} characterizes the behavior of spiked random matrices given our result on spiked $3$-order tensor models from Section \ref{sec:order_3}. The generalization of our results to arbitrary $d$-order tensors is then presented in Section \ref{sec:order_d}. Section \ref{sec:rank_k} discusses the application of our findings to arbitrary rank-$k$ tensors with mutually orthogonal components. Further discussions are presented in Section \ref{sec:discussion}. In Appendix \ref{sec:simus}, we provide some simulations to support our findings and discuss some algorithmic aspects. Finally, Appendix \ref{sec:proofs} provides the proofs of the main developed results.

\section{Random matrix theory tools}\label{sec:rmt_tools}
Before digging into our main findings, we briefly recall some random matrix theory results that are at the heart of our analysis. Specifically, we recall the \textit{resolvent} formalism which allows to assess the spectral behavior of large \textit{symmetric} random matrices. In particular, given a symmetric matrix $\mS\in \sM_n$ and denoting $\lambda_1(\mS)\leq \ldots\leq \lambda_n(\mS)$ its $n$ eigenvalues with corresponding eigenvectors $\vu_i(\mS)$ for $i\in [n]$, its spectral decomposition writes as
\begin{align*}
    \mS = \sum_{i=1}^n \lambda_i(\mS) \vu_i(\mS) \vu_i(\mS)^\top.
\end{align*}
In the sequel we omit the dependence on $\mS$ by simply writing $\lambda_i$ and $\vu_i$ if there is no ambiguity.
\subsection{The resolvent matrix}
We start by defining the \textit{resolvent} of a symmetric matrix and present its main properties. 
\begin{definition}\label{def:resolvent}
Given a symmetric matrix $\mS\in \sM_n$, the resolvent of $\mS$ is defined as
\begin{align*}
    \mR_\mS(z) \equiv \left( \mS - z \mI_n \right)^{-1},\quad z \in \sC \setminus \supp(\mS),
\end{align*}
where $\supp(\mS) = \{\lambda_1, \ldots, \lambda_n \}$ is the spectrum of $\mS$.
\end{definition}
The resolvent is a fundamental object since it retrieves the spectral characteristics (spectrum and eigenvectors) of $\mS$. It particularly verifies the following property which we will use extensively to derive our main results,
\begin{align}\label{eq:resolvent_identity}
    \mR_\mS(z) = -\frac{1}{z} \mI_n + \frac{1}{z} \mS \mR_\mS(z) = -\frac{1}{z} \mI_n + \frac{1}{z} \mR_\mS(z) \mS.
\end{align}
The above identity, coupled with Stein's Lemma (Lemma~\ref{lemma:stein} in Section~\ref{section_Gaussian_calculations} below), is a fundamental tool used to derive fixed point equations that allow the evaluation of functionals of interests involving $\mR_\mS(z)$.
Another interesting property of the resolvent concerns its spectral norm, which we denote by $\Vert \mR_\mS(z) \Vert$. Indeed, if the spectrum of $\mS$ has a bounded support, then the spectral norm of $\mR_\mS(z)$ is bounded. This is a consequence of the inequality
\begin{align}\label{eq:resovent_bounded}
    \Vert \mR_\mS(z) \Vert \leq \frac{1}{\dist(z, \supp(\mS))}
\end{align}
where $\dist(z,\cdot)$ denotes the distance of $z$ to a set.
The resolvent encodes rich information about the behavior of the eigenvalues of $\mS$ through the so-called Stieltjes transform which we describe subsequently.
\subsection{The Stieltjes transform}\label{sec:stieltjes_transform}
Random matrix theory, originally, aims at describing the limiting spectral measure of random matrices when their dimensions grow large. Typically, under certain technical conditions on $\mS$, the \textit{empirical spectral measure} of $\mS$ defined as
\begin{align}
    \nu_\mS \equiv \frac1n \sum_{i=1}^n \delta_{\lambda_i(\mS)},
\end{align}
where $\delta_x$ is a Dirac mass on $x$, converges to a deterministic probability measure $\nu$. To characterize such asymptotic measure, the Stieltjes transform (defined below) approach is a widely used tool.
\begin{definition}[Stieltjes transform]
Given a probability measure $\nu$, the Stieltjes transform of $\nu$ is defined by
\begin{align*}
    g_\nu(z)\equiv \int \frac{d\nu(\lambda)}{\lambda - z}, \quad z\in \sC \setminus \supp(\nu).
\end{align*}
\end{definition}
Particularly, the Stieltjes transform of $\nu_\mS$ is closely related to its associated resolvent $\mR_\mS(z)$ through the algebraic identity
\begin{align*}
    g_{\nu_\mS}(z) = \frac1n \sum_{i=1}^n \int \frac{\delta_{\lambda_i(\mS)}(d\lambda)}{\lambda - z} = \frac1n \sum_{i=1}^n \frac{1}{\lambda_i(\mS) - z} = \frac1n \tr \mR_\mS(z)
\end{align*}
The Stieltjes transform $g_\nu$ has several interesting analytical properties, among which, we have
\begin{enumerate}
    \item $g_\nu$ is complex analytic on its definition domain $\sC\setminus \supp(\nu)$ and $\Im [g_\nu(z)] >0$ if $\Im [z] >0$;
    \item $g_\nu(z)$ is bounded for $z\in \sC \setminus \supp(\nu)$ if $\supp(\nu)$ is bounded, since $\vert g_\nu(z)\vert \leq \dist(z, \supp(\nu))^{-1} $;
    \item Since $g_\nu'(z) = \int (\lambda -z)^{-2}d\nu(\lambda)>0$, $g_\nu$ is monotonously increasing for $z\in \mathbb{R}$. 
\end{enumerate}
The Stieltjes transform $g_\nu$ admits an inverse formula which provides access to the evaluation of the underlying probability measure $\nu$, as per the following theorem.
\begin{theorem}[Inverse formula of the Stieltjes transform]\label{thm:inverse} Let $a,b$ be some continuity points of the probability measure $\nu$, then the segment $[a,b]$ is measurable with $\nu$, precisely
\begin{align*}
    \nu([a,b]) = \frac1\pi \lim_{\epsilon\to 0} \int_a^b \Im[g_\nu(x + i\epsilon)]dx.
\end{align*}
Moreover, if $\nu$ admits a density function $f$ at some point $x$, i.e., $\nu(x)$ is differentiable in a neighborhood of $x$ with $\lim_{\epsilon\to 0}\epsilon^{-1}\nu([x-\epsilon/2, x+\epsilon/2])=f(x)$, then we have the inverse formula
\begin{align*}
    f(x) = \frac1\pi \lim_{\epsilon\to 0} \Im[g_\nu(x + i\epsilon)].
\end{align*}
\end{theorem}
When it comes to large random matrices, the Stieltjes transform admits a continuity property, in the sense that if a sequence of random probability measures converges to a deterministic measure then the corresponding Stieltjes transforms converge almost surely to a deterministic one, and vice-versa. The following theorem from \cite{tao2012topics} states precisely this continuity property.
\begin{theorem}[Stieltjes transform continuity]\label{thm:continuity}
A sequence of random probability measures $\nu_n$, supported on $\sR$ with corresponding Stieltjes transforms $g_{\nu_n}(z)$, converges almost surely weakly to a deterministic measure $\nu$, with corresponding Stieltjes transform $g_\nu$ if and only if $g_{\nu_n}(z)\to g_{\nu}(z)$ almost surely, for all $z\in \sR + i\sR_+$.
\end{theorem}
In particular, Theorems \ref{thm:inverse} and \ref{thm:continuity} combined together allow the description of the spectrum of large symmetric random matrices.

\subsection{Gaussian calculations}
\label{section_Gaussian_calculations}
The following lemma by Stein (also called Stein's identity or Gaussian integration by parts) allows to replace the expectation of the product of a Gaussian variable with a differentiable function $f$ by the variance of that variable times the expectation of $f'$.
\begin{lemma}[Stein's Lemma \cite{stein1981estimation}]\label{lemma:stein}
Let $X\sim \mathcal{N}(0, \sigma^2)$ and $f$ a continuously differentiable function having at most polynomial growth, then
\begin{align*}
    \EE\left[ X f(X) \right] = \sigma^2\EE \left[ f'(X) \right],
\end{align*}
when the above expectations exist.
\end{lemma}

We will further need the Poincaré inequality which allows to control the variance of a functional of $i.i.d.$ Gaussian random variables. 
\begin{lemma}[Poincaré inequality]\label{lem:poincaré} Let $F:\sR^n \to \sR$ a continuously differentiable function having at most polynomial growth and $X_1,\ldots, X_n$ a collection of i.i.d.\@ standard Gaussian variables. Then,
\begin{align*}
\Var F(X_1, \ldots, X_n) \leq \sum_{i=1}^n \EE \left\vert \frac{\partial F}{\partial X_i} \right\vert^2.
\end{align*}
\end{lemma}

\section{The \textit{asymmetric} rank-one spiked tensor model of order $3$}\label{sec:order_3}
To best illustrate our approach, we start by considering the asymmetric rank-one tensor model of order $3$ of the form
\begin{align}\label{eq:spiked_tensor_model}
    \tT = \beta\, \vx\otimes \vy \otimes \vz + \frac{1}{\sqrt{N}} \tX,
\end{align}
where $\tT$ and $\tX$ both have dimensions $m\times n \times p$ and $N = m+n+p$. We assume that $\vx,\vy$ and $\vz$ are on the unit spheres $\sS^{m-1}$, $\sS^{n-1}$ and $\sS^{p-1}$ respectively, $\tX$ is a Gaussian noise tensor with i.i.d.\@ entries $X_{ijk}\sim \mathcal N(0,1)$ independent from $\vx,\vy,\vz$. 
\subsection{Tensor singular value and vectors}
According to \eqref{eq:order_d_singular}, the $\ell_2$-singular value and vectors, corresponding to the best rank-one approximation $\lambda_* \vu_* \otimes \vv_* \otimes \vw_*$ of $\tT$, satisfy the identities 
\begin{align}\label{eq:singular_vectors}
    \tT(\vv)\vw = \lambda \vu \quad \tT(\vu) \vw = \lambda \vv\quad \tT(\vv)^\top \vu = \lambda \vw \quad \text{with} \quad (\vu,\vv,\vw)\,\in\, \sS^{m-1}\times \sS^{n-1}\times \sS^{p-1},
\end{align}
where we denoted $\tT(\vu) = \tT(\vu,:,:)\in \sM_{n,p}, \tT(\vv) = \tT(:,\vv,:)\in \sM_{m,p}$ and $\tT(\vw) = \tT(:,:,\vw)\in \sM_{m,n}$.
Furthermore, the singular value $\lambda$ can be characterized through the contraction of $\tT$ along all its singular vectors $\vu,\vv,\vw$, i.e.
\begin{align}\label{eq:singular_value}
\lambda = \tT(\vu,\vv,\vw) = \sum_{ijk} u_i v_j w_k T_{ijk}.
\end{align}

\subsection{Associated random matrix model}\label{sec:equivalent_order_3}
We follow the approach developed in \cite{de2021random} which consists in studying random matrices that are obtained through contractions of a random tensor model, and extend it to the \textit{asymmetric} spiked tensor model of \eqref{eq:spiked_tensor_model}. Indeed, it has been shown in \cite{de2021random} that the study of a rank-one \textit{symmetric} spiked tensor model $\tY\in \sT_{q}^3$ boils down to the analysis of the \textit{symmetric} random matrix\footnote{The contraction of the tensor $\tY$ with its eigenvector $\va$.} $\tY(\va)\in \sM_q$ where $\va \in \sS^{q-1}$ stands for the eigenvector of $\tY$ corresponding to the best symmetric rank-one approximation of $\tY$, i.e.,
\begin{align}
    (\mu,\va) = \argmin_{\mu\in \sR,\, \va\in \sS^{q-1}} \Vert \tY - \mu \va^{\otimes 3}\Vert_{\text{F}}^2 \quad \Leftrightarrow \quad  \tY(\va) \va = \mu \va.
\end{align}
In the \textit{asymmetric} case of \eqref{eq:spiked_tensor_model}, the choice of a ``relevant'' random matrix to study is not trivial since the corresponding contractions $\tT(\vu), \tT(\vv)$ and $\tT(\vw)$ yield asymmetric random matrices which present more technical difficulties from the random matrix theory perspective, therefore the extension of the approach in \cite{de2021random} to asymmetric tensors is not straightforward. We will see in the following that such a choice of the relevant matrix to study an asymmetric $\tT$ is naturally obtained through the use of the Pastur's Stein approach \cite{marvcenko1967distribution}.

{\color{cblue}
As described in \cite{de2021random} and since the singular vectors $\vu,\vv$ and $\vw$ depend statistically on $\tX$, the first technical challenge consists in expressing the derivatives of the singular vectors $\vu,\vv$ and $\vw$ w.r.t.\@ the entries of the Gaussian noise tensor $\tX$. Indeed, one can show that there exists a differentiable mapping $\mathcal{F}: \sT_{m,n,p} \to \sR^{m+n+p+1} $ that maps $\tX$ to $\mathcal{F}(\tX) = (\lambda(\tX), \vu(\tX), \vv(\tX), \vw(\tX))$ singular-value and vectors of $\tT$, since the components of $\vu(\tX), \vv(\tX)$ and  $\vw(\tX)$ are bounded and $\lambda(\tX)$ has polynomial growth. Indeed, we have the following Lemma which is analog to \cite[Lemma 8]{de2021random} and which justifies the application of Stein's Lemma subsequently.

\begin{lemma}\label{lemma_differentiable}
    There exists an almost everywhere continuously differentiable function $\mathcal{F}: \sT_{m,n,p} \to \sR^{m+n+p+1} $ such that $\mathcal{F}(\tX) = (\lambda(\tX), \vu(\tX), \vv(\tX), \vw(\tX))$ is singular-value and vectors of $\tT$ (for almost every $\tX$).
\end{lemma}
\begin{proof}
    The proof relies on the same arguments as \cite[Lemma 8]{de2021random}. 
\end{proof}
}
Calculus (see details in \ref{proof:derivs}) show that deriving the identities in \eqref{eq:singular_vectors} w.r.t.\@ an entry $X_{ijk}$ of $\tX$ with $(i,j,k)\in [m]\times[n]\times[p]$ results in
\begin{align}\label{eq:deriv}
\begin{bmatrix}
\frac{\partial \vu}{\partial X_{ijk}} \\ 
\frac{\partial \vv}{\partial X_{ijk}} \\ 
\frac{\partial \vw}{\partial X_{ijk}}
\end{bmatrix} 
= - \frac{1}{\sqrt N} \left(  \begin{bmatrix}
\vzero_{m\times m} & \tT(\vw) & \tT(\vv) \\ 
\tT(\vw)^\intercal & \vzero_{n\times n} & \tT(\vu) \\ 
\tT(\vv)^\intercal & \tT(\vu)^\intercal & \vzero_{p\times p}
\end{bmatrix} - \lambda \mI_{N}  \right)^{-1} 
\begin{bmatrix}
 v_j w_k (\ve_i^m - u_i \vu) \\ 
 u_i w_k (\ve_j^n - v_j \vv) \\
 u_i v_j (\ve_k^p - w_k \vw)
\end{bmatrix} \in \sR^{N},
\end{align}
where we recall that $\ve_i^d$ denotes the canonical vector in $\sR^d$ with $[\ve_i^d]_j = \delta_{ij}$. We further have the identity
\begin{align}\label{eq:deriv_lambda}
\frac{\partial \lambda}{\partial X_{ijk}} = \frac{1}{\sqrt N} u_i v_j w_k.
\end{align}
Denoting by $\mM$ the \textit{symmetric} block-wise random matrix which appears in the matrix inverse in \eqref{eq:deriv}, the derivatives of $\vu,\vv$ and $\vw$ are therefore expressed in terms of the \textit{resolvent} $\mR_\mM(z)=\left(\mM - z\mI_N \right)^{-1}$ evaluated on $\lambda$, and we will see subsequently that the assessment of the spectral properties of $\tT$ boils down to the estimation of the normalized trace $\frac1N\tr\mR_\mM(z)$. As such, the matrix $\mM$ provides the associated random matrix model that encodes the spectral properties of $\tT$.
We will henceforth focus our analysis on this random matrix, in order to assess the spectral behavior of $\tT$. More generally, we will be interested in studying random matrices from the $3$-order \textit{block-wise tensor contraction ensemble} $\mathcal{B}_3(\tX)$ for $\tX \sim\sT_{m,n,p}(\mathcal{N}(0,1))$ defined as
\begin{align}
    \mathcal{B}_3(\tX) \equiv \left\lbrace \Phi_3(\tX,\va,\vb,\vc) \,\, \big\vert \,\,  (\va,\vb,\vc) \in \sS^{m-1}\times \sS^{n-1}\times \sS^{p-1}  \right\rbrace,
\end{align}
where $\Phi_3$ is the mapping 
\begin{equation}
    \begin{split}
            \Phi_3: \sT_{m,n,p}\times \sS^{m-1}\times \sS^{n-1}\times \sS^{p-1} & \longrightarrow \sM_{m+n+p} \\
    (\tX, \va, \vb, \vc) & \longmapsto \begin{bmatrix}
\vzero_{m\times m} & \tX(\vc) & \tX(\vb) \\ 
\tX(\vc)^\intercal & \vzero_{n\times n} & \tX(\va) \\ 
\tX(\vb)^\intercal & \tX(\va)^\intercal & \vzero_{p\times p}
\end{bmatrix}.
    \end{split}
\end{equation}

We associate tensor $\tT$ to the random matrix
\begin{align} \label{eq_phi3}
    \Phi_3(\tT, \vu, \vv, \vw ) =  \beta\, \mV \mB \mV^\top + \frac{1}{\sqrt N} \Phi_3(\tX, \vu, \vv, \vw )  \in \sM_N,
\end{align}
where
\begin{align*}
      \mB \equiv \begin{bmatrix}
 0 & \langle \vz, \vw\rangle & \langle \vy, \vv\rangle \\
 \langle \vz, \vw\rangle & 0 & \langle \vx, \vu\rangle \\
 \langle \vy, \vv\rangle & \langle \vx, \vu\rangle & 0
\end{bmatrix}\in \sM_3, \quad \mV \equiv \begin{bmatrix}
 \vx & \vzero_m & \vzero_m \\
 \vzero_n & \vy & \vzero_n \\
 \vzero_p & \vzero_p & \vz
\end{bmatrix}\in \sM_{N,3}.
\end{align*}
and $\vx,\vy$ and $\vz$ are on the unit spheres $\sS^{m-1}$, $\sS^{n-1}$ and $\sS^{p-1}$ respectively. Note that $\Phi_3(\tT, \vu, \vv, \vw )$ is \textit{symmetric} and behaves as a so-called spiked random matrix model where the signal part $\beta\, \mV \mB \mV^\top$ correlates with the true signals $\vx,\vy$ and $\vz$ for a sufficiently large $\beta$. However, the noise part (i.e., the term $\frac{1}{\sqrt N}\Phi_3(\tX, \vu, \vv, \vw )$ in \eqref{eq_phi3}) has a non-trivial structure due to the statistical dependencies between the singular vectors $\vu,\vv,\vw$ and the tensor noise $\tX$. Despite this statistical dependency, we will show in the next subsection that the \textit{spectral measure} of $\Phi_3(\tT, \vu, \vv, \vw )$ (and that of $\frac{1}{\sqrt N}\Phi_3(\tX, \vu, \vv, \vw )$) converges to a deterministic measure (see Theorem \ref{thm:semi-circle_dependent}) that coincides with the limiting spectral measure of $\Phi_3(\tT, \va, \vb, \vc )$ where $\va,\vb,\vc$ are unit vectors and independent of $\tX$.
Furthermore, the matrix $\Phi_3(\tT, \vu, \vv, \vw )$ admits $2\lambda$ as an eigenvalue (regardless of the value of $\beta$), which is a simple consequence of the identities in \eqref{eq:singular_vectors}, since 
\begin{align}\label{eq:2lambda}
    \Phi_3(\tT, \vu, \vv, \vw ) \begin{bmatrix}
     \vu \\
     \vv \\
     \vw
    \end{bmatrix}= \begin{bmatrix}
     \tT(\vw) \vv + \tT(\vv)\vw \\
     \tT(\vw)^\top \vu + \tT(\vu) \vw \\
     \tT(\vv)^\top \vu + \tT(\vu)^\top \vv
    \end{bmatrix} = 2\lambda \begin{bmatrix}
     \vu \\
     \vv \\
     \vw
    \end{bmatrix}.
\end{align}

Note however that the expression in \eqref{eq:deriv} exists only if $\lambda$ is not an eigenvalue of $\Phi_3(\tT, \vu, \vv, \vw)$. As discussed in \cite{de2021random} in the symmetric case, such condition is related to the locality of the maximum ML estimator. Indeed, we first have the following remark that concerns the Hessian of the underlying objective function.

\begin{remark}\label{rem:hessian}
Recall the Lagrangian of the ML problem as $\mathcal{L}(\vu, \vv, \vw) = \tT(\vu, \vv, \vw) - \lambda \left( \Vert \vu \Vert \Vert \vv \Vert \Vert \vw \Vert - 1\right)$ and denote $\vh \equiv \frac{[\vu^\top, \vv^\top, \vw^\top]^\top}{\sqrt{3}} $. If $\lambda$ is a singular value of $\tT$ with associated singular vectors $\vu, \vv, \vw$, then $(\lambda, \vh)$ is an eigenpair of the Hessian $\nabla^2\mathcal{L}\equiv \frac{\partial^2 \mathcal{L}}{\partial \vh^2}$ evaluated at $\vh$. Indeed, $\nabla^2\mathcal{L}(\vh) = \Phi_3(\tT, \vu, \vv, \vw) - \lambda \mI_N$, thus $\nabla^2\mathcal{L} (\vh) \vh = \Phi_3(\tT, \vu, \vv, \vw) \vh - \lambda \vh = \lambda \vh$ since $\Phi_3(\tT, \vu, \vv, \vw)\vh = 2\lambda \vh$ (from \eqref{eq:2lambda}).
\end{remark}

From \cite[Theorem 12.5]{nocedal2006numerical}, a necessary condition for $\vu, \vv, \vw$ to be a local maximum of the ML problem is that 
\begin{align*}
\langle \nabla^2\mathcal{L}(\vh) \vk, \vk \rangle \leq 0, \quad \forall \vk \in \vh^{\perp} \equiv \left\lbrace \vk \in \sR^N \, \mid \, \langle \vh, \vk \rangle = 0 \right\rbrace.
\end{align*}
which yields (by remark \ref{rem:hessian}) the condition 
\begin{align}
\label{eq_locality}
\max_{\vk \in \sS^{N-1} \cap \vh^{\perp} } \langle \Phi_3(\tT, \vu, \vv, \vw)\vk, \vk \rangle \leq \lambda.
\end{align} 
As such, for $\lambda >0$, $2\lambda$ must be the largest eigenvalue of $\Phi_3(\tT, \vu, \vv, \vw)$ as per \eqref{eq:2lambda}, while its second largest eigenvalue cannot exceed $\lambda$ as shown by \eqref{eq_locality}. Thus our analysis applies only to some local optimum of the ML problem for which the corresponding singular value $\lambda$ lies outside the \textit{bulk} of $\Phi_3(\tT, \vu, \vv, \vw)$, namely we suppose\footnote{We will find that such assumption is satisfied provided $\beta > \beta_s$ for some $\beta_s>0$ that we will determine} that there exists a tuple $(\lambda, \vu, \vv, \vw)$ verifying the identities in \eqref{eq:singular_vectors} such that $\lambda$ is not an eigenvalue of $\Phi_3(\tT, \vu, \vv, \vw)$. Moreover, this allows the existence of the matrix inverse in \eqref{eq:deriv} and it is not restrictive in the sense that it is satisfied for a sufficiently large value of the SNR $\beta$.\\
In the sequel, for any $(\lambda, \vu, \vv, \vw)$ verifying the identities in \eqref{eq:singular_vectors}, we find that the largest eigenvalue $2\lambda$ is always isolated from the \textit{bulk} of $\Phi_3(\tT, \vu, \vv, \vw)$ independently of the SNR $\beta$. In addition, there exists $\beta_s>0$ such that for $\beta\leq \beta_s$, the observed spike is not informative (in the sense that the corresponding alignments will be null) and is visible (an isolated eigenvalue appears in the spectrum of $\Phi_3(\tT, \vu, \vv, \vw)$, see Figure \ref{fig:spectrum_dependent}) because of the statistical dependencies between $\vu, \vv, \vw$ and $\tX$. The same phenomenon has been observed in the case of symmetric random tensors by \cite{de2021random}. 
 
\subsection{Limiting spectral measure of block-wise $3$-order tensor contractions}
We start by presenting our first result which characterizes the limiting spectral measure of the ensemble $\mathcal{B}_3(\tX)$ with $\tX \sim \sT_{m,n,p}(\mathcal{N}(0,1))$. We characterize this measure in the limit when the tensor dimensions grow large as in the following assumption.

\begin{assumption}[Growth rate]\label{ass:growth}
As $ m,n,p\to \infty$, the dimension ratios $\frac{m}{m+n+p}\to c_1 \in (0, 1)$, $\frac{n}{m+n+p}\to c_2 \in (0, 1)$ and $\frac{p}{m+n+p}\to c_3=1-(c_1+c_2) \in (0, 1)$.
\end{assumption}

We have the following theorem which characterizes the spectrum of $ \frac{1}{\sqrt N} \Phi_3(\tX, \va, \vb, \vc)$ with arbitrary deterministic unit vectors $\va, \vb, \vc$.

\begin{theorem}\label{thm:semi-circle_independent}
Let $\tX \sim \sT_{m,n,p}(\mathcal{N}(0,1))$ be a sequence of random asymmetric Gaussian tensors and $(\va,\vb,\vc) \in \sS^{m-1}\times \sS^{n-1}\times \sS^{p-1}$ a sequence of deterministic vectors of increasing dimensions, following Assumption~\ref{ass:growth}. Then the empirical spectral measure of $\frac{1}{\sqrt N}  \Phi_3(\tX, \va, \vb, \vc)$ converges weakly almost surely to a deterministic measure $\nu$ whose Stieltjes transform $g(z)$ is defined as the solution to the equation $g(z) = \sum_{i=1}^3 g_i(z)$ such that $\Im[g(z)]>0$ for $\Im[z]>0$ where, for $i\in [3]$ $g_i(z)$ satisfies $g_i^2(z) - (g(z) + z) g_i(z) - c_i = 0$ for $z\in \sC\setminus \mathcal{S}(\nu)$.
\end{theorem}

\begin{proof}
See Appendix \ref{proof:semi-circle_independent}.
\end{proof}

\begin{remark}
\label{remark_polynome} The equation $g(z) = \sum_{i=1}^3 g_i(z)$ yields a polynomial equation of degree $5$ in $g$ which can be solved numerically or via an iteration procedure which converges to a fixed point for $z$ large enough. 
\end{remark}

\begin{figure}[t!]
    \centering
    \includegraphics[width=14cm]{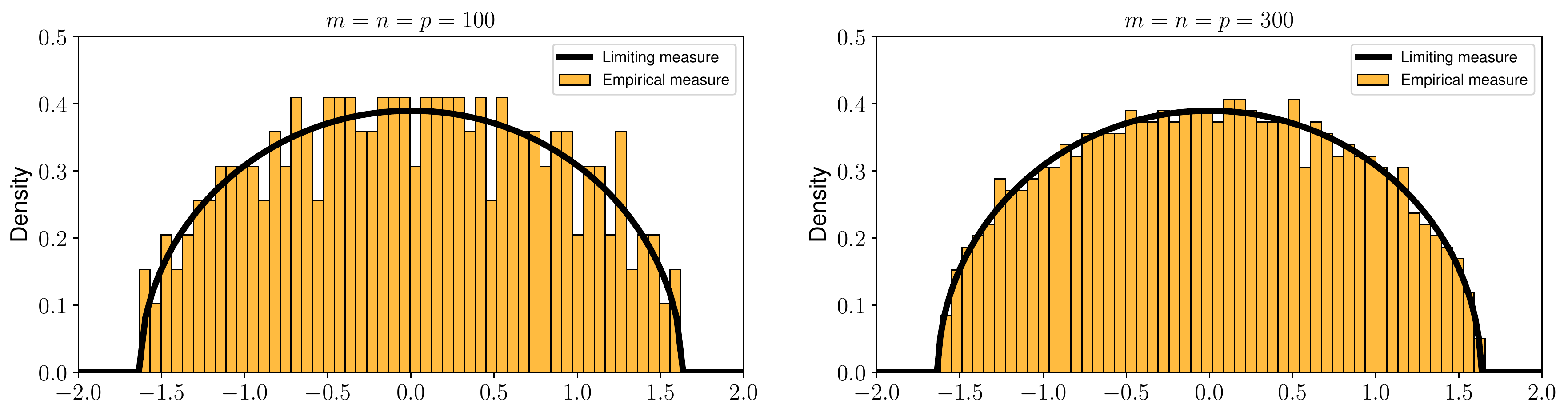}
    \caption{Spectrum of $ \frac{1}{\sqrt N} \Phi_3(\tX, \va, \vb, \vc)$ with $\tX \sim \sT_{m,n,p}(\mathcal{N}(0, 1))$ and $\va,\vb,\vc$ independently sampled from the unit spheres $\sS^{m-1}, \sS^{n-1}$, $\sS^{p-1}$ respectively. In black the semi-circle law as per Theorem \ref{thm:semi-circle_independent}.}
    \label{fig:spectrum_independent}
\end{figure}

In particular, in the cubic case when $c_1=c_2=c_3=\frac13$, the empirical spectral measure of $\frac{1}{\sqrt N}  \Phi_3(\tX, \va, \vb, \vc)$ converges to a semi-circle law as precisely stated by the following Corollary of Theorem \ref{thm:semi-circle_independent}.

\begin{corollary}\label{cor:semi-circle-cubic} Given the setting of Theorem \ref{thm:semi-circle_independent} with $c_1=c_2=c_3=\frac13$, the empirical spectral measure of $\frac{1}{\sqrt N}  \Phi_3(\tX, \va, \vb, \vc)$ converges weakly almost surely to the semi-circle distribution supported on $\mathcal{S}(\nu)\equiv\left[-2\sqrt{\frac23}, 2\sqrt{\frac23}\right]$, whose density and Stieltjes transform write respectively as
\begin{align*}
\nu(dx)=\frac{3}{4\pi}\sqrt{\left( \frac{8}{3} - x^2 \right)^+}dx, \quad g(z)\equiv \frac{-3z + 3\sqrt{z^2 - \frac83}}{4},\quad \text{where} \quad z\in \sC\setminus \mathcal{S}(\nu).
\end{align*}
\end{corollary}

\begin{proof}
See Appendix \ref{proof:semi-circle-cubic}.
\end{proof}

Figure \ref{fig:spectrum_independent} depicts the spectrum of $ \frac{1}{\sqrt N} \Phi_3(\tX, \va, \vb, \vc)$ with $\tX \sim \sT_{m,n,p}(\mathcal{N}(0, 1))$ and independent unit vectors $\va,\vb, \vc$, it particularly illustrates the convergence in law of this spectrum when the dimensions $m,n,p$ grow large. 

\begin{remark} More generally, the spectral measure of $\frac{1}{\sqrt N} \Phi_3(\tX, \va, \vb, \vc)$ converges to a deterministic measure with connected support if $\frac{1}{\sqrt N} \Phi_3(\tX, \va, \vb, \vc)$ is almost surely full rank, i.e., if $\max(m,n)-\min(m,n)\leq p \leq m+n$ since the rank of $\frac{1}{\sqrt N} \Phi_3(\tX, \va, \vb, \vc)$ is $\min(m,n+p) + \min(n,m+p) + \min(p,m+n)$. In contrast if it is not full rank, its spectral measure converges to a deterministic measure with unconnected support (see the case of matrices in Corollary \ref{cor:spectrum_matrices} subsequently).
\end{remark}

The analysis of the tensor $\tT$ relies on describing the spectrum of $ \Phi_3(\tT, \vu, \vv, \vw)$ where the singular vectors $\vu,\vv, \vw$ depend statistically on $\tX$ (the noise part of $\tT$). Despite these dependencies, it turns out that the spectrum of $\Phi_3(\tT, \vu, \vv, \vw)$ converges in law to the same deterministic measure described by Theorem \ref{thm:semi-circle_independent}. Besides, we need a further technical assumption on the singular value and vectors of $\tT$.
{\color{cblue}
\begin{assumption}\label{ass:lambda_outside}
    We assume that there exists a sequence of critical points $(\lambda_*, \vu_*,\vv_*, \vw_*)$ satisfying \eqref{eq:singular_vectors} such that $\lambda_* \asto \lambda^\infty(\beta)$, $\vert\langle \vu_*, \vx \rangle\vert \asto a_x^\infty(\beta )$, $\vert\langle \vv_*, \vy \rangle\vert \asto a_y^\infty(\beta )$, $\vert\langle \vw_*, \vz \rangle\vert \asto a_z^\infty(\beta )$ with $\lambda^\infty(\beta ) \notin \mathcal{S}(\nu) $ and $a_x^\infty(\beta ), a_y^\infty(\beta ), a_z^\infty(\beta ) > 0$.
\end{assumption}
}
We precisely have the following result.

\begin{theorem}\label{thm:semi-circle_dependent}
Let $\tT$ be a sequence of random tensors defined as in \eqref{eq:spiked_tensor_model}. Under Assumptions \ref{ass:growth} and \ref{ass:lambda_outside}, the empirical spectral measure of $\Phi_3(\tT, \vu_*, \vv_*, \vw_*)$ converges weakly almost surely to a deterministic measure $\nu$ whose Stieltjes transform $g(z)$ is defined as the solution to the equation $g(z) = \sum_{i=1}^3 g_i(z)$ such that $\Im[g(z)]>0$ for $\Im[z]>0$ where, for $i\in [3]$ $g_i(z)$ satisfies $g_i^2(z) - (g(z) + z) g_i(z) - c_i = 0$ for $z\in \sC\setminus \mathcal{S}(\nu)$.
\end{theorem}

\begin{proof}
See Appendix \ref{proof:semi-circle_dependent}.
\end{proof}

\begin{figure}[t!]
    \centering
    \includegraphics[width=14cm]{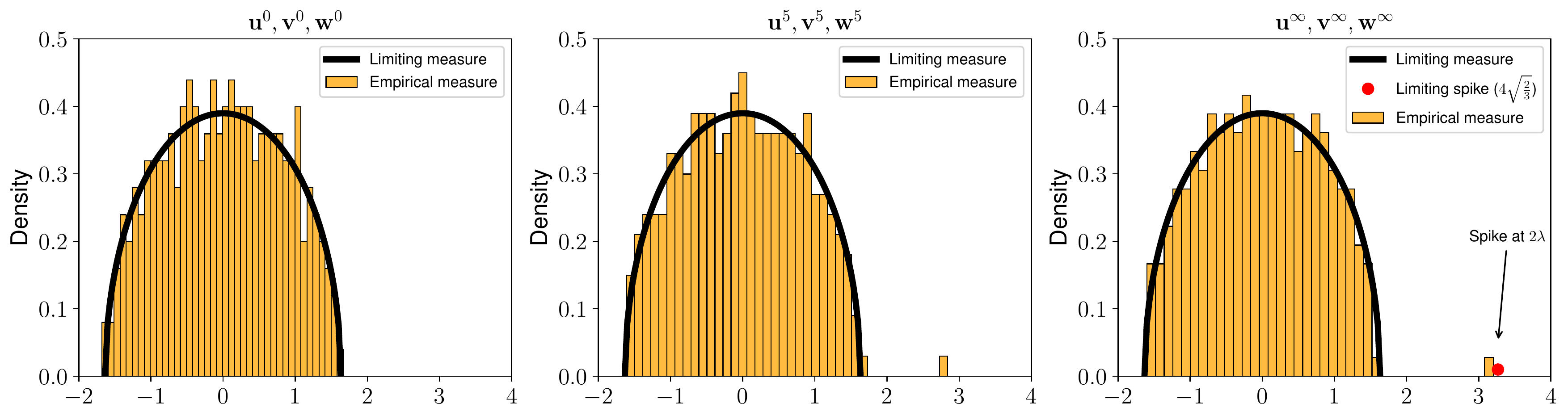}
    \caption{Spectrum of $  \Phi_3(\tT, \vu, \vv, \vw)$ at iterations $0, 5$ and $\infty$ of the power method (see Algorithm \ref{alg:power_method}) applied on $\tT$, for $m=n=p=100$ and $\beta = 0$.}
    \label{fig:spectrum_dependent}
\end{figure}

However, the statistical dependency between $\vu,\vv,\vw$ and $\tX$ exhibits an isolated eigenvalue in the spectrum of $\Phi_3(\tT, \vu, \vv, \vw)$ at the value $2\lambda$ independently of the value of the SNR $\beta$, which is a consequence of \eqref{eq:2lambda}. Figure \ref{fig:spectrum_dependent} depicts iterations of the power iteration method where the leftmost histogram corresponds to the fixed point solution, where one sees that the spike converges to the value $2\lambda$.

\begin{remark}
Note that, from Assumption \ref{ass:lambda_outside}, the almost sure limit $\lambda^\infty$ of $\lambda$ must lie outside the support $\mathcal{S}(\nu)$ of the deterministic measure $\nu$ described by Theorem \ref{thm:semi-circle_dependent}. The same phenomenon was noticed in \cite{de2021random} where it is assumed that the almost sure limit $\mu^\infty$ of $\mu$ must satisfy $\mu^\infty>2\sqrt{\frac{2}{3}}$ in the case of symmetric tensors.
\end{remark}

Let us denote the blocks of the resolvent of $\Phi_3(\tT, \vu, \vv, \vw)$ as
\begin{align}
\mR_{\Phi_3(\tT, \vu, \vv, \vw)}(z) = \begin{bmatrix}
\mR^{11}(z) & \mR^{12}(z) & \mR^{13}(z)\\
\mR^{12}(z)^\top & \mR^{22}(z) & \mR^{23}(z)\\
\mR^{13}(z)^\top & \mR^{23}(z)^\top & \mR^{33}(z)
\end{bmatrix},
\end{align}
where $\mR^{11}(z)\in \sM_m, \mR^{22}(z)\in \sM_n$ and $\mR^{33}(z)\in \sM_p$. The following corollary of Theorem \ref{thm:semi-circle_dependent} will be useful subsequently.
\begin{corollary}\label{cor:tracesofR} Recall the setting and notations of Theorem \ref{thm:semi-circle_dependent}. For all $z\in \sC\setminus \mathcal{S}(\nu)$, we have
\begin{align*}
\frac1N \tr \mR^{11}(z) \asto g_1(z), \quad \frac1N \tr \mR^{22}(z) \asto g_2(z), \quad \frac1N \tr \mR^{33}(z) \asto g_3(z).
\end{align*}
\end{corollary}

\subsection{Concentration of the singular value and the alignments}\label{sec:concentration}
When the dimensions of $\tT$ grow large under Assumption \ref{ass:growth}, the singular value $\lambda$ and the alignments $\langle \vu, \vx\rangle, \langle \vv, \vy\rangle$ and $\langle \vw, \vz\rangle$ converge almost surely to some deterministic limits. This can be shown by controlling the variances of these quantities using the Poincaré's inequality (Lemma \ref{lem:poincaré}). Precisely, for $\lambda$, invoking \eqref{eq:deriv_lambda} we have
\begin{align*}
\Var \lambda \leq \sum_{ijk} \EE \left \vert \frac{\partial \lambda}{\partial X_{ijk}} \right \vert^2 = \frac1N \sum_{ijk} u_i^2 v_j^2 w_k^2 = \frac1N.
\end{align*}
Bounding higher order moments of $\lambda$ similarly allows to obtain the concentration of $\lambda$, e.g. by Chebyshev's inequality, we have for all $t>0$
\begin{align*}
\mathbb{P} \left( \vert \lambda - \EE\lambda \vert \geq t \right) \leq \frac{1}{N\, t^2}.
\end{align*}
Similarly, with \eqref{eq:deriv}, there exists $C>0$ such that for all $t>0$
\begin{align*}
\mathbb{P} \left( \vert \langle \vu, \vx\rangle - \EE\langle \vu, \vx\rangle \vert \geq t \right) \leq \frac{C}{N\, t^2},\,\, \mathbb{P} \left( \vert \langle \vv, \vy\rangle - \EE\langle \vv, \vy\rangle \vert \geq t \right) \leq \frac{C}{N\, t^2}, \,\, \mathbb{P} \left( \vert \langle \vw, \vz\rangle - \EE\langle \vw, \vz\rangle \vert \geq t \right) \leq \frac{C}{N\, t^2}.
\end{align*}
For the remainder of the manuscript, we denote the almost sure limits of $\lambda$ and the of alignments $\langle \vu, \vx\rangle, \langle \vv, \vy\rangle$ and $\langle \vw, \vz\rangle$ respectively as 
\begin{align}\label{eq:definition_ax_ay_az}
\lambda^\infty(\beta) \equiv \lim_{N\to \infty} \lambda,\quad a_\vx^\infty(\beta) \equiv \lim_{N\to \infty} \langle \vu, \vx\rangle,\quad a_\vy^\infty(\beta) \equiv \lim_{N\to \infty} \langle \vv, \vy\rangle,\quad a_\vz^\infty(\beta) \equiv \lim_{N\to \infty} \langle \vw, \vz\rangle.
\end{align}
\subsection{Asymptotic singular value and alignments}
Having the concentration result from the previous subsection, it remains to estimate the expectations $\EE\lambda, \EE\langle \vu, \vx\rangle, \EE\langle \vv, \vy\rangle$ and $\EE\langle \vw, \vz\rangle$. Usually, the evaluation of these quantities using tools from random matrix theory relies on computing Cauchy integrals involving $\mR_{\Phi_3(\tT, \vu, \vv, \vw)}(z)$ (the resolvent of $\Phi_3(\tT, \vu, \vv, \vw)$). Here, we take a different (yet analytically simpler) approach by taking directly the expectation of the identities in \eqref{eq:singular_vectors} and \eqref{eq:singular_value}, then applying Stein's Lemma (Lemma~\ref{lemma:stein}) and Lemma \ref{lemma_differentiable}. For instance, for $\lambda$, we have
\begin{align*}
\EE \lambda = \frac{1}{\sqrt N} \sum_{ijk} \EE \left[ v_j w_k \frac{\partial u_i }{\partial X_{ijk}} \right] + \EE \left[ u_i w_k \frac{\partial v_j }{\partial X_{ijk}} \right] + \EE \left[ u_i v_j \frac{\partial w_k }{\partial X_{ijk}} \right] + \beta \EE \left[ \langle \vu, \vx\rangle \langle \vv, \vy\rangle \langle \vw, \vz\rangle  \right].
\end{align*}
From \eqref{eq:deriv}, when $N\to \infty$, it turns out that the only contributing terms\footnote{Yielding non-vanishing terms in the expression of $\lambda^\infty(\beta)$.} of the derivatives $\frac{\partial u_i }{\partial X_{ijk}}, \frac{\partial v_j }{\partial X_{ijk}}$ and $\frac{\partial w_k }{\partial X_{ijk}}$ in the above sum are respectively
\begin{align*}
\frac{\partial u_i }{\partial X_{ijk}} \simeq -\frac{1}{\sqrt N} v_j w_k R^{11}_{ii}(\lambda), \quad \frac{\partial v_j }{\partial X_{ijk}} \simeq -\frac{1}{\sqrt N} u_i w_k R^{22}_{jj}(\lambda), \quad \frac{\partial w_k }{\partial X_{ijk}} \simeq -\frac{1}{\sqrt N} u_i v_j R^{33}_{kk}(\lambda).
\end{align*} 
This yields
\begin{align*}
\EE \lambda = -\frac1N \left( \tr \mR^{11}(\lambda) + \tr \mR^{22}(\lambda) + \tr \mR^{33}(\lambda) \right)+ \beta \EE \left[ \langle \vu, \vx\rangle \langle \vv, \vy\rangle \langle \vw, \vz\rangle  \right] + \mathcal{O}\left( N^{-1} \right).
\end{align*}
Therefore, the almost sure limit $\lambda^\infty(\beta)$ as $N\to \infty$ of $\lambda$ satisfies 
\begin{align*}
\lambda^{\infty}(\beta) + g(\lambda^{\infty}(\beta)) = \beta a_\vx^\infty(\beta) a_\vy^\infty(\beta) a_\vz^\infty(\beta),
\end{align*}
where $a_\vx^\infty(\beta), a_\vy^\infty(\beta)$ and $a_\vz^\infty(\beta)$ are defined in \eqref{eq:definition_ax_ay_az}. From \eqref{eq:singular_vectors}, proceeding similarly as above with the identities
\begin{align*}
\vx^\top \tT(\vv)\vw = \lambda \langle \vu, \vx\rangle, \quad \vy^\top \tT(\vu)\vw = \lambda \langle \vv, \vy\rangle, \quad \vz^\top \tT(\vv)^\top \vu = \lambda \langle \vw, \vz\rangle, 
\end{align*}
we obtain the following result.
\begin{theorem}\label{thm:asymptotics}
Recall the notations in Theorem \ref{thm:semi-circle_dependent}. Under Assumptions \ref{ass:growth} and \ref{ass:lambda_outside}, there exists $\beta_s>0$ such that for $\beta > \beta_s$,
\begin{align*}
\begin{cases}
\lambda \asto \lambda^\infty(\beta),\\
\left\vert \langle \vu , \vx \rangle \right\vert \asto \frac{1}{\sqrt{ \alpha_2(\lambda^\infty(\beta)) \alpha_3(\lambda^\infty(\beta)) }}, \\ \left\vert \langle \vv , \vy \rangle \right\vert \asto \frac{1}{\sqrt{ \alpha_1(\lambda^\infty(\beta)) \alpha_3(\lambda^\infty(\beta)) }}, \\\left\vert \langle \vw , \vz \rangle \right\vert \asto \frac{1}{\sqrt{ \alpha_1(\lambda^\infty(\beta)) \alpha_2(\lambda^\infty(\beta)) }},
\end{cases}
\end{align*}
where $\alpha_i(z) \equiv \frac{\beta}{ z + g(z) - g_i(z)   } $ and $\lambda^\infty(\beta)$ satisfies $f(\lambda^\infty(\beta), \beta) = 0$ with $f(z, \beta) = z + g(z) - \frac{\beta}{\alpha_1(z) \alpha_2(z) \alpha_3(z)}$. Besides, for $\beta\in [0, \beta_s]$, $\lambda^\infty$ is bounded\footnote{Such bound might be computed numerically and corresponds to the right edge of the limiting spectral measure in Theorem \ref{thm:semi-circle_dependent}.} by an order one constant and $\left\vert \langle \vu , \vx \rangle \right\vert, \left\vert \langle \vv , \vy \rangle \right\vert, \left\vert \langle \vw , \vz \rangle \right\vert \asto 0$. 
\end{theorem}

\begin{proof}
See Appendix \ref{proof:asymptotics}.
\end{proof}
{\color{cblue}
\begin{remark}
    A more compact expression for the alignments is provided in Theorem \ref{thm:asymptotics_general}.
\end{remark}
}
\begin{figure}[t!]
    \centering
    \includegraphics[width=14cm]{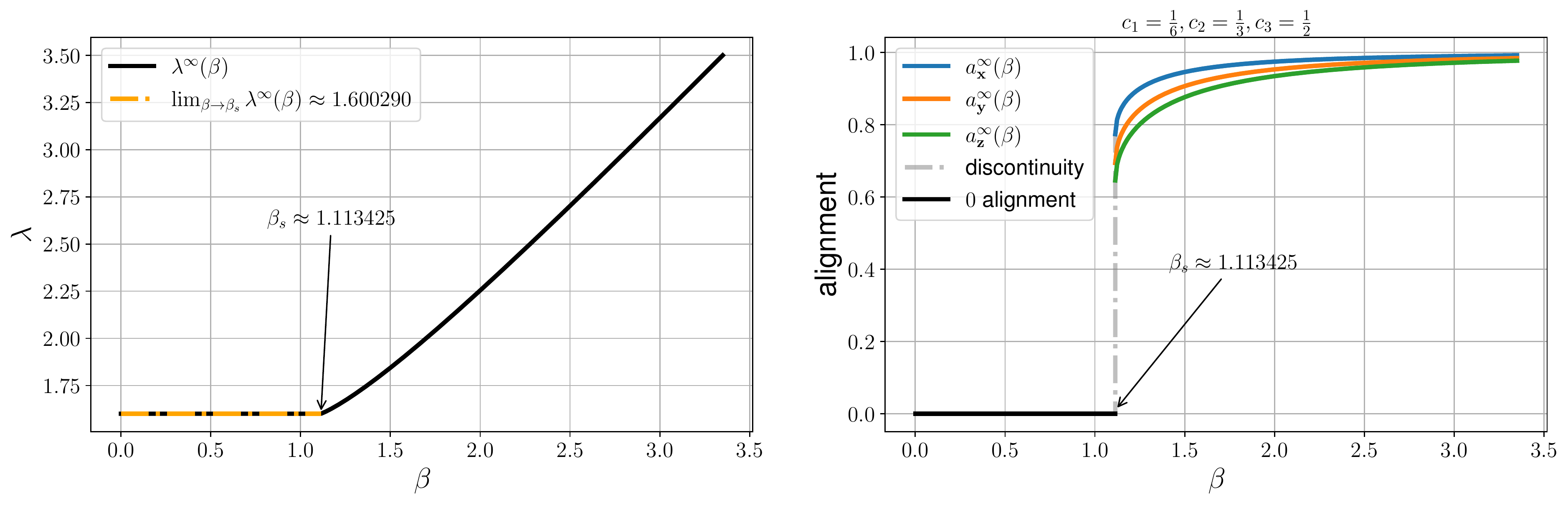}
    \caption{ Asymptotic singular value and alignments of $\tT$ as predicted by Theorem \ref{thm:asymptotics} for $c_1=\frac16, c_2=\frac13$ and $c_3=\frac12$.}
    \label{fig:rectangular_order3}
\end{figure}

Figure \ref{fig:rectangular_order3} depicts the predicted asymptotic dominant singular value $\lambda^\infty(\beta)$ of $\tT$ and the corresponding alignments from Theorem \ref{thm:asymptotics}. As we can see, the result of Theorem \ref{thm:asymptotics} predicts that a non-zero correlation between the population signals $\vx,\vy,\vz$ and their estimated counterparts is possible only when $\beta> \beta_s \approx 1.1134$, which corresponds to the value after which $\lambda^\infty(\beta)$ starts to increase with $\beta$.

\begin{remark} Note that, given $g(z)$, the inverse formula expressing $\beta$ in terms of $\lambda^\infty$ is explicit. Specifically, we have $\beta(\lambda^\infty) = \sqrt{ \frac{ \prod_{i=1}^3 (\lambda^\infty + g(\lambda^\infty) - g_i(\lambda^\infty)) }{ \lambda^\infty + g(\lambda^\infty)} }$. In particular, this inverse formula provides an estimator for the SNR $\beta$ given the largest singular value $\lambda$ of $\tT$.
\end{remark}

\subsection{Cubic $3$-order tensors: case $c_1=c_2=c_3=\frac13$}
In this section, we study the particular case where all the tensor dimensions are equal. As such, the three alignments $\langle \vu, \vx\rangle, \langle \vv, \vy\rangle, \langle \vw, \vz\rangle$ converge almost surely to the same quantity. In this case, the almost sure limits of $\lambda$ and $\langle \vu, \vx\rangle, \langle \vv, \vy\rangle, \langle \vw, \vz\rangle$ can be obtained explicitly in terms of the signal strength $\beta$ as per the following corollary of Theorem \ref{thm:asymptotics}. 

\begin{corollary}\label{cor:cubic} Under Assumptions \ref{ass:lambda_outside} and \ref{ass:growth} with $c_1=c_2=c_3=\frac{1}{3}$, for $\beta > \beta_s = \frac{2\sqrt 3}{3}$
\begin{align*}
\begin{cases}
\lambda \asto \lambda^\infty(\beta) \equiv \sqrt{\frac{\beta^{2}}{2} + 2 + \frac{\sqrt{3} \sqrt{\left(3 \beta^{2} - 4\right)^{3}}}{18 \beta}},\\
\left\vert \langle \vu , \vx \rangle \right\vert , \left\vert \langle \vv , \vy \rangle \right\vert , \left\vert \langle \vw , \vz \rangle \right\vert \asto  \frac{\sqrt{9 \beta^{2} - 12 + \frac{\sqrt{3} \sqrt{\left(3 \beta^{2} - 4\right)^{3}}}{\beta}} + \sqrt{9 \beta^{2} + 36 + \frac{\sqrt{3} \sqrt{\left(3 \beta^{2} - 4\right)^{3}}}{\beta}}}{6\sqrt{2} \beta}.
\end{cases}
\end{align*}
Besides, for $\beta\in\left[0, \frac{2\sqrt3}{3}\right]$, $\lambda \asto \lambda^\infty \leq 2\sqrt{\frac23}$ and $\left\vert \langle \vu , \vx \rangle \right\vert, \left\vert \langle \vv , \vy \rangle \right\vert, \left\vert \langle \vw , \vz \rangle \right\vert \asto 0$.
\end{corollary}
\begin{proof}
See Appendix \ref{proof:spike_align_cubic}.
\end{proof}

\begin{figure}[t!]
    \centering
    \includegraphics[width=14cm]{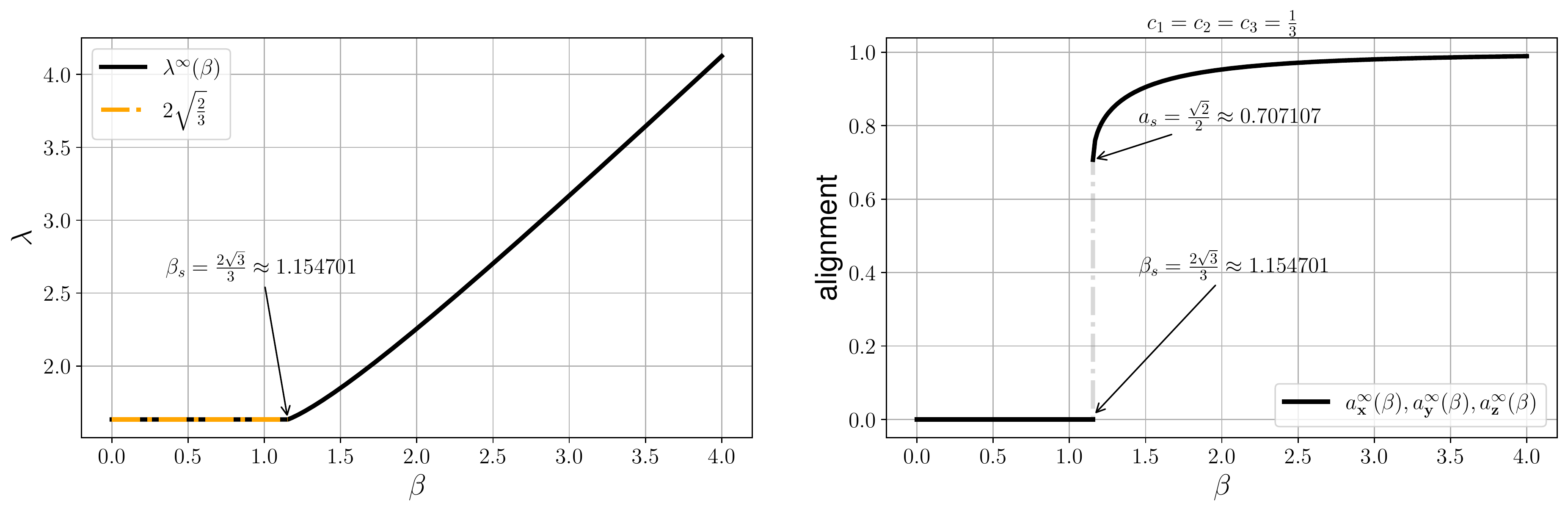}
    \caption{ Asymptotic singular value and alignments of $\tT$ for $c_1=c_2=c_3=\frac13$ as per Corollary \ref{cor:cubic}.}
    \label{fig:cubic_order3}
\end{figure}

Figure \ref{fig:cubic_order3} provides plots of the almost sure limits of the singular value and alignments when $\tT$ is cubic as per Corollary \ref{cor:cubic} (see Subsection \ref{sec:simus_order_3} for simulations supporting the above result). In particular, this result predicts a possible correlation between the singular vectors and the underlying signal components above the value $\beta_s = \frac{2\sqrt{3}}{3} \approx 1.154$ with corresponding singular value $\lambda_s = 2\sqrt{\frac{2}{3}}\approx 1.633$ with an alignment $a_s = \frac{\sqrt 2}{2}\approx 0.707$. In addition, we can easily check from the formulas above that $\frac{\lambda^\infty(\beta)}{\beta} \to 1$ and $ \left\vert \langle \vu , \vx \rangle \right\vert \to 1  $ for large values of $\beta$. Besides, for values of $\beta$ around the value $\beta_s$, the expression of $\lambda^\infty(\beta)$ admits the following expansion
\begin{align*}
\lambda^\infty(\beta) = 2\sqrt{\frac23} + \frac{\sqrt2}{4} \left(\beta - \beta_s \right) + \frac{\sqrt2\, 3^{\frac14}}{4} \left(\beta - \beta_s \right)^{\frac32} + \frac{3\sqrt6}{64} \left(\beta - \beta_s \right)^2 + o(\left(\beta - \beta_s \right)^2),
\end{align*}
whereas the corresponding alignment expends as
\begin{align*}
\frac{\sqrt2}{2} + \frac{\sqrt2\, 3^{\frac14}}{4} \sqrt{\beta - \beta_s} - \frac{\sqrt6}{16} \left(\beta - \beta_s \right) - \frac{\sqrt2 \, 3^{\frac34}}{16} \left(\beta - \beta_s \right)^{\frac32} + \frac{21\sqrt2}{256} \left(\beta - \beta_s \right)^2 + o(\left(\beta - \beta_s \right)^2).
\end{align*}

\section{Random tensors meet random matrices}\label{sec:tensors2matrices}
In this section, we investigate the application of Theorem \ref{thm:asymptotics} to the particular case of spiked random matrices. Indeed, for instance when $p=1$ the spiked tensor model in \eqref{eq:spiked_tensor_model} becomes a spiked matrix model which we will now denote as
\begin{align}\label{eq:spiked_random_matrix}
\mM = \beta \vx \vy^\top + \frac{1}{\sqrt{N}} \mX,
\end{align}
where again $\vx$ and $\vy$ are on the unit spheres $\sS^{m-1}$ and $\sS^{n-1}$ respectively, $N=m+n$ and $\mX$ is a Gaussian noise matrix with i.i.d. entries $X_{ij}\sim \mathcal{N}(0, 1)$.

Our approach does not apply directly to the matrix $\mM$ above. Indeed, the singular vectors $\vu, \vv$ of $\mM$ corresponding to its largest singular value $\lambda$ satisfy the identities
\begin{align}
\mM \vv = \lambda \vu, \quad \mM^\top \vu = \lambda \vv, \quad \lambda = \vu^\top \mM \vv,
\end{align}
and deriving $\vu, \vv$ w.r.t.\@ an entry $X_{ij}$ of $\mX$ will result in
\begin{align}
\left( \begin{bmatrix}
\vzero_{m\times m} & \mM \\
\mM^\top & \vzero_{n\times n}
\end{bmatrix} - \lambda \mI_N \right) \begin{bmatrix}
\frac{\partial \vu}{\partial X_{ij}} \\
\frac{\partial \vv}{\partial X_{ij}}
\end{bmatrix} = - \frac{1}{\sqrt N} \begin{bmatrix}
 v_j \left( \ve_i^m - u_i \vu \right) \\
u_i \left( \ve_j^n - v_j \vv \right)
\end{bmatrix}.
\end{align}
Now since $\lambda$ is an eigenvalue of $\begin{bmatrix}
\vzero_{m\times m} & \mM \\
\mM^\top & \vzero_{n\times n}
\end{bmatrix}$, the matrix $\left( \begin{bmatrix}
\vzero_{m\times m} & \mM \\
\mM^\top & \vzero_{n\times n}
\end{bmatrix} - \lambda \mI_N \right)$ is not invertible, which makes our approach inapplicable for $d=2$. However, we can retrieve the behavior of the spiked matrix model in \eqref{eq:spiked_random_matrix} by considering an order $3$ tensor and setting for instance $c_3\to 0$ as we will discuss subsequently. 

\subsection{Limiting spectral measure}
Given the result of Theorem \ref{thm:asymptotics}, the asymptotic singular value and alignments for the spiked matrix model in \eqref{eq:spiked_random_matrix} correspond to the particular case $c_3\to 0$. We start by characterizing the corresponding limiting spectral measure, we have the following corollary of Theorem \ref{thm:semi-circle_dependent} when $c_3\to0$. 

\begin{corollary}\label{cor:spectrum_matrices} Given the setting and notations of Theorem \ref{thm:semi-circle_independent}, under Assumptions \ref{ass:lambda_outside} and \ref{ass:growth} with $c_1=c, c_2=1-c$ for $c\in (0, 1)$, the empirical spectral measure of $\frac{1}{\sqrt{N}}\Phi_3(\tX, \va, \vb, \vc)$ converges weakly almost surely to a deterministic measure $\nu$ defined on the support $\mathcal{S}(\nu) \equiv \left[-\sqrt{1 + 2\sqrt{\eta}},  -\sqrt{1 - 2\sqrt{\eta}} \right]\cup \left[\sqrt{1 - 2\sqrt{\eta}},  \sqrt{1 + 2\sqrt{\eta}} \right]$ with $\eta = c(1-c)$, whose density function writes as
\begin{align*}
\nu(dx) = \frac{1}{\pi x} \sin\left( \frac{\arctan_2(0, q_c(x))}{2} \right) \sqrt{ \vert q_c(x) \vert } \sign \left( \frac{ \sin\left( \frac{\arctan_2(0, q_c(x))}{2} \right) }{x} \right)dx + \left( 1 - 2\min(c,1-c) \right)\delta(x),
\end{align*}
where $q_c(x) = (x^2 - 1)^2 + 4c(c-1)$. And the corresponding Stieltjes transform writes as
\begin{align*}
g(z) = -z + \frac{\sqrt{ q_c(z) }}{z}, \quad \text{for}\quad z\in \sC\setminus \mathcal{S}(\nu).
\end{align*}
\end{corollary}
\begin{proof}
See Appendix \ref{proof:spectrum_matrices}.
\end{proof}

\begin{figure}[t!]
    \centering
    \includegraphics[width=12cm]{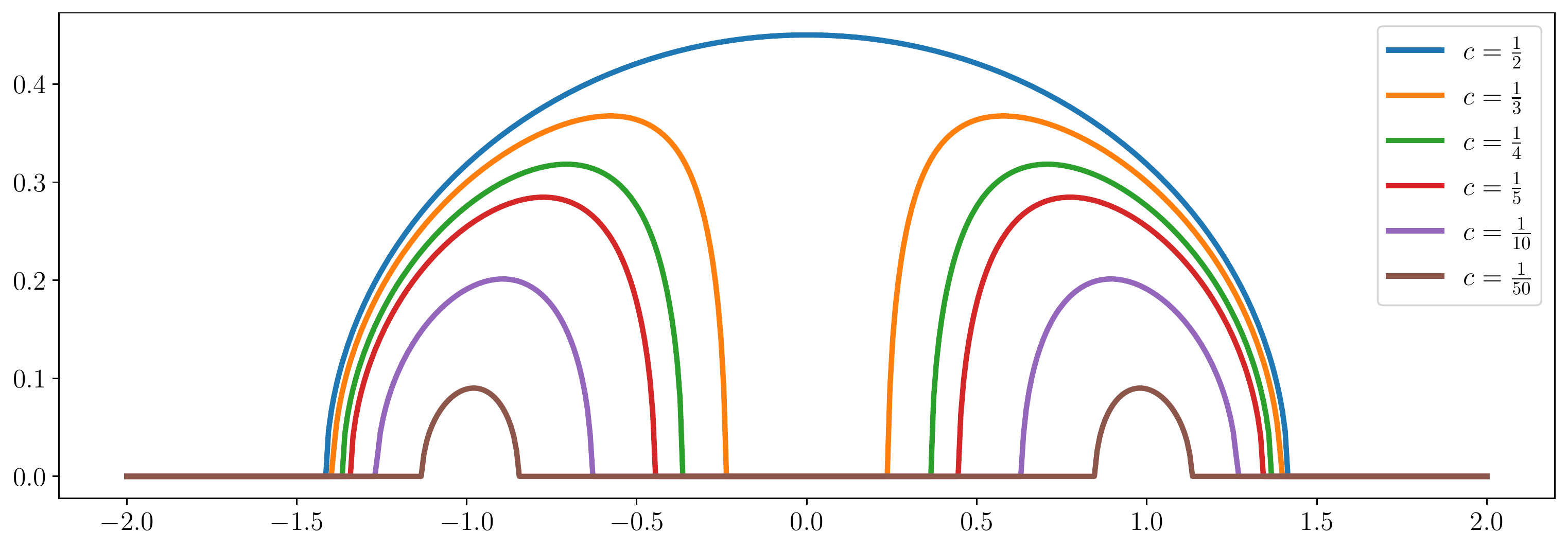}
    \caption{ Limiting spectral measure of $\frac{1}{\sqrt{N}}\Phi_3(\tX, \va, \vb, \vc)$ for $c_1=c, c_2=1-c$ as per Corollary \ref{cor:spectrum_matrices}.}
    \label{fig:matrix_density}
\end{figure}

\begin{remark}
Corollary \ref{cor:spectrum_matrices} describes the limiting spectral measure of $\frac{1}{\sqrt{m+n}}\Phi_3(\tX, \va, \vb, \vc)$ for $\tX\sim \sT_{m,n,1}(\mathcal{N}(0, 1))$, $(\va, \vb)\in \sS^{m-1}\times\sS^{n-1}$ and $\vc=1$ (a scalar), which is equivalent to random matrices of the form $\frac{1}{\sqrt{m+n}}\begin{bmatrix}
\vzero_{m\times m} & \mX \\
\mX^\top & \vzero_{n\times n}
\end{bmatrix}$ with $\mX\sim \sM_{m,n}(\mathcal{N}(0, 1))$.
\end{remark}
Figure \ref{fig:matrix_density} depicts the limiting spectral measures as per Corollary \ref{cor:spectrum_matrices} for various values of $c$ (see also simulations in Appendix \ref{sec:simus_matrix}). In particular, the limiting measure is a semi-circle law for square matrices (i.e., $c=\frac12$), while it decomposes into a two-mode distribution\footnote{With a Dirac $\delta(x)$ at $0$ weighted by $1-2\min(c, 1-c)$, not shown in Figure \ref{fig:matrix_density}. } for $c\in (0, 1)\setminus\{\frac12\}$ due to the fact that the underlying random matrix model is not full rank ($\Phi_3(\tX, \va, \vb, 1)$ if of rank $2\min(m,n)$).

\subsection{Limiting singular value and alignments}
Given the limiting Stieltjes transform from the previous subsection (Corollary \ref{cor:spectrum_matrices}), the limiting largest singular value of $\mM$ (in \eqref{eq:spiked_random_matrix}) and the alignments of the corresponding singular vectors are obtained thanks to Theorem \ref{thm:asymptotics}, yielding the following corollary.
\begin{corollary}\label{cor:spike_align_matrices} Under Assumption \ref{ass:growth} with $c_1=c, c_2=1-c$, we have for $\beta> \beta_s = \sqrt[4]{c(1-c)} $
\begin{align*}
\begin{cases}
\lambda \asto \lambda^\infty(\beta) = \sqrt{\beta^2 + 1 + \frac{c(1-c)}{\beta^2}},\\
\left\vert \langle \vu, \vx \rangle\right\vert \asto \frac{1}{\kappa(\beta, c)},\quad 
\left\vert \langle \vv, \vy \rangle\right\vert \asto \frac{1}{\kappa(\beta, 1 - c)},
\end{cases}
\end{align*}
where $\vu\in \sS^{m-1}, \vv\in \sS^{n-1}$ are the singular vectors of $\mM$ corresponding to its largest singular value $\lambda$ and $\kappa(\beta, c)$ is given by
\begin{align*}
\kappa(\beta, c) = \beta \sqrt{\frac{   \beta^{2} \left(\beta^{2} + 1\right) - c \left(c - 1\right)  }{ (\beta^4 + c(c-1)) \left( \beta^{2} + 1 -  c \right)}}.
\end{align*}
Besides, for $\beta\in[0, \beta_s]$, $\lambda\asto \sqrt{1 + 2\sqrt{c(1-c)}} $ and $\left\vert \langle \vu, \vx \rangle\right\vert, \left\vert \langle \vv, \vy \rangle\right\vert \asto 0$.
\end{corollary}
\begin{proof}
See Appendix \ref{proof:spike_align_matrices}.
\end{proof}

\begin{figure}[t!]
    \centering
    \includegraphics[width=14cm]{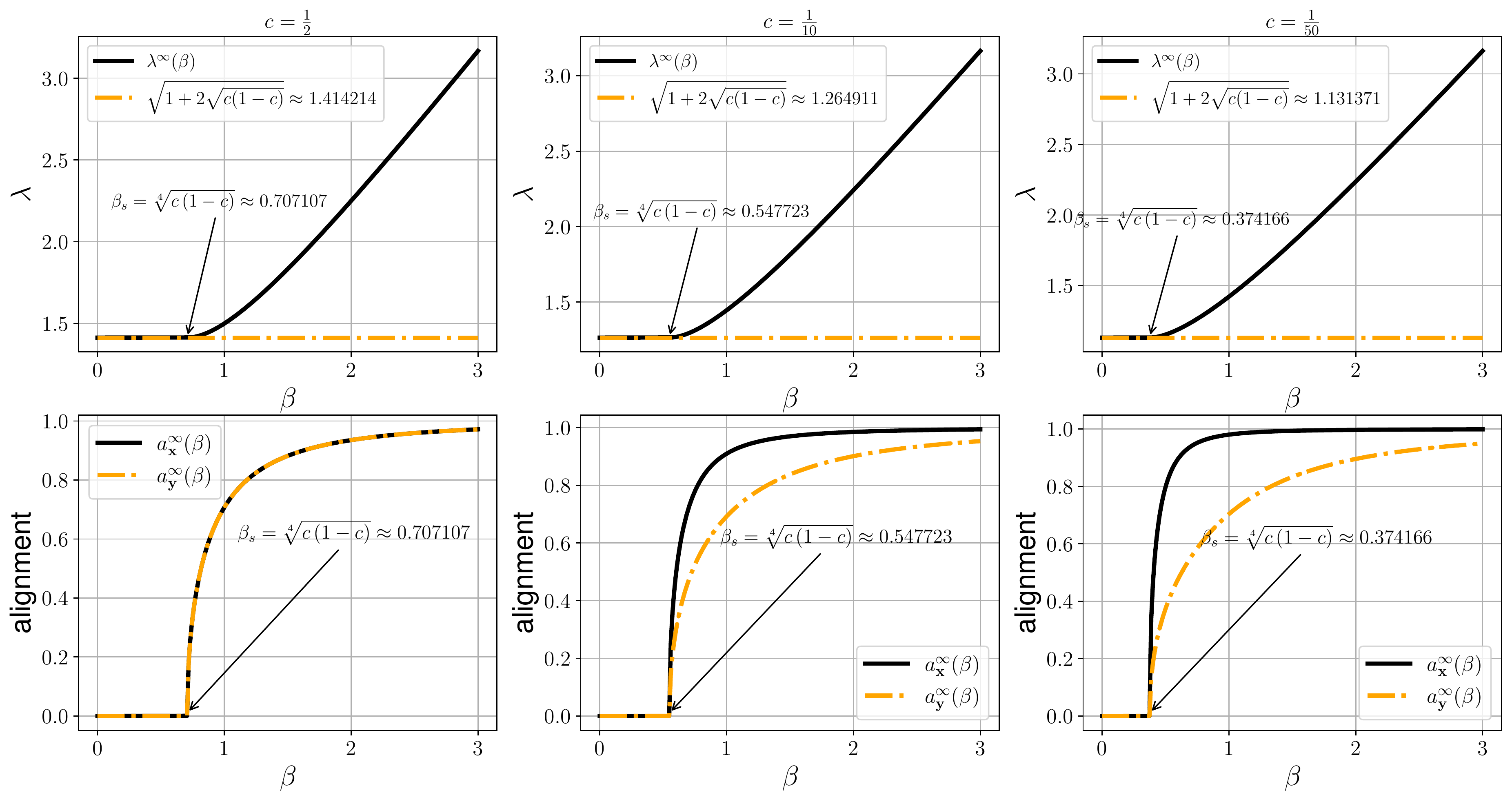}
    \caption{ Asymptotic singular value and alignments of $\tT$ (and of $\mM$ in \eqref{eq:spiked_random_matrix}) for $c_1=c$ and $c_2=1-c$ as per Corollary \ref{cor:spike_align_matrices}, for different values of $c\in \left\{\frac{1}{2}, \frac{1}{10}, \frac{1}{50} \right\}$.}
    \label{fig:spike_align_matrices}
\end{figure}

Figure \ref{fig:spike_align_matrices} provides the curves of the asymptotic singular value and alignments of Corollary \ref{cor:spike_align_matrices} (see Subsection \ref{sec:simus_matrix} for simulations supporting the above result). Unlike the tensor case from the previous section, we see that the asymptotic alignments $a_\vx^\infty(\beta), a_\vy^\infty(\beta)$ are continuous and a positive alignment is observed for $\beta> \beta_s$ which corresponds to the classical BBP phase transition of spiked random matrix models \cite{baik2005phase, benaych2011eigenvalues, peche2006largest}.

{\color{cblue}
\begin{remark}[Application to tensor unfolding]\label{remark:unfolding} The tensor unfolding method consists in estimating the spike components of a given tensor $\tT = \gamma \vx^{(1)} \otimes \cdots \otimes \vx^{(d)} + \frac{1}{\sqrt n}\tX \in \sT_n^d$ with $\tX\sim\sT_n^d(\mathcal{N}(0, 1))$ by applying an SVD to its unfolded matrices (for $i\in [d]$) $\mat_i(\tT) = \gamma \vx^{(i)} (\vy^{(i)})^\top + \frac{1}{\sqrt n}\mX_i $ of size $n\times n^{d-1}$. Applying Corollary \ref{cor:spike_align_matrices} to this case predicts a phase transition at $\beta_s = \sqrt[4]{c(1-c)}$ with\footnote{\color{cblue} In our assumptions, we assume that $c$ is some constant that does not depend on the tensor dimensions. Still, we believe that it can be relaxed in the same vein as in \cite{arous2019landscape}.} $c=\frac{n}{n + n^{d-1}} = \frac{1}{1 + n^{d-2}} \to 0$, which yields $\beta_s = \sqrt[4]{ \frac{n^{d-2}}{(1 + n^{d-2})^2} }$. After re-scaling the noise component (by multiplying \eqref{eq:spiked_random_matrix} by $\sqrt{1 + \frac{m}{n}}$ with $m=n^{d-1}$) this yields $\gamma > n^{\frac{d-2}{4}} $, i.e., the phase transition for tensor unfolding obtained by \cite{arous2021long} (see Theorem 3.3 therein). More generally, for any order-$d$ tensor $\tT$ of arbitrary dimensions $n_1\times \cdots \times n_d$, the tensor unfolding method succeeds provided that $\beta \geq \left( \prod_{i=1}^d n_i \right)^{1/4} / \sqrt{\sum_{i=1}^d n_i}$.
\end{remark}
}
\section{Generalization to arbitrary $d$-order tensors}\label{sec:order_d}
We now show that our approach can be generalized straightforwardly to the $d$-order spiked tensor model of \eqref{eq_asymmetric_spike_model}. Indeed, from \eqref{eq:order_d_singular}, the $\ell_2$-singular value $\lambda$ and vectors $\vu^{(1)}, \ldots, \vu^{(d)} \in \sS^{n_1-1}\times \cdots \times \sS^{n_d-1}$, corresponding to the best rank-one approximation $\lambda \vu^{(1)} \otimes \ldots \otimes \vu^{(d)}$ of the $d$-order tensor $\tT$ in \eqref{eq_asymmetric_spike_model}, satisfy the identities
\begin{align}\label{eq:singular_vectors_id_order_d}
\begin{cases}
\tT(\vu^{(1)}, \ldots, \vu^{(i-1)}, :, \vu^{(i+1)}, \ldots, \vu^{(d)}) = \lambda \vu^{(i)},\\
\lambda = \tT(\vu^{(1)}, \ldots, \vu^{(d)}) = \sum_{i_1, \ldots, i_d} u_{i_1}^{(1)}\ldots u_{i_d}^{(d)} T_{i_1, \ldots, i_d}.
\end{cases}
\end{align}
\subsection{Associated random matrix ensemble}
Let $\tT^{ij}\in \sM_{n_i, n_j}$ denote the matrix obtained by contracting the tensor $\tT$ with the singular vectors in $\{ \vu^{(1)}, \ldots, \vu^{(d)} \} \setminus \{ \vu^{(i)}, \vu^{(j)} \}  $, i.e.,
\begin{align}
\tT^{ij} \equiv \tT( \vu^{(1)}, \ldots, \vu^{(i-1)}, :, \vu^{(i+1)}, \ldots, \vu^{(j-1)}, :, \vu^{(j+1)}, \ldots, \vu^{(d)} ).
\end{align}
As in the order $3$ case, from \eqref{eq:singular_vectors_id_order_d}, the derivatives of the singular vectors $\vu^{(1)}, \ldots, \vu^{(d)}$ with respect to the entry $X_{i_1, \ldots, i_d}$ of the noise tensor $\tX$ express as
\begin{scriptsize}
\begin{align}\label{eq:deriv_order_d}
\begin{bmatrix}
\frac{\partial \vu^{(1)}}{\partial X_{i_1, \ldots, i_d}} \\
\vdots \\
\frac{\partial \vu^{(d)}}{\partial X_{i_1, \ldots, i_d}}
\end{bmatrix} = - \frac{1}{\sqrt{N}} \left( \begin{bmatrix}
\vzero_{n_1\times n_1} & \tT^{12} & \tT^{13} & \cdots & \tT^{1d}\\
(\tT^{12})^\top & \vzero_{n_2\times n_2} & \tT^{23} & \cdots & \tT^{2d}\\
(\tT^{13})^\top & (\tT^{23})^\top & \vzero_{n_3\times n_3} & \ldots & \tT^{3d}\\
\vdots & \vdots & \vdots & \ddots & \vdots \\
(\tT^{1d})^\top & (\tT^{2d})^\top & (\tT^{3d})^\top & \cdots & \vzero_{n_d\times n_d}
\end{bmatrix} - \lambda \mI_N\right)^{-1} \begin{bmatrix}
\prod_{\ell \in \{2, \ldots, d \}} u_{i_\ell}^{(\ell)} (\ve_{i_1}^{n_1} - u_{i_1}^{(i_1)} \vu^{(i_1)})  \\
\vdots \\
\prod_{\ell \in \{1, \ldots, d-1 \}} u_{i_\ell}^{(\ell)} (\ve_{i_d}^{n_d} - u_{i_d}^{(i_d)} \vu^{(i_d)})
\end{bmatrix}
\end{align}
\end{scriptsize}
where $N = \sum_{i \in [d]} n_i$, and the derivative of $\lambda$ w.r.t.\@ $X_{i_1, \ldots, i_d}$ writes as
\begin{align}
\frac{\partial \lambda}{\partial X_{i_1, \ldots, i_d}} = \frac{1}{\sqrt N} \prod_{\ell \in [d]} u_{i_\ell}^{(\ell)}.
\end{align}
As such, the \textit{associated} random matrix model of $\tT$ is the matrix appearing in the \textit{resolvent} in \eqref{eq:deriv_order_d}. More generally, the $d$-order \textit{block-wise tensor contraction ensemble} $\mathcal{B}_d(\tX)$ for $\tX \sim \sT_{n_1, \ldots, n_d}(\mathcal{N}(0, 1))$ is defined as
\begin{align}
\mathcal{B}_d(\tX) \equiv \left\lbrace \Phi_d(\tX, \va^{(1)}, \ldots , \va^{(d)}) \,\, \big\vert \,\,  (\va^{(1)},\ldots,\va^{(d)}) \in \sS^{n_1-1}\times \cdots \times \sS^{n_d-1} \right\rbrace
\end{align}
where $\Phi_d$ is the mapping 
\begin{equation}\label{eq:phi_d}
    \begin{split}
            \Phi_d: \sT_{n_1,\ldots,n_d}\times \sS^{n_1-1}\times \cdots\times \sS^{n_d-1} & \longrightarrow \sM_{\sum_i n_i} \\
    (\tX, \va^{(1)}, \ldots , \va^{(d)}) & \longmapsto \begin{bmatrix}
\vzero_{n_1\times n_1} & \tX^{12} & \tX^{13} & \cdots & \tX^{1d}\\
(\tX^{12})^\top & \vzero_{n_2\times n_2} & \tX^{23} & \cdots & \tX^{2d}\\
(\tX^{13})^\top & (\tX^{23})^\top & \vzero_{n_3\times n_3} & \ldots & \tX^{3d}\\
\vdots & \vdots & \vdots & \ddots & \vdots \\
(\tX^{1d})^\top & (\tX^{2d})^\top & (\tX^{3d})^\top & \cdots & \vzero_{n_d\times n_d}
\end{bmatrix},
    \end{split}
\end{equation}
with $\tX^{ij} \equiv \tX( \va^{(1)}, \ldots, \va^{(i-1)}, :, \va^{(i+1)}, \ldots, \va^{(j-1)}, :, \va^{(j+1)}, \ldots, \va^{(d)} ) \in \sM_{n_i, n_j}$.

{\color{cblue}
\begin{remark}
    As in the order $3$ case, to ensure the existence of the matrix inverse in \eqref{eq:deriv_order_d}, we need to suppose that there exists a tuple $(\lambda_*, \vu_*^{(1)},\ldots,\vu_*^{(d)})$ verifying the identities in \eqref{eq:singular_vectors_id_order_d} such that $\lambda_*$ is not an eigenvalue of $\Phi_d(\tT, \vu_*^{(1)},\ldots,\vu_*^{(d)})$.
\end{remark}
}
\subsection{Limiting spectral measure of block-wise $d$-order tensor contractions}
In this section, we characterize the limiting spectral measure of the ensemble $\mathcal{B}_d(\tX)$ for $\tX \sim \sT_{n_1,\ldots, n_d}(\mathcal{N}(0, 1))$ in the limit when all tensor dimensions grow as per the following assumption.
\begin{assumption}\label{ass:growth_order_d} For all $i\in[d]$, assume that $n_i\to \infty$ with $\frac{n_i}{\sum_j n_j} \to c_i \in (0, 1)$.
\end{assumption}

We thus have the following result which characterizes the spectrum of $\frac{1}{\sqrt N}\Phi_d(\tX , \va^{(1)}, \ldots , \va^{(d)})$ for any deterministic unit norm vectors $\va^{(1)}, \ldots , \va^{(d)}$.

\begin{theorem}\label{thm:semi-circle_independent_d_order}
Let $\tX \sim \sT_{n_1,\ldots, n_d}(\mathcal{N}(0,1))$ be a sequence of random asymmetric Gaussian tensors and $(\va^{(1)}, \ldots , \va^{(d)}) \in \sS^{n_1-1}\times \cdots\times \sS^{n_d-1}$ a sequence of deterministic vectors of increasing dimensions, following Assumption~\ref{ass:growth_order_d}. Then the empirical spectral measure of $\frac{1}{\sqrt N}\Phi_d(\tX , \va^{(1)}, \ldots , \va^{(d)})$ converges weakly almost surely to a deterministic measure $\nu$ whose Stieltjes transform $g(z)$ is defined as the solution to the equation $g(z) = \sum_{i=1}^d g_i(z)$ such that $\Im[g(z)]>0$ for $\Im[z]>0$ where, for $i\in [d]$ $g_i(z)$ satisfies $g_i^2(z) - (g(z) + z) g_i(z) - c_i = 0$ for $z\in \sC\setminus \mathcal{S}(\nu)$.
\end{theorem}

\begin{proof}
See Appendix \ref{proof:circle_independent_d_order}.
\end{proof}

\begin{algorithm}[t!]
\caption{Fixed point equation to compute the Stieltjes transform in Theorem \ref{thm:semi-circle_independent_d_order}}\label{alg:stieltjes_transform}
\begin{algorithmic}
\Require $z\in \sR\setminus \mathcal{S}(\nu)$ and tensor dimension ratios $\vc = [c_1,\ldots, c_d]^\top\in (0, 1)^d$
\Ensure $ g\in \sR,\quad \vg = [g_1, \ldots, g_d]^\top \in \sR^d $\\
\textbf{Repeat}
\State $\quad\vg \gets \frac{g+z}{2} - \frac{\sqrt{4\vc + (g+z)^2}}{2}$ \Comment{Element-wise vector operation.}
\State $\quad g \gets \sum_{i=1}^d g_i$\\
\textbf{until convergence of} $g$
\end{algorithmic}
\end{algorithm}

Algorithm \ref{alg:stieltjes_transform} provides a pseudo-code to compute the Stieltjes transform $g(z)$ in Theorem \ref{thm:semi-circle_independent_d_order} through an iterative solution to the fixed point equation $g(z) = \sum_{i=1}^d g_i(z)$. In particular, as for the order $3$ case, when all the tensor dimensions are equal (i.e., $c_i = \frac1d$ for all $i\in [d]$), the spectral measure of $\frac{1}{\sqrt N}\Phi_d(\tX , \va^{(1)}, \ldots , \va^{(d)})$ converges to a semi-circle law. We have the following corollary of Theorem \ref{thm:semi-circle_independent_d_order}.

\begin{corollary}\label{cor:semi_circle_order_d} With the setting of Theorem \ref{thm:semi-circle_independent_d_order}. Under Assumption \ref{ass:growth_order_d} with $c_i=\frac1d$ for all $i\in[d]$, the empirical spectral measure of $\frac{1}{\sqrt N}\Phi_d(\tX , \va^{(1)}, \ldots , \va^{(d)})$ converges weakly almost surely to the semi-circle distribution support $\mathcal{S}(\nu)\equiv \left[ -2\sqrt{\frac{d-1}{d}} ,2\sqrt{\frac{d-1}{d}}  \right]$, whose density and Stieltjes transform write respectively as
\begin{align*}
\nu(dx) = \frac{d}{2(d-1)\pi} \sqrt{\left( \frac{4(d-1)}{d} -x^2  \right)^+},\quad g(z) \equiv \frac{-zd + d\sqrt{z^2 - \frac{4(d-1)}{d} }}{2(d-1)}, \quad \text{where} \quad z\in \sC\setminus \mathcal{S}(\nu).
\end{align*}
\end{corollary}

\begin{proof}
See Appendix \ref{proof:semi_circle_order_d}.
\end{proof}

\begin{figure}[t!]
    \centering
    \includegraphics[width=14cm]{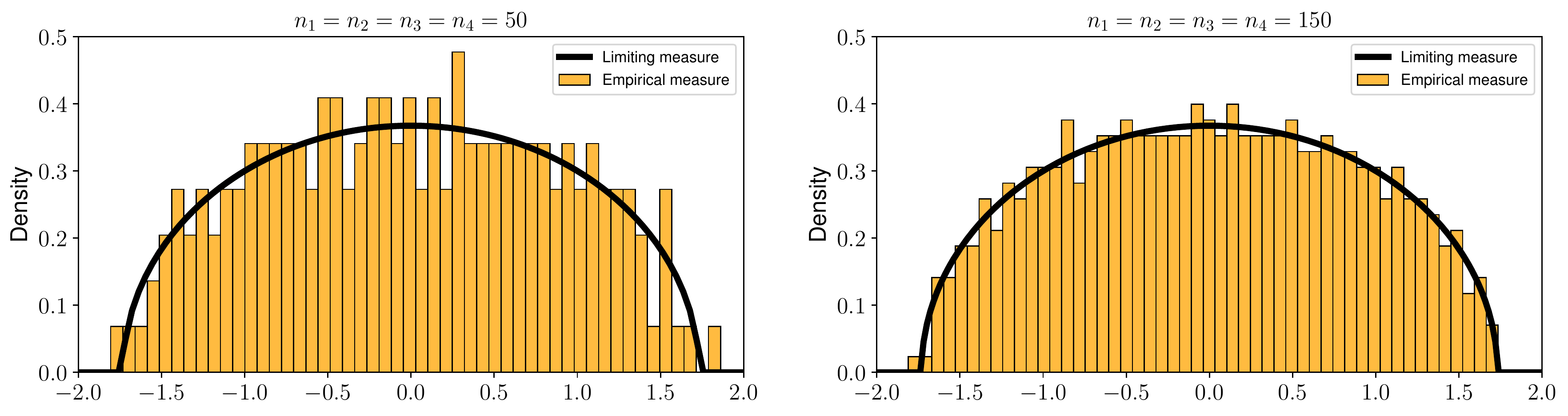}
    \caption{Spectrum of $ \frac{1}{\sqrt N} \Phi_4(\tX, \va, \vb, \vc, \vd)$ with $\tX \sim \sT_{n_1, \ldots, n_4}(\mathcal{N}(0, 1))$ and $\va,\vb,\vc,\vd$ independently sampled from the unit spheres $\sS^{n_1-1}, \ldots, \sS^{n_4-1}$ respectively. In black the semi-circle law as per Theorem \ref{thm:semi-circle_independent_d_order}.}
    \label{fig:spectrum_independent_order_4}
\end{figure}

Figure \ref{fig:spectrum_independent_order_4} provides an illustration of the convergence in law (when the dimension $n_i$ grow large) of the spectrum of $ \frac{1}{\sqrt N} \Phi_4(\tX, \va, \vb, \vc, \vd)$ with a $4$th-order tensor $\tX \sim \sT_{n_1, \ldots, n_4}(\mathcal{N}(0, 1))$ and $\va,\vb,\vc,\vd$ independently sampled from the unit spheres $\sS^{n_1-1}, \ldots, \sS^{n_4-1}$ respectively.

\begin{remark} $\frac{1}{\sqrt N}\Phi_d(\tX , \va^{(1)}, \ldots , \va^{(d)})$ is almost surely of rank $\sum_i \min(n_i, \sum_{j\neq i} n_j)$. As in the order $3$ case, when $\frac{1}{\sqrt N}\Phi_d(\tX , \va^{(1)}, \ldots , \va^{(d)})$ is almost surely full rank, its spectral measure converges to a semi-circle law with connected support, while in general the limiting spectral measure has unconnected support. 
\end{remark}

The result of Theorem \ref{thm:semi-circle_independent_d_order} still holds for $\Phi_d(\tT, \vu^{(1)}, \ldots, \vu^{(d)})$ where $\vu^{(1)}, \ldots, \vu^{(d)}$ stand for the singular vectors of $\tT$ defined through \eqref{eq:singular_vectors_id_order_d}. Indeed, the statistical dependencies between these singular vectors and the noise tensor $\tX$ do not affect the convergence of the spectrum of $\Phi_d(\tT, \vu^{(1)}, \ldots, \vu^{(d)})$ to the limiting measure described by Theorem \ref{thm:semi-circle_independent_d_order}. We further need the following assumption.
{\color{cblue}
\begin{assumption}\label{ass:lambda_outside_order_d}
    We assume that there exists a sequence of critical points $(\lambda_*, \vu_*^{(1)},\ldots, \vu_*^{(d)})$ satisfying \eqref{eq:singular_vectors_id_order_d} such that $\lambda_* \asto \lambda^\infty(\beta)$, $\vert\langle \vu_*^{(i)}, \vx^{(i)} \rangle\vert \asto a_{\vx^{(i)}}^\infty(\beta )$ with $\lambda^\infty(\beta ) \notin \mathcal{S}(\nu) $ and $a_{\vx^{(i)}}^\infty(\beta ) > 0$.
\end{assumption}
}
\begin{theorem}\label{thm:semi-circle_dependent_order_d}
Let $\tT$ be a sequence of spiked random tensors defined as in \eqref{eq_asymmetric_spike_model}. Under Assumptions \ref{ass:growth_order_d} and \ref{ass:lambda_outside_order_d}, the empirical spectral measure of $\Phi_d(\tT, \vu_*^{(1)}, \ldots, \vu_*^{(d)})$ converges weakly almost surely to a deterministic measure $\nu$ whose Stieltjes transform $g(z)$ is defined as the solution to the equation $g(z) = \sum_{i=1}^d g_i(z)$ such that $\Im[g(z)]>0$ for $\Im[z]>0$ where, for $i\in [d]$ $g_i(z)$ satisfies $g_i^2(z) - (g(z) + z) g_i(z) - c_i = 0$ for $z\in \sC\setminus \mathcal{S}(\nu)$.
\end{theorem}

\begin{proof}
See Appendix \ref{proof:semi-circle_dependent_order_d}.
\end{proof}

\begin{figure}[t!]
    \centering
    \includegraphics[width=14cm]{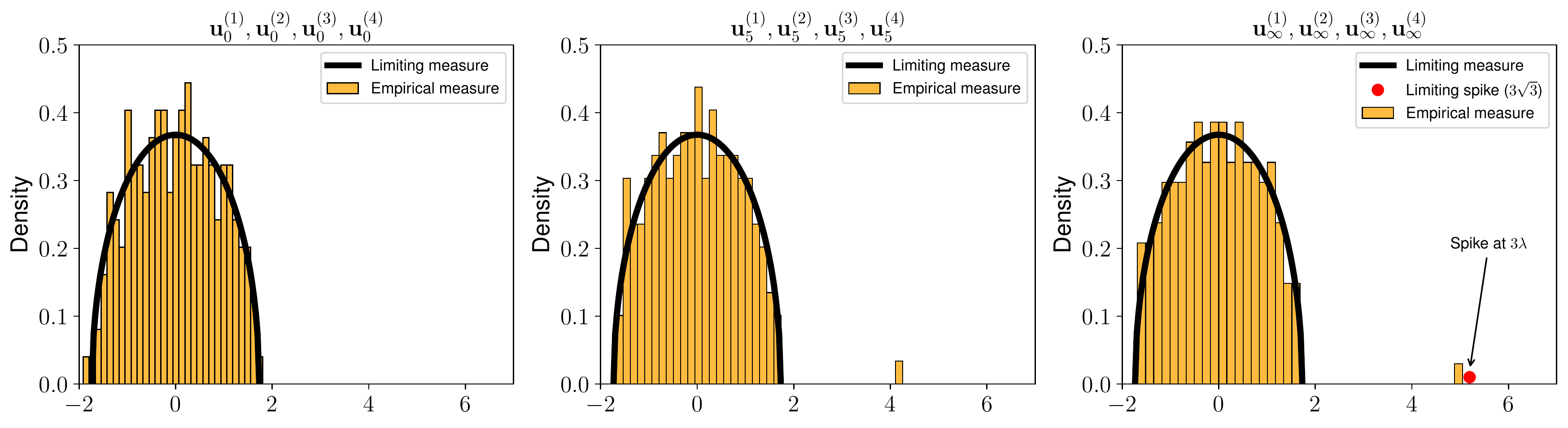}
    \caption{Spectrum of $\Phi_4(\tT, \vu^{(1)}, \ldots, \vu^{(4)})$ at iterations $0,5$ and at convergence of the power method (see Algorithm \ref{alg:power_method}) applied to $\tT$, for $n_1=n_2 =n_3 = n_4 = 50$ and $\beta = 0$.}
    \label{fig:spectrum_dependent_order_4}
\end{figure}

Figure \ref{fig:spectrum_dependent_order_4} depicts the spectrum of $\Phi_4(\tT, \vu^{(1)}, \ldots, \vu^{(4)})$ for an order $4$ tensor $\tT$ with $\beta = 0$. As we saw previously, an isolated eigenvalue pops out from the continuous bulk of $\Phi_4(\tT, \vu^{(1)}, \ldots, \vu^{(4)})$ because of the statistical dependencies between the tensor noise $\tX$ and the singular vectors $\vu^{(1)}, \ldots, \vu^{(4)}$. More generally, for an order-$d$ tensor $\tT$, the spectrum of $\Phi_d(\tT, \vu^{(1)}, \ldots, \vu^{(d)})$ admits an isolated spike at the value $(d-1)\lambda$ independently of the signal strength $\beta$. Indeed, $(d-1)\lambda$ is an eigenvalue of $\Phi_d(\tT, \vu^{(1)}, \ldots, \vu^{(d)})$ with corresponding eigenvector the concatenation of the singular vectors $\vu^{(1)}, \ldots, \vu^{(d)}$, i.e.,
\begin{footnotesize}
\begin{align}
\Phi_d(\tT, \vu^{(1)}, \ldots, \vu^{(d)}) \begin{bmatrix}
\vu^{(1)} \\
\vdots \\
\vu^{(d)}
\end{bmatrix} = \begin{bmatrix}
\sum_{i\neq 1} \tT(:, \vu^{(2)}, \ldots, \vu^{(i-1)}, :, \vu^{(i+1)}, \ldots, \vu^{(d)}) \vu^{(i)} \\
\vdots \\
\sum_{i\neq d} \tT(\vu^{(1)}, \ldots, \vu^{(i-1)}, :, \vu^{(i+1)}, \ldots, \vu^{(d-1)}, :)^\top \vu^{(i)}
\end{bmatrix} = (d-1)\lambda \begin{bmatrix}
\vu^{(1)} \\
\vdots \\
\vu^{(d)}
\end{bmatrix}.
\end{align}
\end{footnotesize}

\subsection{Asymptotic singular value and alignments of hyper-rectangular tensors}
Similarly to the $3$-order case studied previously, when the dimensions of $\tT$ grow large at a same rate, its singular value $\lambda$ and the corresponding alignments $\langle \vu^{(i)}, \vx^{(i)} \rangle$ for $i\in [d]$ concentrate almost surely around some deterministic quantities which we denote $\lambda^\infty(\beta)\equiv \lim_{N\to \infty} \lambda$ and $a_{\vx^{(i)}}^\infty(\beta) \equiv \lim_{N\to \infty} \vert \langle \vu^{(i)}, \vx^{(i)} \rangle \vert $ respectively. Applying again Stein's Lemma (Lemma \ref{lemma:stein}) to the identities in \eqref{eq:singular_vectors_id_order_d}, we obtain the following theorem which characterizes the aforementioned deterministic limits. 
{\color{cblue}
\begin{theorem}\label{thm:asymptotics_general}
Recall the notations in Theorem \ref{thm:semi-circle_dependent_order_d}. Under Assumptions \ref{ass:growth_order_d} and \ref{ass:lambda_outside_order_d}, for $d\geq 3$, there exists $\beta_s>0$ such that for $\beta> \beta_s$,
\begin{align*}
\begin{cases}
\lambda \asto \lambda^\infty(\beta),\\
\left\vert \langle \vx^{(i)}, \vu^{(i)} \rangle \right\vert \asto q_i\left( \lambda^\infty(\beta) \right),
\end{cases}
\end{align*}
where $q_i(z)$ is given by $ q_i(z) = \sqrt{ 1 - \frac{g_i^2(z)}{c_i} }$ and $\lambda^\infty(\beta)$ satisfies $f(\lambda^\infty(\beta),\beta)=0$ with $f(z, \beta) = z + g(z) - \beta \prod_{i=1}^d q_i(z)$. Besides, for $\beta \in [0, \beta_s]$, $\lambda^\infty$ is bounded (in particular when $c_i=\frac1d$ for all $i\in[d]$, $\lambda\asto \lambda^\infty \leq 2 \sqrt{\frac{d-1}{d}}$) and $\left\vert \langle \vx^{(i)}, \vu^{(i)} \rangle \right\vert \asto 0$.
\end{theorem}
}
\begin{proof}
See Appendix \ref{proof:asymptotics_general}.
\end{proof}

\begin{remark} As in the order $3$ case, note that the inverse formula expressing $\beta$ in terms of $\lambda^\infty$ is explicit. Specifically, we have $\beta(\lambda^\infty) =  \frac{ \lambda^\infty + g(\lambda^\infty) }{ \prod_{i=1}^d q_i(\lambda^\infty) }$. In particular, this inverse formula provides an estimator for the SNR $\beta$ given the largest singular value $\lambda$ of $\tT$. Algorithm \ref{alg:alignments} provides a pseudo-code to compute the asymptotic alignments.
\end{remark}

\begin{algorithm}[t!]
\caption{Compute alignments as per Theorem \ref{thm:asymptotics_general}}\label{alg:alignments}
\begin{algorithmic}
\Require SNR $\beta\in\sR_+$ and tensor dimension ratios $\vc = [c_1,\ldots, c_d]^\top\in (0, 1)^d$
\Ensure Asymptotic singular value $\lambda^\infty$ and corresponding alignments $\va = \left[ a_{\vx^{(1)}}^\infty, \ldots, a_{\vx^{(d)}}^\infty \right]^\top\in [0, 1]^d$
\State Set $\lambda^\infty$ as the solution of $f(z, \beta)=0$ in $z$.
\State Set the alignments as $a_{\vx^{(i)}}^\infty \gets \sqrt{ 1 - \frac{g_i^2(\lambda^\infty)}{c_i} }$ with $g_i(z)$ computed by Algorithm \ref{alg:stieltjes_transform} for $z=\lambda^\infty$.
\end{algorithmic}
\end{algorithm}

\begin{figure}[t!]
    \centering
    \includegraphics[width=14cm]{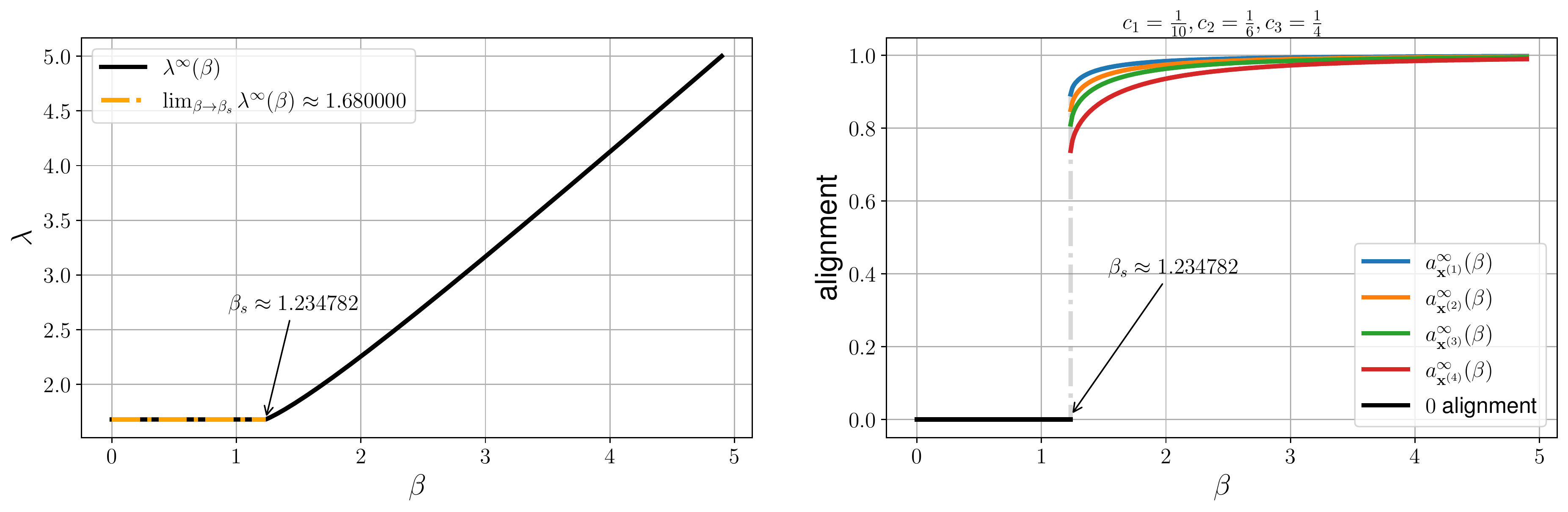}
    \caption{ Asymptotic singular value and alignments of $\tT$ as predicted by Theorem \ref{thm:asymptotics_general} for $c_1=\frac{1}{10}, c_2=\frac16, c_3=\frac14$ and $c_4 = 1-(c_1+c_2+c_3)$.}
    \label{fig:rectangular_order4}
\end{figure}

\begin{figure}[t!]
    \centering
    \includegraphics[width=14cm]{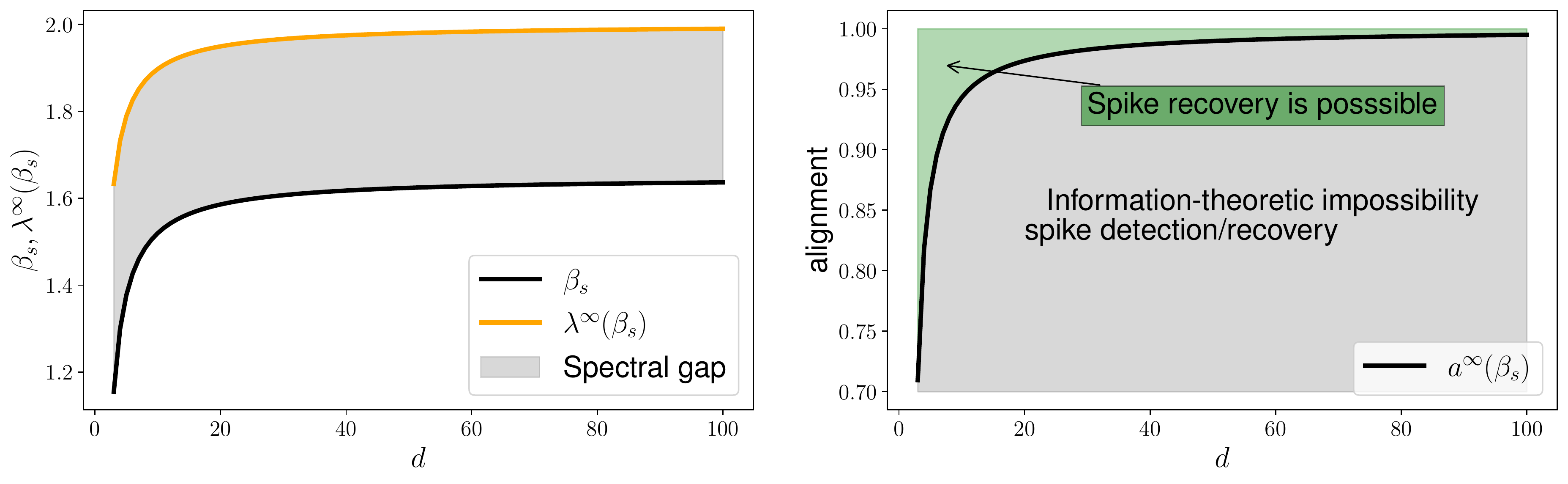}
    \caption{Minimal signal strength $ \beta_s$ and the corresponding singular value and alignments in terms of the tensor order $d$, when all the tensor dimensions are equal.}
    \label{fig:bc_ac_orders}
\end{figure}

Figure \ref{fig:rectangular_order4} depicts the asymptotic singular value and alignments for an order $4$ tensor as per Theorem \ref{thm:asymptotics_general}. As discussed previously, the predicted alignments are discontinuous and a strictly positive correlation between the singular vectors and the underlying signals is possible for the considered local optimum of the ML estimator only above the minimal signal strength $\beta_s\approx 1.234$ in the shown example. In the case of hyper-cubic tensors, i.e., all the dimensions $n_i$ are equal, Figure \ref{fig:bc_ac_orders} depicts the minimal SNR values in terms of the tensor order $d$ and the corresponding asymptotic singular value and alignments. As such, when the tensor order increases the minimal SNR and singular value converge respectively to $\approx 1.6$ and $2$, while the corresponding alignment\footnote{corresponding to the minimal theoretical SNR $\beta_s$.} gets closer to $1$. The expressions of $\beta_s$ and $a^\infty(\beta_s)$ for hyper-cubic tensors of order $d$ are explicitly given as
\begin{align}\color{cblue}
    \beta_s = \sqrt{ \frac{d-1}{d} } \left( 
\frac{d-2}{d-1} \right)^{1 - \frac{d}{2} }, \quad \lim_{\beta \to \beta_s } a^\infty(\beta) = \sqrt{ \frac{d-2}{d-1} }.
\end{align}

\section{Generalization to rank $r$ tensor with orthogonal components}\label{sec:rank_k}
In this section we discuss the generalization of our previous findings to rank $r$ spiked tensor model with $r>1$ of the form
\begin{align}
\tT = \sum_{\ell=1}^r \beta_\ell \, \vx^{(1)}_\ell \otimes \cdots \otimes \vx^{(d)}_\ell + \frac{1}{\sqrt N} \tX,
\end{align}
where $\beta_1> \cdots > \beta_r$ are the signal strengths, $(\vx^{(1)}_\ell , \ldots , \vx^{(d)}_\ell)\in \sS^{n_1 - 1} \times \cdots \times \sS^{n_d - 1}$ for $\ell\in [r]$, $\tX \sim \sT_{n_1,\ldots, n_d}(\mathcal{N}(0,1))$ and $N = \sum_{i=1}^d n_i$. Supposing that the components $\vx^{(i)}_1 , \ldots , \vx^{(i)}_r$ are mutually orthogonal, i.e., $\langle \vx^{(i)}_\ell, \vx^{(i)}_{\ell'} \rangle = 0$ for $\ell \neq \ell'$, the above $r$-rank spiked tensor model can be treated through equivalent rank-one tensors defined for each component $\beta_\ell \, \vx^{(1)}_\ell \otimes \cdots \otimes \vx^{(d)}_\ell$ independently, by applying the result of the rank-one case established in the previous section. Precisely, the best rank-$r$ approximation of $\tT$ corresponds to 
\begin{align} \label{eq_variational_definition_rank_r}
    \argmin_{\lambda_\ell\,\in\, \sR,\, (\vu^{(1)}_\ell,\ldots,\vu^{(d)}_\ell)\,\in\, \sS^{n_1 - 1} \times \cdots \times \sS^{n_d - 1}} \Vert \tT - \sum_{\ell=1}^r \lambda_\ell \, \vu^{(1)}_\ell \otimes  \cdots \otimes \vu^{(d)}_\ell \Vert_{\text{F}}^2,
\end{align}
and the $\vu^{(i)}_\ell$'s correlate with $\vx^{(i)}_\ell$'s as a result of the uniqueness of orthogonal tensor decomposition (see Theorem 4.1 in \cite{anandkumar2014tensor}). Therefore, the study of $\tT$ boils down to the study of the random matrices $\Phi_d(\tT, \vu^{(1)}_\ell,\ldots,\vu^{(d)}_\ell)$ for $\ell \in [r]$ which behave as the rank-one case treated previously with signal strength $\beta_\ell$ respectively. Indeed, since $\vu^{(i)}_\ell \in \sS^{n_i - 1}$ by definition and given the orthogonality condition on the $\vx^{(i)}_\ell$'s, for $\ell'\neq \ell$ the inner product $\langle \vu^{(i)}_\ell, \vx^{(i)}_{\ell'} \rangle \asto 0$ in high dimension, as such 
\begin{align}
\Vert \Phi_d(\tT, \vu^{(1)}_\ell,\ldots,\vu^{(d)}_\ell) - \Phi_d(\tT_\ell, \vu^{(1)}_\ell,\ldots,\vu^{(d)}_\ell) \Vert \asto 0,
\end{align}
where $\tT_\ell\equiv \beta_\ell \, \vx^{(1)}_\ell \otimes \cdots \otimes \vx^{(d)}_\ell + \frac{1}{\sqrt N} \tX$. Meaning, that the study of $\Phi_d(\tT, \vu^{(1)}_\ell,\ldots,\vu^{(d)}_\ell)$ is equivalent to consider the rank-one spiked tensor model $\tT_\ell$.

\section{Discussion}\label{sec:discussion}
In this work, we characterized the asymptotic behavior of spiked asymmetric tensors by mapping them to equivalent (in spectral sense) random matrices. Our starting point is mainly the identities in \eqref{eq:singular_vectors_id_order_d} which are verified by all critical points of the ML problem. Quite surprisingly and as also discussed in \cite{de2021random} for symmetric tensors, we found that our asymptotic equations describe precisely the maximum of the ML problem which correlates with the true spike. Extrapolating the findings from \cite{jagannath2020statistical, de2021random} in the symmetric case, we conjuncture the existence of an order one threshold $\beta_c$ above which our equations describe the behavior of the global maximum of the ML problem. Unfortunately, it is still unclear how we can characterize such $\beta_c$ with our present approach, which remains an open question. The same question concerns the characterization of the algorithmic threshold $\beta_a$ which is more interesting from a practical standpoint as computing the ML solution is NP-hard for $\beta<\beta_a$.

In the present work, our results were derived under a Gaussian assumption on the tensor noise components. We believe that the derived formulas are universal in the sense that they extend to other distributions provided that the fourth order moment is finite (as assumed by \cite{arous2021long} for long random matrices). Other extensions concern the generalization to higher-ranks with arbitrary components and possibly correlated noise components since the present RMT-tools are more flexible than the use of tools form statistical physics.

\begin{acks}[Acknowledgments]
This work was supported by the MIAILargeDATA Chair at University Grenoble Alpes led by R. Couillet and the UGA-HUAWEI LarDist project led by M. Guillaud. We would like to thank Henrique Goulart, Pierre Common and Gérard Ben-Aroud for valuable discussions on the topic of random tensors.
\end{acks}

\appendix
\section{Simulations}\label{sec:simus}
In this section we provide simulations to support our findings.
\subsection{Matrix case}\label{sec:simus_matrix}
\begin{figure}[b!]
    \centering
    \includegraphics[width=14cm]{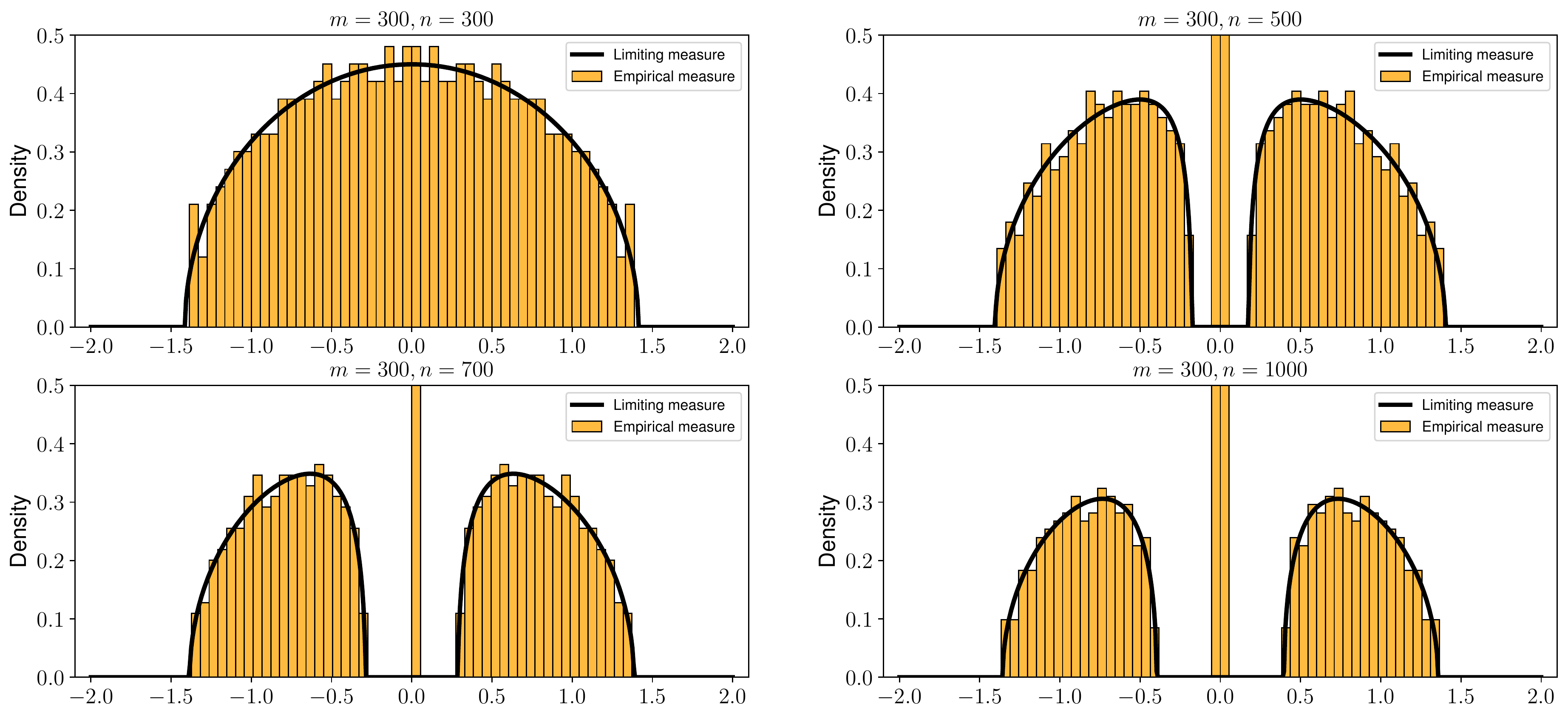}
    \caption{Spectrum of $\begin{bmatrix}
\vzero_{m\times m} & \mM\\
\mM^\top & \vzero_{n\times n} 
\end{bmatrix}$ for $\beta=0$ with different values of the dimensions $m,n$ and the limiting spectral measure in Corollary \ref{cor:spectrum_matrices}.}
    \label{fig:spectrum_matrix}
\end{figure}

\begin{figure}[t!]
    \centering
    \includegraphics[width=14cm]{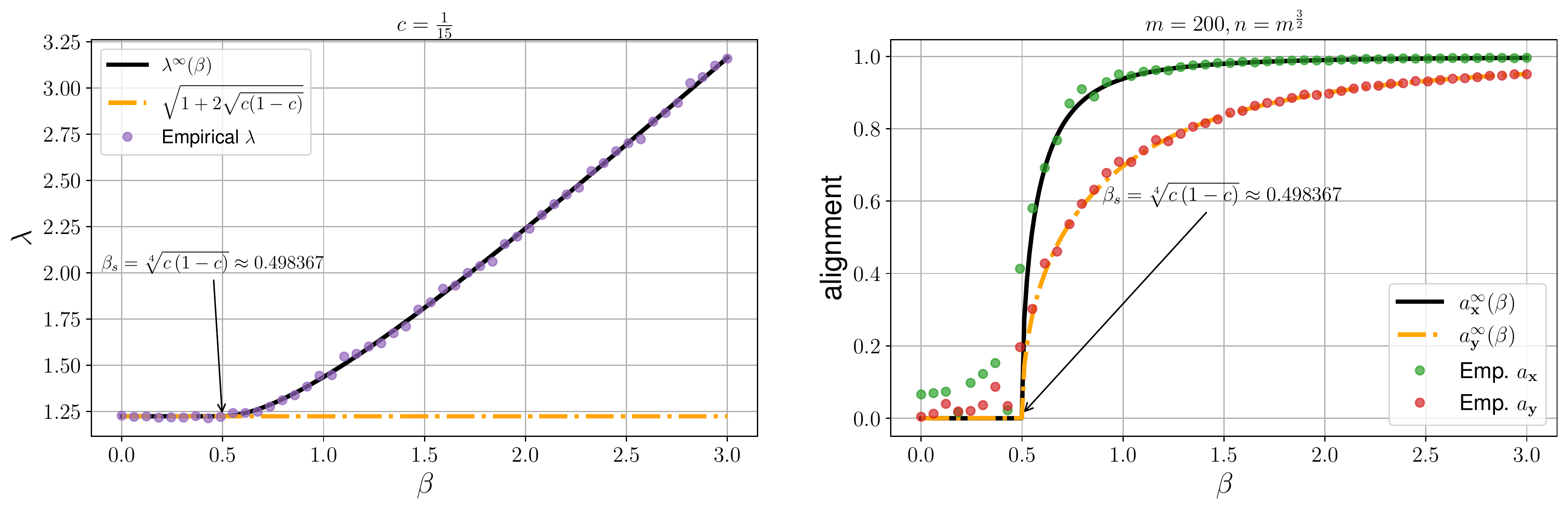}
    \caption{Asymptotic singular value and alignments of the spiked random matrix in \eqref{eq:spiked_random_matrix_simus} as per Corollary \ref{cor:spike_align_matrices}, and their simulated counterparts in the long matrix regime ($m=200, n=m^{\frac32}$) which simulated tensor unfolding as in Remark \ref{remark:unfolding}. Simulations are performed through SVD applied to $\mM$.}
    \label{fig:long_matrix}
\end{figure}

We start by considering the spiked random matrix model of the form
\begin{align}\label{eq:spiked_random_matrix_simus}
\mM = \beta \vx \vy^\top + \frac{1}{\sqrt{m+n}} \mX, \quad \text{with}\quad \mX \sim \sM_{m,n}(\mathcal{N}(0, 1)),
\end{align}
and $\vx, \vy$ are unitary vectors of dimensions $m$ and $n$ respectively.
Figure \ref{fig:spectrum_matrix} depicts the spectrum of $\begin{bmatrix}
\vzero_{m\times m} & \mM\\
\mM^\top & \vzero_{n\times n} 
\end{bmatrix}$ for $\beta = 0$ for different values of the dimensions $m,n$, and the predicted limiting spectral measure as per Corollary \ref{cor:spectrum_matrices}.
Figure \ref{fig:simus_matrix} shows a comparison between the asymptotic singular value and alignments obtained in Corollary \ref{cor:spike_align_matrices} and their simulated counterparts through SVD applied on $\mM$. Figure \ref{fig:long_matrix} further provides comparison between theory and simulations in the case of long random matrices ($m=200, n = m^{\frac32}$), which corresponds to tensor unfolding as per Remark \ref{remark:unfolding}, where a perfect matching is also observed between theory and simulations.

\begin{figure}[t!]
    \centering
    \includegraphics[width=14cm]{figs/simus_matrix.pdf}
    \caption{Asymptotic singular value and alignments of the spiked random matrix in \eqref{eq:spiked_random_matrix_simus} as per Corollary \ref{cor:spike_align_matrices}, and their simulated counterparts for different dimensions $m,n$. Simulations are performed through SVD applied to $\mM$.}
    \label{fig:simus_matrix}
\end{figure}

\subsection{Order $3$ tensors}\label{sec:simus_order_3}
Now we consider a $3$-order random tensor model of the form
\begin{align}\label{eq:spiked_tensor_model_simus}
\tT = \beta \vx \otimes \vy \otimes \vz + \frac{1}{\sqrt{m+n+p}} \tX, \quad \text{with}\quad \tX\sim\sT_{m,n,p} (\mathcal{N}(0, 1)),
\end{align}
and $\vx, \vy, \vz$ are unitary vectors of dimensions $m,n$ and $p$ respectively. In our simulations the estimation of the singular vectors $\vu, \vv, \vw$ of $\tT$ is performed using the power iteration method described by Algorithm \ref{alg:power_method}. We consider three initialization strategies for Algorithm \ref{alg:power_method}:
\begin{itemize}
\item[(i)] Random initialization by randomly sampling $\vu_0,\vv_0,\vw_0$ from the unitary spheres $\sS^{m-1}, \sS^{n-1}$ and $\sS^{p-1}$ respectively. We refer to this strategy in the figures legends by ``Random init.''.
\item[(ii)] Initialization with the true components $\vx,\vy,\vz$. We refer to this strategy in the figures legends by ``Init. with $\vx, \vy, \vz$''.
\item[(iii)] We start by running Algorithm \ref{alg:power_method} with strategy (i) for a high value of SNR $\beta \gg 1$ then progressively diminishing $\beta$ and initializing Algorithm \ref{alg:power_method} with the components obtained for the precedent value of $\beta$. We refer to this strategy in the figures legends by ``Init. with $\vu, \vv, \vw$''.
\end{itemize}

\begin{algorithm}
\caption{Tensor power method \cite{anandkumar2014tensor}}\label{alg:power_method}
\begin{algorithmic}
\Require An order $d$ tensor $\tT \in \sT_{n_1,\ldots, n_d}$ and initial components $(\vu^{(1)}_0,\ldots, \vu^{(d)}_0) \in \sS^{n_1 - 1} \times \cdots \times \sS^{n_d - 1}$
\Ensure Estimate of $ \argmin_{\lambda\,\in\, \sR,\, (\vu^{(1)},\ldots,\vu^{(d)})\,\in\, \sS^{n_1 - 1} \times \cdots \times \sS^{n_d - 1}} \Vert \tT - \lambda \, \vu^{(1)} \otimes  \cdots \otimes \vu^{(d)} \Vert_{\text{F}}^2 $\\
\textbf{Repeat}
\For{$i\in[d]$}
	\State $\vu^{(i)} \gets \frac{\tT(\vu^{(1)}, \ldots, \vu^{(i-1)}, :, \vu^{(i+1)}, \ldots, \vu^{(d)})}{\Vert \tT(\vu^{(1)}, \ldots, \vu^{(i-1)}, :, \vu^{(i+1)}, \ldots, \vu^{(d)}) \Vert} $ \Comment{Contract $\tT$ on all the $\vu^{(j)}$'s for $j\neq i$}
\EndFor\\
\textbf{until convergence of the $\vu^{(i)}$'s}
\end{algorithmic}
\end{algorithm}

\begin{figure}[h!]
    \centering
    \includegraphics[width=14cm]{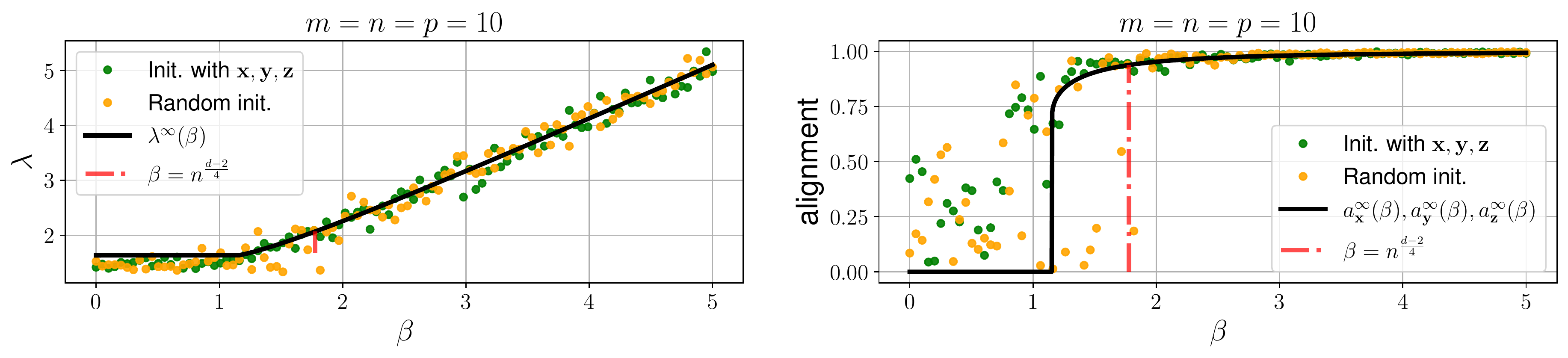}
    \includegraphics[width=14cm]{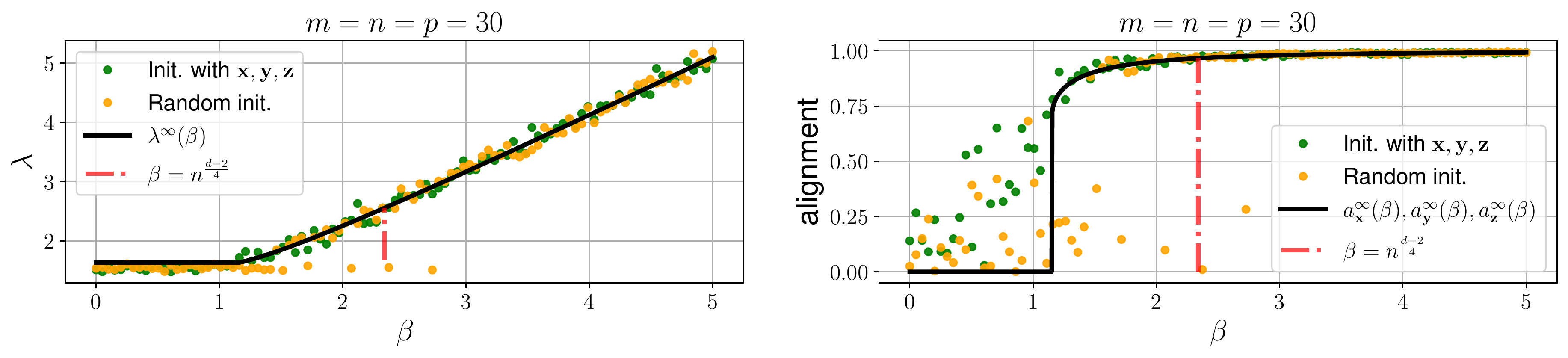}
    \includegraphics[width=14cm]{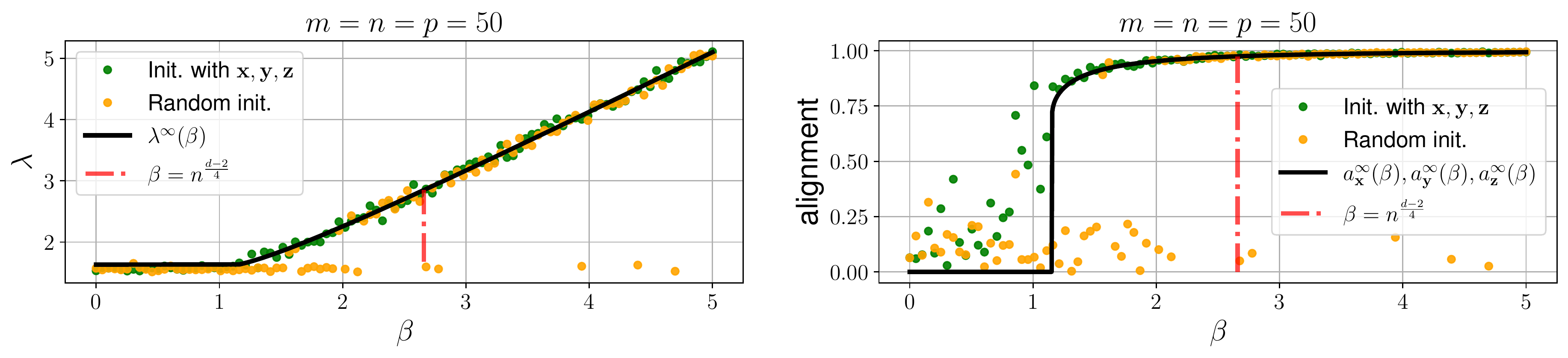}
    \includegraphics[width=14cm]{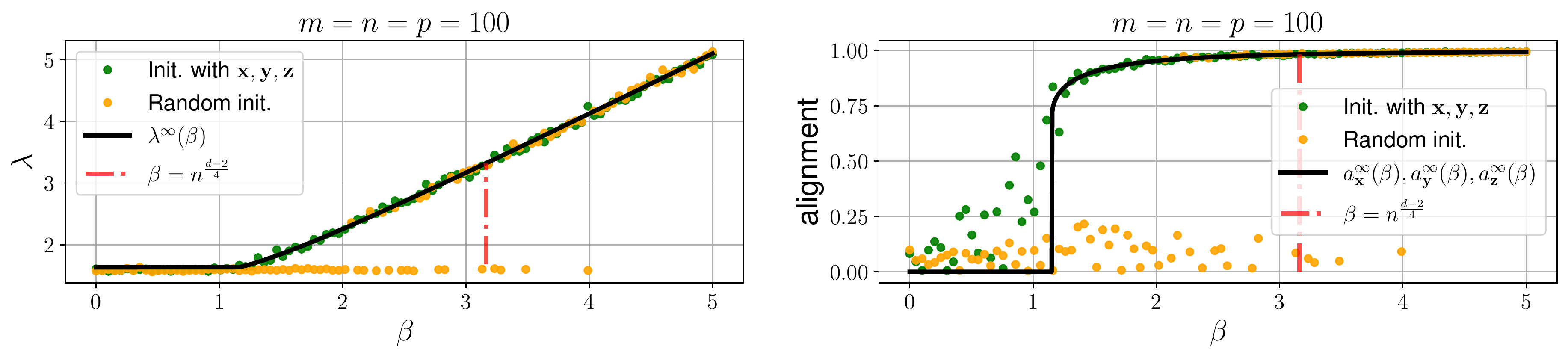}
    \includegraphics[width=14cm]{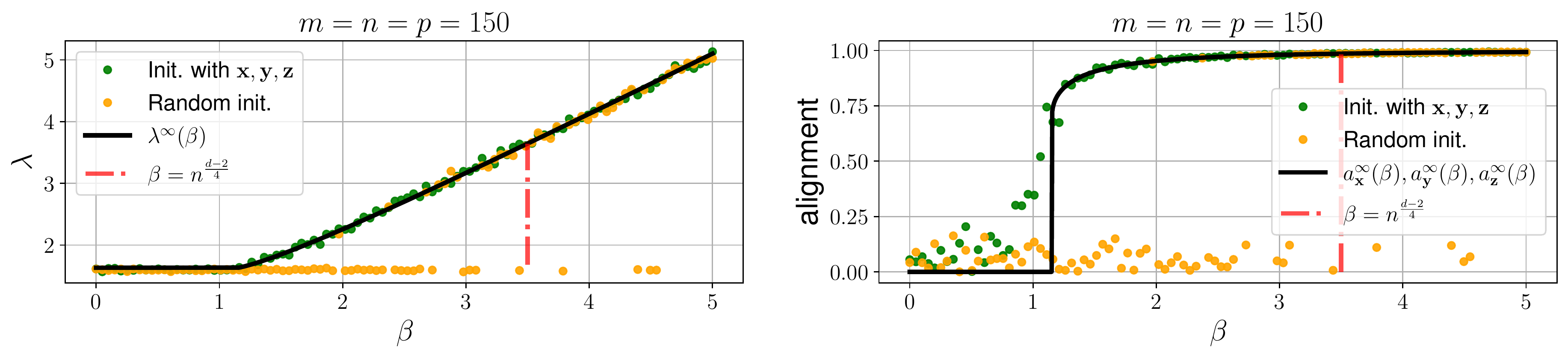}
    \caption{Asymptotic singular value and alignments of the order $3$ tensor in \eqref{eq:spiked_tensor_model_simus} having equal dimensions as per Corollary \ref{cor:cubic}, and their simulated counterparts for different dimensions $\{10, 30, 50, 100, 150\}$. The yellow dots correspond to a random initialization of the power iteration method strategy (i), while the green dots correspond to an initialization with the true spike components $\vx, \vy$ and $\vz$ strategy (ii).}
    \label{fig:simus_order_3}
\end{figure}

\begin{figure}[h!]
    \centering
    \includegraphics[width=14cm]{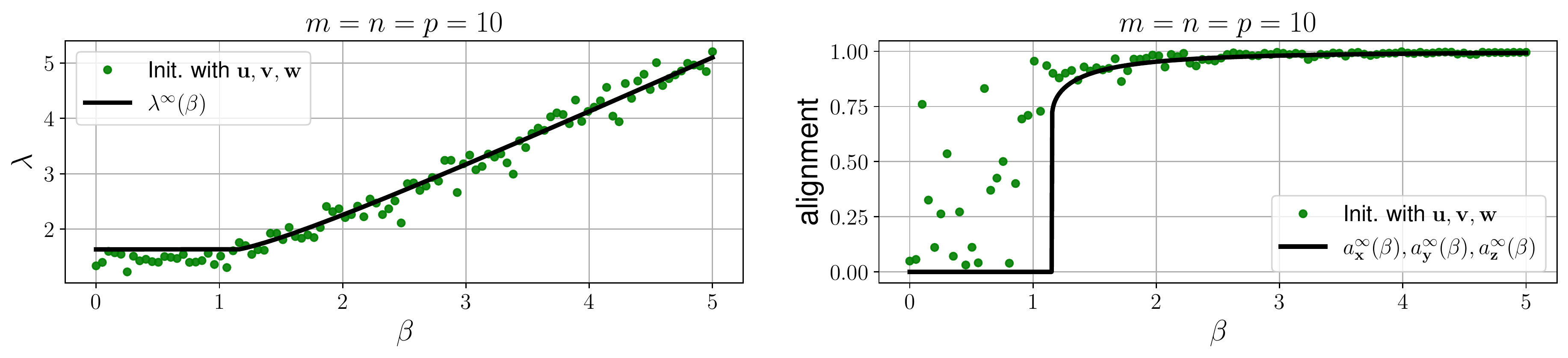}
    \includegraphics[width=14cm]{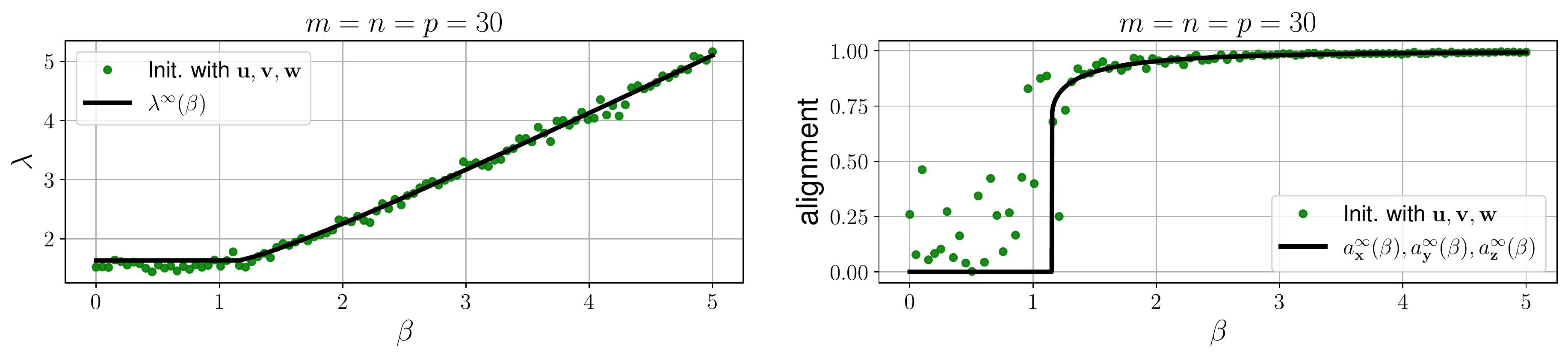}
    \includegraphics[width=14cm]{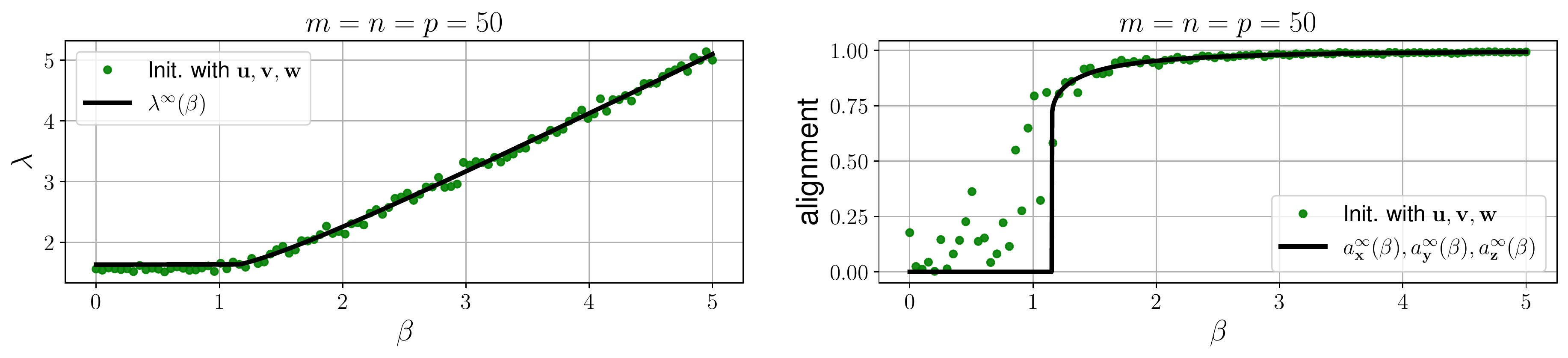}
    \includegraphics[width=14cm]{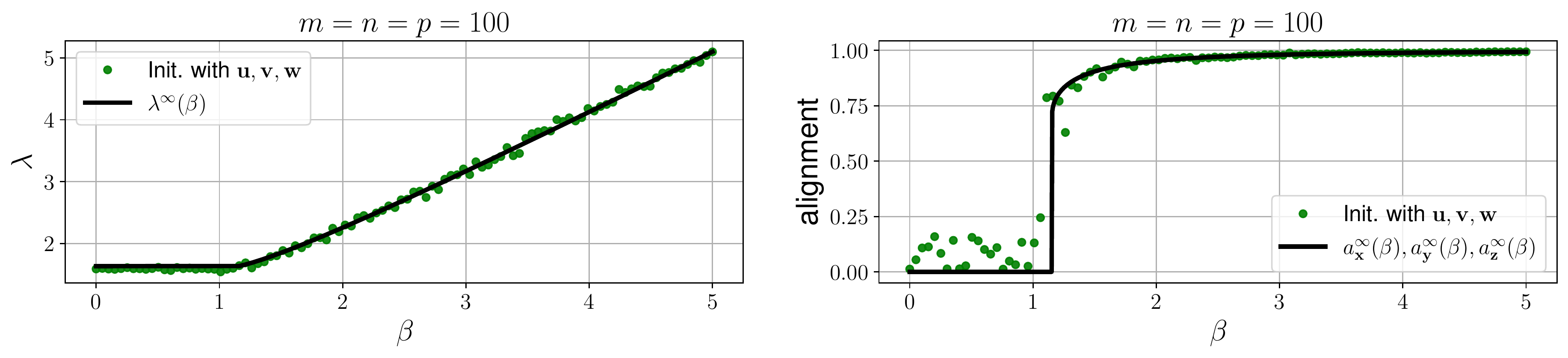}
    \includegraphics[width=14cm]{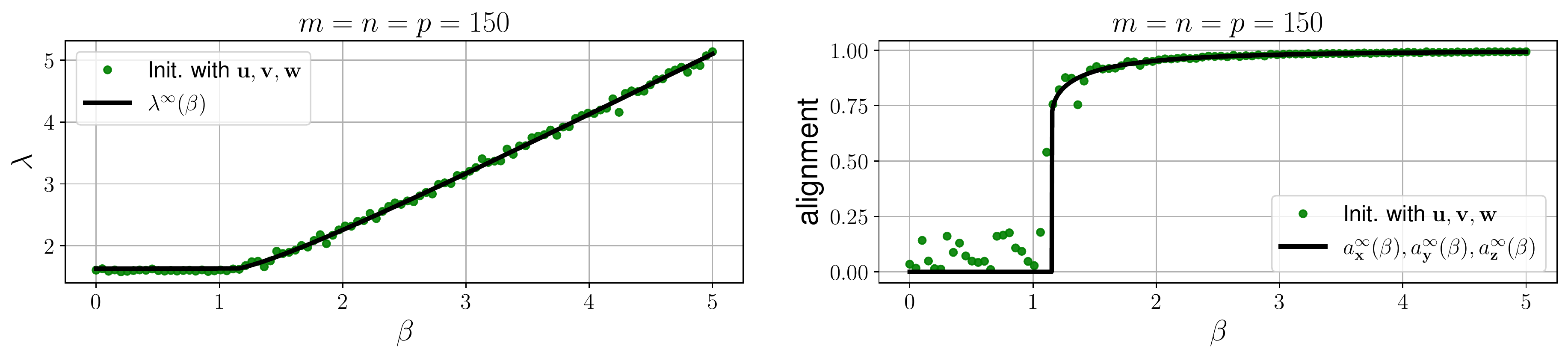}
    \caption{Asymptotic singular value and alignments of the order $3$ tensor in \eqref{eq:spiked_tensor_model_simus} having equal dimensions as per Corollary \ref{cor:cubic}, and their simulated counterparts for different dimensions $\{10, 30, 50, 100, 150\}$. The green dots correspond the initialization with strategy (iii).}
    \label{fig:simus_order_3_bis}
\end{figure}

Figure \ref{fig:simus_order_3} provides a comparison between the asymptotic singular value and alignments of $\tT$ (obtained by Corollary \ref{cor:cubic}) with their simulation counterparts where the singular vectors are estimated by Algorithm \ref{alg:power_method} with the initialization strategies (i) in yellow dots and (ii) in green dots. As we can see, as the tensor dimensions grow, the numerical estimates approach their asymptotic counterparts for (ii). Besides, the random initialization (i) yields poor convergence for $\beta$ around its minimal value $\beta_s$ when the tensor dimension is large enough (see $m=n=p=150$). This phenomenon is related to the algorithmic time complexity of the power method which is known to succeed in recovering the underlying signal in polynomial time provided that $\beta \gtrsim n^{\frac{d-2}{4}}$ \cite{montanari2014statistical, biroli2020iron}, which is also noticeable in our simulations (the algorithmic phase transition seems to grow with the tensor dimensions). Figure \ref{fig:simus_order_3_bis} further depicts comparison between theory and simulations when Algorithm \ref{alg:power_method} is initialized following strategy (iii) which allows to follow the trajectory of the global maximum of the ML problem, where we can notice a good matching between the asymptotic curves and their simulated counterparts.   

\section{Proofs}\label{sec:proofs}
\subsection{Derivative of tensor singular value and vectors}\label{proof:derivs}
Deriving the identities in \eqref{eq:singular_vectors} and \eqref{eq:singular_value} w.r.t.\@ the entry $X_{ijk}$ of the tensor noise $\tX$, we obtain the following set of equations
\begin{align*}
\begin{cases}
\tT(\vw) \frac{\partial \vv}{\partial X_{ijk}} + \tT(\vv) \frac{\partial \vw}{\partial X_{ijk}} + \frac{1}{\sqrt N} v_j w_k \ve_i^m = \frac{\partial \lambda}{\partial X_{ijk}} \vu + \lambda \frac{\partial \vu}{\partial X_{ijk}},\\
\tT(\vw)^\top \frac{\partial \vu}{\partial X_{ijk}} + \tT(\vu) \frac{\partial \vw}{\partial X_{ijk}} + \frac{1}{\sqrt N} u_i w_k \ve_j^n = \frac{\partial \lambda}{\partial X_{ijk}} \vv + \lambda \frac{\partial \vv}{\partial X_{ijk}},\\
\tT(\vv)^\top \frac{\partial \vu}{\partial X_{ijk}}  + \tT(\vu)^\top \frac{\partial \vv}{\partial X_{ijk}} + \frac{1}{\sqrt N} u_i v_j \ve_k^p = \frac{\partial \lambda}{\partial X_{ijk}} \vw + \lambda \frac{\partial \vw}{\partial X_{ijk}},\\
\frac{\partial \lambda}{\partial X_{ijk}} = \tT\left( \frac{\partial \vu}{\partial X_{ijk}}, \vv, \vw \right) + \tT\left(\vu,  \frac{\partial \vv}{\partial X_{ijk}}, \vw \right) + \tT\left(\vu, \vv,  \frac{\partial \vw}{\partial X_{ijk}} \right) + \frac{1}{\sqrt N} u_i v_j w_k.
\end{cases}
\end{align*}
Writing $\tT\left( \frac{\partial \vu}{\partial X_{ijk}}, \vv, \vw \right)$ as $\tT\left( \frac{\partial \vu}{\partial X_{ijk}}, \vv, \vw \right) = \left(\frac{\partial \vu}{\partial X_{ijk}}\right)^\top \tT(\vv)\vw $, we can apply again the identities in \eqref{eq:singular_vectors} which results in $\tT\left( \frac{\partial \vu}{\partial X_{ijk}}, \vv, \vw \right) = \lambda  \left(\frac{\partial \vu}{\partial X_{ijk}}\right)^\top \vu$. Doing similarly with $\tT\left(\vu,  \frac{\partial \vv}{\partial X_{ijk}}, \vw \right)$ and $\tT\left(\vu, \vv,  \frac{\partial \vw}{\partial X_{ijk}} \right)$, we have
\begin{align*}
\frac{\partial \lambda}{\partial X_{ijk}} = \lambda \left(  \left(\frac{\partial \vu}{\partial X_{ijk}}\right)^\top \vu +  \left(\frac{\partial \vv}{\partial X_{ijk}}\right)^\top \vv +  \left(\frac{\partial \vw}{\partial X_{ijk}}\right)^\top \vw  \right) + \frac{1}{\sqrt N} u_i v_j w_k.
\end{align*}
Furthermore, since $\vu^\top \vu = \vv^\top \vv = \vw^\top \vw = 1$, we have 
\begin{align*}
\left(\frac{\partial \vu}{\partial X_{ijk}}\right)^\top \vu = \left(\frac{\partial \vv}{\partial X_{ijk}}\right)^\top \vv = \left(\frac{\partial \vw}{\partial X_{ijk}}\right)^\top \vw = 0.
\end{align*}
Thus the derivative of $\lambda$ writes simply as
\begin{align*}
\frac{\partial \lambda}{\partial X_{ijk}} = \frac{1}{\sqrt N} u_i v_j w_k.
\end{align*}
Hence, we find that
\begin{align*}
\lambda \begin{bmatrix}
\frac{\partial \vu}{\partial X_{ijk}} \\
\frac{\partial \vv}{\partial X_{ijk}} \\
\frac{\partial \vw}{\partial X_{ijk}}
\end{bmatrix} = \frac{1}{\sqrt N} \begin{bmatrix}
v_j w_k (\ve_i^m - u_i \vu) \\
u_i w_k (\ve_j^n - v_j \vv) \\
u_i v_j (\ve_k^p - w_k \vw)
\end{bmatrix} + \Phi_3(\tT, \vu, \vv, \vw) \begin{bmatrix}
\frac{\partial \vu}{\partial X_{ijk}} \\
\frac{\partial \vv}{\partial X_{ijk}} \\
\frac{\partial \vw}{\partial X_{ijk}}
\end{bmatrix}.
\end{align*}
Yielding the expression in \eqref{eq:deriv}. The same calculations apply to the more general $d$-order tensor case yielding the identity in \eqref{eq:deriv_order_d}.

\subsection{Proof of Theorem \ref{thm:semi-circle_independent}}\label{proof:semi-circle_independent}
Denote the matrix model as
\begin{align*}
\mN \equiv \frac{1}{\sqrt N} \Phi_3 (\tX, \va, \vb, \vc),
\end{align*}
where we recall $\tX\sim \sT_{m,n,p}(\mathcal{N}(0, 1))$ and $(\va, \vb, \vc)\in \sS^{m-1}\times \sS^{m-1}\times \sS^{p-1}$ are independent of $\tX$. We further denote the resolvent matrix of $\mN$ as 
\begin{align*}
\mQ(z) \equiv \left( \mN - z \mI_N \right)^{-1} = \begin{bmatrix}
\mQ^{11}(z) & \mQ^{12}(z) & \mQ^{13}(z)\\
\mQ^{12}(z)^\top & \mQ^{22}(z) & \mQ^{23}(z) \\
\mQ^{13}(z)^\top & \mQ^{23}(z)^\top & \mQ^{33}(z)
\end{bmatrix}.
\end{align*}
In order to characterize the limiting Stieltjes transform $g(z)$ of $\mN$, we need to estimate the quantity $\frac1N \tr \mQ(z) \asto g(z)$ (as a consequence of Theorem \ref{thm:continuity}). We further introduce the following limits
\begin{align*}
\frac1N \tr \mQ^{11}(z) \asto g_1(z), \quad \frac1N \tr \mQ^{22}(z) \asto g_2(z), \quad \frac1N \tr \mQ^{33}(z) \asto g_3(z).
\end{align*}
From the identity in \eqref{eq:resolvent_identity}, we have
\begin{align*}
\mN \mQ(z) - z \mQ(z) = \mI_N,
\end{align*}
from which we particularly have
\begin{align*}
\frac{1}{\sqrt N} \left[ \tX(\vc) \mQ^{12}(z)^\top \right]_{ii} + \frac{1}{\sqrt N} \left[ \tX(\vb) \mQ^{13}(z)^\top \right]_{ii} - z Q^{11}_{ii}(z) = 1,
\end{align*} 
or 
\begin{align}\label{eq:traceQ11}
\frac{1}{N \sqrt N} \sum_{i=1}^m \left[ \tX(\vc) \mQ^{12}(z)^\top \right]_{ii} + \frac{1}{N \sqrt N} \sum_{i=1}^m \left[ \tX(\vb) \mQ^{13}(z)^\top \right]_{ii} - \frac{z}{N} \tr \mQ^{11}(z) = \frac{m}{N}.
\end{align} 
We thus need to compute the expectations of $\frac{1}{N \sqrt N} \sum_{i=1}^m \left[ \tX(\vc) \mQ^{12}(z)^\top \right]_{ii}$ and $ \frac{1}{N \sqrt N} \sum_{i=1}^m \left[ \tX(\vb) \mQ^{13}(z)^\top \right]_{ii}$, which develop as 
\begin{align*}
 \frac{1}{N \sqrt N} \sum_{i=1}^m \EE \left[ \tX(\vc) \mQ^{12}(z)^\top \right]_{ii} = \frac{1}{N\sqrt N} \sum_{i=1}^m \sum_{j=1}^n \sum_{k=1}^p c_k \EE \left[ X_{ijk} Q^{12}_{ij} \right] = \frac{1}{N\sqrt N} \sum_{i=1}^m \sum_{j=1}^n \sum_{k=1}^p c_k \EE \left[ \frac{\partial Q^{12}_{ij}}{\partial X_{ijk}} \right],
\end{align*} 
where the last equality is obtained by applying Stein's lemma (Lemma \ref{lemma:stein}). For continuing the derivations, we need to express the derivative of the resolvent $\mQ(z)$ with respect to an entry $X_{ijk}$ of the tensor noise $\tX$. Indeed, since $\mN \mQ(z) - z \mQ(z) = \mI_N$, we have
\begin{align*}
\frac{\partial \mN}{\partial X_{ijk}} \mQ(z) + \mN \frac{\partial \mQ(z)}{\partial X_{ijk}} - z \frac{\partial \mQ(z)}{\partial X_{ijk}} = \vzero_{N\times N} \quad \Rightarrow\quad  (\mN - z\mI_N) \frac{\partial \mQ(z)}{\partial X_{ijk}} = - \frac{\partial \mN}{\partial X_{ijk}} \mQ(z),
\end{align*}
from which we get
\begin{align*}
\frac{\partial \mQ(z)}{\partial X_{ijk}} = - \mQ(z)  \frac{\partial \mN}{\partial X_{ijk}} \mQ(z),
\end{align*}
where
\begin{align*}
\frac{\partial \mN}{\partial X_{ijk}} = \frac{1}{\sqrt N} \begin{bmatrix}
\vzero_{m\times m} & c_k \ve_i^m (\ve_j^n)^\top & b_j \ve_i^m (\ve_k^p)^\top \\
c_k \ve_j^n (\ve_i^m)^\top & \vzero_{n\times n} & a_i \ve_j^n (\ve_k^p)^\top \\
b_j \ve_k^p (\ve_i^m )^\top & a_i \ve_k^p (\ve_j^n )^\top & \vzero_{p\times p}
\end{bmatrix},
\end{align*}
and we finally obtain the following derivatives
\begin{align*}
\frac{\partial Q_{ab}^{11}}{\partial X_{ijk}} &= - \frac{1}{\sqrt N} \left[ a_i (Q_{aj}^{12} Q_{bk}^{13} + Q_{ak}^{13} Q_{bj}^{12}) + b_j (Q_{ai}^{11} Q_{bk}^{13} + Q_{ak}^{13} Q_{ib}^{11}) + c_k (Q_{ai}^{11} Q_{bj}^{12} + Q_{aj}^{12} Q_{ib}^{11}) \right],\\
\frac{\partial Q_{ab}^{12}}{\partial X_{ijk}} &= - \frac{1}{\sqrt N} \left[ a_i (Q_{aj}^{12} Q_{bk}^{23} + Q_{ak}^{13} Q_{jb}^{22}) + b_j (Q_{ai}^{11} Q_{bk}^{23} + Q_{ak}^{13} Q_{ib}^{12}) + c_k (Q_{ai}^{11} Q_{jb}^{22} + Q_{aj}^{12} Q_{ib}^{12}) \right],\\
\frac{\partial Q_{ab}^{13}}{\partial X_{ijk}} &= - \frac{1}{\sqrt N} \left[ a_i (Q_{aj}^{12} Q_{kb}^{33} + Q_{ak}^{13} Q_{jb}^{23}) + b_j (Q_{ai}^{11} Q_{kb}^{33} + Q_{ak}^{13} Q_{ib}^{13}) + c_k (Q_{ai}^{11} Q_{jb}^{23} + Q_{aj}^{12} Q_{ib}^{13}) \right],\\
\frac{\partial Q_{ab}^{22}}{\partial X_{ijk}} &= - \frac{1}{\sqrt N} \left[ a_i (Q_{aj}^{22} Q_{bk}^{23} + Q_{ak}^{23} Q_{jb}^{22}) + b_j (Q_{ia}^{12} Q_{bk}^{23} + Q_{ak}^{23} Q_{ib}^{12}) + c_k (Q_{ia}^{12} Q_{jb}^{22} + Q_{aj}^{22} Q_{ib}^{12}) \right],\\
\frac{\partial Q_{ab}^{23}}{\partial X_{ijk}} &= - \frac{1}{\sqrt N} \left[ a_i (Q_{aj}^{22} Q_{kb}^{33} + Q_{ak}^{23} Q_{jb}^{23}) + b_j (Q_{ia}^{12} Q_{kb}^{33} + Q_{ak}^{23} Q_{ib}^{13}) + c_k (Q_{ia}^{12} Q_{jb}^{23} + Q_{aj}^{22} Q_{ib}^{13}) \right],\\
\frac{\partial Q_{ab}^{33}}{\partial X_{ijk}} &= - \frac{1}{\sqrt N} \left[ a_i (Q_{ja}^{23} Q_{kb}^{33} + Q_{ak}^{33} Q_{jb}^{23}) + b_j (Q_{ia}^{13} Q_{kb}^{33} + Q_{ak}^{33} Q_{ib}^{13}) + c_k (Q_{ia}^{13} Q_{jb}^{23} + Q_{ja}^{23} Q_{ib}^{13}) \right].
\end{align*}
In particular, 
\begin{align*}
\frac{\partial Q_{ij}^{12}}{\partial X_{ijk}} &= - \frac{1}{\sqrt N} \left[ a_i (Q_{ij}^{12} Q_{jk}^{23} + Q_{ik}^{13} Q_{jj}^{22}) + b_j (Q_{ii}^{11} Q_{jk}^{23} + Q_{ik}^{13} Q_{ij}^{12}) + c_k (Q_{ii}^{11} Q_{jj}^{22} + Q_{ij}^{12} Q_{ij}^{12}) \right].
\end{align*}
Going back to the computation of $A \equiv \frac{1}{N\sqrt N} \sum_{i=1}^m \sum_{j=1}^n \sum_{k=1}^p c_k \EE \left[ \frac{\partial Q^{12}_{ij}}{\partial X_{ijk}} \right]$, we will see that the only contributing term in the derivative $\frac{\partial Q_{ij}^{12}}{\partial X_{ijk}}$ is $- \frac{1}{\sqrt N} c_k Q_{ii}^{11} Q_{jj}^{22}$. Indeed,
\begin{align*}
&A = -\frac{1}{N^2} \sum_{ijk} c_k \EE \left[ a_i (Q_{ij}^{12} Q_{jk}^{23} + Q_{ik}^{13} Q_{jj}^{22}) + b_j (Q_{ii}^{11} Q_{jk}^{23} + Q_{ik}^{13} Q_{ij}^{12}) + c_k (Q_{ii}^{11} Q_{jj}^{22} + Q_{ij}^{12} Q_{ij}^{12}) \right],\\
&= -\frac{1}{N^2} \EE \left[ \va^\top \mQ^{12} \mQ^{23} \vc + \va^\top \mQ^{13}\vc \tr\mQ^{22} + \tr \mQ^{11} \vb^\top \mQ^{23}\vc + \vb^\top (\mQ^{12})^\top \mQ^{13}\vc + \tr \mQ^{11}\tr\mQ^{22} + \tr(\mQ^{12}(\mQ^{12})^\top) \right].
\end{align*}
Now, since the vectors $\va, \vb, \vc$ are of bounded norms and assuming $\mQ(z)$ is of bounded spectral norm (see condition in \eqref{eq:resovent_bounded}), under Assumption \ref{ass:growth} (as $N\to \infty$), the terms $\frac{1}{N^2} \va^\top \mQ^{12} \mQ^{23} \vc$, $\frac{1}{N^2}\va^\top \mQ^{13}\vc \tr\mQ^{22}$, $\frac{1}{N^2} \tr \mQ^{11} \vb^\top \mQ^{23}\vc$, $\frac{1}{N^2} \vb^\top (\mQ^{12})^\top \mQ^{13}\vc $ and $\frac{1}{N^2} \tr(\mQ^{12}(\mQ^{12})^\top) $ are vanishing almost surely. As such, we find that
\begin{align*}
\frac{1}{N \sqrt N} \sum_{i=1}^m \left[ \tX(\vc) \mQ^{12}(z)^\top \right]_{ii}= -\frac1N \tr \mQ^{11}(z)\frac1N \tr\mQ^{22}(z) + \mathcal{O}(N^{-1})\asto -g_1(z)g_2(z) + \mathcal{O}(N^{-1}).
\end{align*}
Similarly, we find that
\begin{align*}
\frac{1}{N \sqrt N} \sum_{i=1}^m \left[ \tX(\vb) \mQ^{13}(z)^\top \right]_{ii} = -\frac1N \tr \mQ^{11}(z)\frac1N \tr\mQ^{33}(z) + \mathcal{O}(N^{-1})\asto -g_1(z)g_3(z) + \mathcal{O}(N^{-1}).
\end{align*}
From \eqref{eq:traceQ11}, $g_1(z)$ satisfies
\begin{align*}
-g_1(z) \left( g_2(z) + g_3(z) \right) - z g_1(z) = c_1,
\end{align*}
where we recall $c_1 = \lim \frac{m}{N}$. Similarly, $g_2(z)$ and $g_3(z)$ satisfy
\begin{align*}
\begin{cases}
-g_2(z) (g_1(z) + g_3(z)) - z g_2(z) = c_2,\\
-g_3(z) (g_1(z) + g_2(z)) - z g_3(z) = c_3,
\end{cases}
\end{align*}
where we recall again $c_2 = \lim \frac{n}{N}$ and $c_3 = \lim \frac{p}{N}$. Moreover, by definition, $g(z) = \sum_{i=1}^3 g_i(z)$, thus we have for each $i\in[3]$
\begin{align*}
g_i(z)(g(z) - g_i(z)) + z g_i(z) + c_i = 0,
\end{align*}
yielding
\begin{align*}
g_i(z) = \frac{g(z) + z}{2} - \frac{\sqrt{4c_i + (g(z)+z)^2}}{2},
\end{align*}
with $g(z)$ solution of the equation $g(z) = \sum_{i=1}^3 g_i(z)$ satisfying $\Im[g(z)]>0$ for $z\in \sC$ with $\Im[z]>0$ (see Property 2 of the Stieltjes transform in Subsection \ref{sec:stieltjes_transform}).

\subsection{Proof of Corollary \ref{cor:semi-circle-cubic}}\label{proof:semi-circle-cubic}
Given the result of Theorem \ref{thm:semi-circle_independent}, setting $c_1=c_2=c_3=\frac13$, we have for all $i\in [3]$
\begin{align*}
g_i(z) = \frac{g(z) + z}{2} - \frac{\sqrt{\frac43 + (g(z) + z)^2}}{2},
\end{align*}
where $g(z)$ satisfies $g(z) = \sum_{i=1}^3 g_i(z)$, thus $g(z)$ is the solution to
\begin{align*}
z + \frac{g(z) + z}{2} - \frac{3 \sqrt{ \frac43 + (g(z) + z)^2 }}{2} = 0,
\end{align*}
solving in $g(z)$ yields
\begin{align*}
g(z) \in \left\{- \frac{3 z}{4} - \frac{\sqrt{3} \sqrt{3 z^{2} - 8}}{4}, - \frac{3 z}{4} + \frac{\sqrt{3} \sqrt{3 z^{2} - 8}}{4}\right\},
\end{align*}
and the limiting Stieltjes transform (with $\Im [g(z)]>0$ for $z$ with $\Im(z)>0$ by property 2 of the Stieltjes transform, see Subsection \ref{sec:stieltjes_transform}) is therefore
\begin{align*}
g(z) = - \frac{3 z}{4} + \frac{\sqrt{3} \sqrt{3 z^{2} - 8}}{4}.
\end{align*}
\subsection{Proof of Theorem \ref{thm:semi-circle_dependent}}\label{proof:semi-circle_dependent}
Given the random tensor model in \eqref{eq:spiked_tensor_model} and its singular vectors characterized by \eqref{eq:singular_vectors}, we denote the associated random matrix model as
\begin{align*}
\mT \equiv \Phi_3(\tT, \vu, \vv, \vw) = \beta \mV \mB \mV^\top + \mN,
\end{align*}
where
\begin{align*}
     \mN = \frac{1}{\sqrt N} \Phi_3(\tX, \vu, \vv, \vw),\quad \mB \equiv \begin{bmatrix}
 0 & \langle \vz, \vw\rangle & \langle \vy, \vv\rangle \\
 \langle \vz, \vw\rangle & 0 & \langle \vx, \vu\rangle \\
 \langle \vy, \vv\rangle & \langle \vx, \vu\rangle & 0
\end{bmatrix}\in \sM_3, \quad \mV \equiv \begin{bmatrix}
 \vx & \vzero_m & \vzero_m \\
 \vzero_n & \vy & \vzero_n \\
 \vzero_p & \vzero_p & \vz
\end{bmatrix}\in \sM_{N,3}.
\end{align*}
We further denote the resolvents of $\mT$ and $\mN$ respectively as
\begin{align*}
\mR(z)& = \left( \mT - z \mI_N \right)^{-1} =  \begin{bmatrix}
\mR^{11}(z) & \mR^{12}(z) & \mR^{13}(z)\\
\mR^{12}(z)^\top & \mR^{22}(z) & \mR^{23}(z) \\
\mR^{13}(z)^\top & \mR^{23}(z)^\top & \mR^{33}(z)
\end{bmatrix},\\ \mQ(z) &= \left( \mN - z \mI_N \right)^{-1} =  \begin{bmatrix}
\mQ^{11}(z) & \mQ^{12}(z) & \mQ^{13}(z)\\
\mQ^{12}(z)^\top & \mQ^{22}(z) & \mQ^{23}(z) \\
\mQ^{13}(z)^\top & \mQ^{23}(z)^\top & \mQ^{33}(z)
\end{bmatrix}.
\end{align*}
By Woodbury matrix identity (Lemma \ref{lem:woodbury}), we have
\begin{align}\label{eq:R2Q}
\mR(z) = \mQ(z) - \mQ(z) \mV \left( \frac{1}{\beta}\mB^{-1} + \mV^\top \mQ(z) \mV \right)^{-1} \mV^\top \mQ(z), 
\end{align}
In particular, taking the normalized trace operator, we get
\begin{align*}
\frac1N \tr \mR(z) &= \frac1N \tr \mQ(z) - \frac1N \tr \left[ \left( \frac{1}{\beta}\mB^{-1} + \mV^\top \mQ(z) \mV \right)^{-1} \mV^\top \mQ^2(z) \mV  \right],\\
&= \frac1N \tr \mQ(z) + \mathcal{O}(N^{-1}),
\end{align*}
since the matrix $ \left( \frac{1}{\beta}\mB^{-1} + \mV^\top \mQ(z) \mV \right)^{-1} \mV^\top \mQ^2(z) \mV$ is of bounded spectral norm (assuming $\Vert \mQ(z) \Vert$ is bounded, see condition in \eqref{eq:resovent_bounded}) and being of finite size ($3\times 3$ matrix). As such, the asymptotic spectral measure of $\mT$ is the same as the one of $\mN$ which can be estimated though $\frac1N \tr \mQ(z)$. Comparing to the result from Appendix \ref{proof:semi-circle_independent}, now the singular vectors $\vu,\vv, \vw$ depend statistically on the tensor noise $\tX$ which needs to be handled. 

From the identity \eqref{eq:resolvent_identity}, we have $\mN \mQ(z) - z \mQ(z) = \mI_N$, from which
\begin{align}\label{eq:traceQ11_dependent}
\frac{1}{N \sqrt N} \sum_{i=1}^m \left[ \tX(\vw) \mQ^{12}(z)^\top \right]_{ii} + \frac{1}{N \sqrt N} \sum_{i=1}^m \left[ \tX(\vv) \mQ^{13}(z)^\top \right]_{ii} - \frac{z}{N} \tr Q^{11}(z) = \frac{m}{N}.
\end{align} 
We thus need to compute the expectations of $\frac{1}{N \sqrt N} \sum_{i=1}^m \left[ \tX(\vw) \mQ^{12}(z)^\top \right]_{ii}$ and $\frac{1}{N \sqrt N} \sum_{i=1}^m \left[ \tX(\vv) \mQ^{13}(z)^\top \right]_{ii} $. In particular, 
\begin{align*}
A\equiv \frac{1}{N \sqrt N} \sum_{i=1}^m \EE \left[ \tX(\vw) \mQ^{12}(z)^\top \right]_{ii} = \frac{1}{N\sqrt N}\sum_{ijk} \EE \left[ X_{ijk} w_k Q^{12}_{ij} \right] = \frac{1}{N\sqrt N}\sum_{ijk} \EE \left[ \frac{\partial (w_k Q^{12}_{ij})}{\partial X_{ijk}} \right],
\end{align*}
where the last equality is obtained by Stein's lemma (Lemma \ref{lemma:stein}). Due to the statistical dependency between $\vw$ and $\tX$, the above sum decomposes in two terms which are
\begin{align*}
A = \frac{1}{N\sqrt N}\sum_{ijk} \EE \left[ w_k \frac{\partial Q^{12}_{ij}}{\partial X_{ijk}} \right] + \frac{1}{N\sqrt N}\sum_{ijk} \EE \left[ \frac{\partial w_k}{\partial X_{ijk}} Q^{12}_{ij} \right] = A_1 + A_2,
\end{align*}
where the first term $A_1$ has already been handled in the previous subsection (if replacing $\vc$ with $\vw$). We now show that the second term $A_2$ is asymptotically vanishing under Assumption \ref{ass:growth}. Indeed, by \eqref{eq:deriv}, we have
\begin{align*}
\frac{\partial w_k}{\partial X_{ijk}} &= - \frac{1}{\sqrt N} \left( v_j w_k (R_{ik}^{13}(\lambda) - u_i \vu^\top\mR_{:,k}^{13}(\lambda)  ) \right) \\
 &- \frac{1}{\sqrt N} \left(  u_i w_k (R_{jk}^{23}(\lambda) - v_j \vv^\top \mR_{:,k}^{23}(\lambda) ) \right) - \frac{1}{\sqrt N} \left(  u_i v_j (R_{kk}^{33}(\lambda) - w_k \vw^\top \mR_{:,k}^{33}(\lambda) ) \right).
\end{align*}
As such $A_2$ decomposes in three terms $A_2 = A_{21} + A_{22} + A_{23}$, where
\begin{align*}
A_{21} &= - \frac{1}{N^2} \sum_{ijk} \EE \left[ v_j w_k R_{ik}^{13}(\lambda) Q_{ij}^{12} \right] + \frac{1}{N^2} \sum_{ijkl} \EE \left[ u_i v_j w_k u_l R_{lk}^{13}(\lambda) Q_{ij}^{12} \right],\\
&= -\frac{1}{N^2} \EE \left[ \vv^\top \mQ^{12}(z)^\top \mR^{13}(\lambda) \vw \right] + \frac{1}{N^2} \EE \left[ \vu^\top \mQ^{12}(z) \vv \vu^\top \mR^{13}(\lambda) \vw \right] \to 0,
\end{align*}
as $N\to \infty$, since the singular vectors $\vu, \vv,\vw$ are of bounded norms and assuming the resolvent $\mQ(z)$ and $\mR(\lambda)$ are of bounded spectral norms ($\mR(\lambda)$ has bounded spectral norm by Assumption \ref{ass:lambda_outside} and \eqref{eq:resovent_bounded}). Similarly, we further have
\begin{align*}
A_{22} &= -\frac{1}{N^2} \sum_{ijk} \EE \left[ u_i w_k R^{23}_{jk}(\lambda) Q_{ij}^{12} \right] + \frac{1}{N^2} \sum_{ijkl} \EE \left[ u_i v_j w_k v_l R^{23}_{lk}(\lambda) Q_{ij}^{12} \right],\\
&= - \frac{1}{N^2} \EE \left[ \vu^\top \mQ^{12}(z) \mR^{23}(\lambda)\vw \right] + \frac{1}{N^2} \EE \left[ \vu^\top \mQ^{12}(z) \vv \vv^\top \mR^{23}(\lambda)\vw \right]\to 0.
\end{align*}
And finally,
\begin{align*}
A_{23} &= - \frac{1}{N^2} \sum_{ijk} \EE \left[ u_i v_j R^{33}_{kk}(\lambda) Q_{ij}^{12} \right] + \frac{1}{N^2} \sum_{ijkl} \EE \left[ u_i v_j w_k w_l R^{33}_{lk}(\lambda) Q_{ij}^{12} \right],\\
&= - \frac1N \EE \left[ \vu^\top \mQ^{12}(z) \vv \frac1N \tr \mR^{33}(\lambda) \right] + \frac{1}{N^2} \EE \left[ \vu^\top \mQ^{12}(z) \vv \vw^\top \mR^{33}(\lambda) \vw \right]\to 0.
\end{align*}
Therefore, 
\begin{align*}
A =  \frac{1}{N\sqrt N}\sum_{ijk} \EE \left[ w_k \frac{\partial Q^{12}_{ij}}{\partial X_{ijk}} \right] + \mathcal{O}(N^{-1}).
\end{align*}
As in the previous subsection, the derivative of $\mQ(z)$ w.r.t.\@ the entry $X_{ijk}$ expresses as $\frac{\partial \mQ(z)}{\partial X_{ijk}} = - \mQ(z) \frac{\partial \mN}{\partial X_{ijk}} \mQ(z)$ but now with
\begin{align*}
\frac{\partial \mN}{\partial X_{ijk}} = \frac{1}{\sqrt N} \begin{bmatrix}
\vzero_{m\times m} & w_k \ve_i^m (\ve_j^n)^\top & v_j \ve_i^m (\ve_k^p)^\top \\
w_k \ve_j^n (\ve_i^m)^\top & \vzero_{n\times n} & u_i \ve_j^n (\ve_k^p)^\top \\
v_j \ve_k^p (\ve_i^m )^\top & u_i \ve_k^p (\ve_j^n )^\top & \vzero_{p\times p}
\end{bmatrix} + \frac{1}{\sqrt N} \Phi_3\left( \tX, \frac{\partial \vu}{\partial X_{ijk}}, \frac{\partial \vv}{\partial X_{ijk}}, \frac{\partial \vw}{\partial X_{ijk}} \right),
\end{align*}
where $\mO\equiv\frac{1}{\sqrt N} \Phi_3\left( \tX, \frac{\partial \vu}{\partial X_{ijk}}, \frac{\partial \vv}{\partial X_{ijk}}, \frac{\partial \vw}{\partial X_{ijk}} \right)$ is of vanishing spectral norm. Indeed, from \eqref{eq:deriv}, there exists $C>0$ independent of $N$ such that $\Vert \frac{\partial \vu}{\partial X_{ijk}} \Vert, \Vert \frac{\partial \vv}{\partial X_{ijk}} \Vert, \Vert \frac{\partial \vw}{\partial X_{ijk}} \Vert \leq \frac{C}{\sqrt N}$ yielding that the spectral norm of $\mO$ is bounded by $\frac{C'}{N^{\frac{3}{2}}}$ for some constant $C'>0$ independent of $N$. Therefore, we find that $A \to - g_1(z) g_2(z) + \mathcal{O}(N^{-1})$ (with $g_1(z), g_2(z)$ the almost sure limits of $\frac1N\tr\mQ^{11}(z)$ and $\frac1N\tr\mQ^{22}(z)$ respectively), thus yielding the same limiting Stieltjes transform as the one obtained in Appendix \ref{proof:semi-circle_independent}. 

\subsection{Proof of Theorem \ref{thm:asymptotics}}\label{proof:asymptotics}
Given the identities in \eqref{eq:singular_vectors}, we have
\begin{align}\label{eq:estimation_ax}
\lambda \langle \vu, \vx \rangle = \vx^\top \tT(\vv)\vw = \beta \langle \vv, \vy \rangle \langle \vw, \vz \rangle + \frac{1}{\sqrt N} \vx^\top \tX(\vv) \vw,
\end{align}
with $\lambda, \langle \vu, \vx \rangle, \langle \vv, \vy \rangle$ and $\langle \vw, \vz \rangle$ converging almost surely to their asymptotic limits $\lambda^\infty(\beta),a_\vx^\infty(\beta), a_\vy^\infty(\beta)$ and $a_\vz^\infty(\beta)$ respectively given the concentration properties in Subsection \ref{sec:concentration}. To characterize such limits we need to evaluate the expectation of $\frac{1}{\sqrt N} \vx^\top \tX(\vv) \vw$. Indeed,
\begin{align*}
\frac{1}{\sqrt N} \EE \left[ \vx^\top \tX(\vv)\vw \right] = \frac{1}{\sqrt N} \sum_{ijk} x_i \EE \left[ v_j w_k X_{ijk}\right] = \frac{1}{\sqrt N} \sum_{ijk} x_i \EE \left[  \frac{\partial v_j }{ \partial X_{ijk}} w_k \right] + \frac{1}{\sqrt N} \sum_{ijk} x_i \EE \left[   v_j \frac{\partial w_k}{\partial X_{ijk}} \right]
\end{align*}
where the last equality is obtained by Stein's lemma (Lemma \ref{lemma:stein}). From \eqref{eq:deriv}, we have
\begin{align*}
\frac{\partial v_j}{\partial X_{ijk}} &= - \frac{1}{\sqrt N} \left( v_j w_k (R_{ij}^{12}(\lambda) - u_i \vu^\top\mR_{:,j}^{12}(\lambda)  ) \right) \\
 &- \frac{1}{\sqrt N} \left(  u_i w_k (R_{jj}^{22}(\lambda) - v_j \vv^\top \mR_{:,j}^{22}(\lambda) ) \right) - \frac{1}{\sqrt N} \left(  u_i v_j (R_{jk}^{23}(\lambda) - w_k \vw^\top \mR_{:,j}^{32}(\lambda) ) \right).
\end{align*}
Hence $A \equiv \frac{1}{\sqrt N} \sum_{ijk} x_i \EE \left[  \frac{\partial v_j }{ \partial X_{ijk}} w_k \right] = A_1 + A_2 + A_3$ decomposes in three terms $A_1, A_2$ and $A_3$. The terms $A_1$ and $A_3$ will be vanishing asymptotically and only $A_2$ contains non-vanishing terms. Indeed, 
\begin{align*}
A_1 &= -\frac1N \sum_{ijk} \EE \left[ x_i w_k v_j w_k R^{12}_{ij}(\lambda) \right] + \frac1N \sum_{ijkl} \EE \left[ x_i w_k v_j w_k u_i u_l R^{12}_{lj}(\lambda) \right],\\
&= - \frac1N \EE\left[ \vx^\top \mR^{12}(\lambda)\vv \right] + \frac1N \EE \left[ \langle \vx, \vu \rangle \vu^\top \mR^{12}(\lambda)\vv \right] \to 0,
\end{align*}
as $N\to \infty$ since $\vx,\vu$ and $\vv$ are of bounded norms and $\mR(\lambda)$ being of bounded spectral norm for $\lambda$ outside the support of $\frac{1}{\sqrt N}\Phi_3(\tX, \vu,\vv, \vw)$ (through the identity in \eqref{eq:R2Q}). Similarly with $A_3$, we have
\begin{align*}
A_3 &= -\frac1N \sum_{ijk} \EE \left[ x_i w_k u_i v_j R^{23}_{jk}(\lambda) \right] + \frac1N \sum_{ijkl} \EE \left[ x_i w_k u_i v_j w_k w_l R_{jl}^{23}(\lambda) \right],\\
&= -\frac1N \EE \left[ \langle \vx, \vu \rangle \vv^\top \mR^{23}(\lambda) \vw \right] + \frac1N \EE \left[ \langle \vx, \vu \rangle  \vv^\top \mR^{23}(\lambda) \vw   \right]\to 0.
\end{align*}
Now $A_2$ is not vanishing, precisely, 
\begin{align*}
A_2 &= -\frac1N \sum_{ijk} \EE \left[ x_i w_k u_i w_k R^{22}_{jj}(\lambda) \right] + \frac1N \sum_{ijkl} \EE \left[ x_i w_k u_i w_k v_j v_l R^{22}_{jl}(\lambda) \right],\\
&= -\frac1N \EE \left[ \langle \vx, \vu \rangle \tr \mR^{22}(\lambda)  \right] + \frac1N \EE \left[ \langle \vx, \vu \rangle \vv^\top \mR^{22}(\lambda)\vv \right],\\
&\to - g_2(\lambda) \EE \left[ \langle \vx, \vu \rangle \right] + \mathcal{O}(N^{-1}),
\end{align*}
where the last line results from the fact that $\frac1N\tr\mR^{22}(\lambda)\asto g_2(\lambda)$ as we saw in the previous subsection. Similarly, we find that 
\begin{align*}
B\equiv \frac{1}{\sqrt N} \sum_{ijk} x_i \EE \left[   v_j \frac{\partial w_k}{\partial X_{ijk}} \right] \to - g_3(\lambda) \EE \left[ \langle \vx, \vu \rangle \right] + \mathcal{O}(N^{-1}).
\end{align*}
Therefore, by \eqref{eq:estimation_ax}, the almost sure limits $\lambda^\infty(\beta),a_\vx^\infty(\beta), a_\vy^\infty(\beta)$ and $a_\vz^\infty(\beta)$ satisfy the equation
\begin{align*}
\lambda^\infty(\beta) a_\vx^\infty(\beta) = \beta a_\vy^\infty(\beta) a_\vz^\infty(\beta) - [g_2(\lambda^\infty(\beta)) + g_3(\lambda^\infty(\beta))] a_\vx^\infty(\beta).
\end{align*}
Hence,
\begin{align*}
a_\vx^\infty(\beta) = \alpha_1(\lambda^\infty(\beta)) a_\vy^\infty(\beta) a_\vz^\infty(\beta),
\end{align*}
with $\alpha_1(z) \equiv \frac{\beta}{z + g(z) - g_1(z)}$. Similarly, we find that
\begin{align*}
\begin{cases}
a_\vy^\infty(\beta)  = \alpha_2(\lambda^\infty(\beta)) a_\vx^\infty(\beta)a_\vz^\infty(\beta),\\
a_\vz^\infty(\beta) = \alpha_3(\lambda^\infty(\beta)) a_\vx^\infty(\beta) a_\vy^\infty(\beta) ,
\end{cases}
\end{align*}
with $\alpha_i(z) \equiv \frac{\beta}{z + g(z) - g_i(z)}$. Solving the above system of equations provides the final asymptotic alignments. Proceeding similarly with \eqref{eq:singular_value} we obtain an estimate of the asymptotic singular value $\lambda^\infty(\beta)$, thereby ending the proof.

\subsection{Proof of Corollary \ref{cor:cubic}}\label{proof:spike_align_cubic}
By Corollary \ref{cor:semi-circle-cubic}, the limiting Stieltjes transform is given by
\begin{align*}
g(z) = -\frac{3z}{4} + \frac{\sqrt{3} \sqrt{3z^2 - 8} }{4},
\end{align*}
and since $g(z) = \sum_{i=1}^3 g_i(z)$ with all $g_i(z)$ equal, then $g_i(z) = \frac{g(z)}{3}$ for all $i\in [3]$. Hence, for each $i\in [3]$, $\alpha_i(z)$ defined in Theorem \ref{thm:asymptotics} writes as
\begin{align*}
\alpha_i(z) = \frac{\beta}{\frac{z}{2} + \frac{\sqrt{3} \sqrt{3 z^{2} - 8}}{6}},
\end{align*}
and $\lambda^\infty(\beta)$ is solution to the equation
\begin{align*}
\frac{z}{4} + \frac{\sqrt{3} \sqrt{3 z^{2} - 8}}{4} - \frac{\left(\frac{z}{2} + \frac{\sqrt{3} \sqrt{3 z^{2} - 8}}{6}\right)^{3}}{\beta^{2}} = 0.
\end{align*}
First we compute the critical value of $\beta$ by solving the above equation in $\beta$ and taking the limit when $z$ tends to the right edge of the semi-circle law (i.e., $z\to 2\sqrt{\frac23}^+$). Indeed, solving the above equation in $\beta$ yields
\begin{align*}
\beta(z) = \sqrt{\frac{2 z^{3} + \frac{2 z^{2} \sqrt{9 z^{2} - 24}}{3} - 4 z - \frac{4 \sqrt{9 z^{2} - 24}}{9}}{z + \sqrt{3} \sqrt{3 z^{2} - 8}}},
\end{align*}
hence
\begin{align*}
\beta_s = \lim_{z\to 2\sqrt{\frac23}^+} \beta(z) = \frac{2\sqrt{3}}{3}.
\end{align*}
Now to express $\lambda^\infty(\beta)$ in terms of $\beta$, we solve the equation $\beta - \beta(z)=0$ in $z$ and choose the unique non-decreasing (in $\beta$) and positive solution, which yields
\begin{align*}
\lambda^\infty(\beta) = \sqrt{\frac{\beta^{2}}{2} + 2 + \frac{\sqrt{3} \sqrt{\left(3 \beta^{2} - 4\right)^{3}}}{18 \beta}}.
\end{align*}
Plugging the above expression of $\lambda^\infty(\beta)$ into the expressions of the asymptotic alignments in Theorem \ref{thm:asymptotics}, we obtain for all $i\in[3]$
\begin{align*}
\alpha_i( \lambda^\infty(\beta) ) = \frac{6 \sqrt{2} \beta}{\sqrt{9 \beta^{2} - 12 + \frac{\sqrt{3} \sqrt{\left(3 \beta^{2} - 4\right)^{3}}}{\beta}} + \sqrt{9 \beta^{2} + 36 + \frac{\sqrt{3} \sqrt{\left(3 \beta^{2} - 4\right)^{3}}}{\beta}}},
\end{align*}
yielding the final result.

\subsection{Proof of Corollary \ref{cor:spectrum_matrices}}\label{proof:spectrum_matrices}
Setting $c_1=c, c_2 = 1-c$ and $c_3 = 0$, we get
\begin{align*}
\begin{cases}
g_1(z) = \frac{g(z) + z}{2} - \frac{\sqrt{4c + (g(z) + z)^2}}{2},\\
g_2(z) = \frac{g(z) + z}{2} - \frac{\sqrt{4(1-c) + (g(z) + z)^2}}{2},\\
g_3(z) = 0.
\end{cases}
\end{align*}
And since $g(z) = \sum_{i=1}^3 g_i(z)$, then $g(z)$ satisfies the equation
\begin{align*}
z - \frac{\sqrt{4c + (g(z) + z)^2}}{2} - \frac{\sqrt{4(1-c) + (g(z) + z)^2}}{2} = 0,
\end{align*}
the solution of which belongs to
\begin{align*}
g(z)\in \left\{ -z - \frac{\sqrt{4c(c-1) + (z^2-1)^2}}{z}, -z + \frac{\sqrt{4c(c-1) + (z^2-1)^2}}{z}   \right\},
\end{align*}
thus the limiting Stieltjes transform (having $\Im[g(z)]>0$) is given by
\begin{align*}
g(z) = -z + \frac{\sqrt{4c(c-1) + (z^2-1)^2}}{z}.
\end{align*}
In particular, the edges of the support of the corresponding limiting distribution $\nu$ are the roots of $4c(c-1)+(z^2-1)^2$, yielding
\begin{align*}
\mathcal{S}(\nu) = \left[ -\sqrt{1 + 2\sqrt{c(1-c)}}, -\sqrt{1 - 2\sqrt{c(1-c)}} \right] \cup \left[ \sqrt{1 - 2\sqrt{c(1-c)}}, \sqrt{1 + 2\sqrt{c(1-c)}} \right].
\end{align*}
And the density function of $\nu$ is obtained by computing the limit $\lim_{\varepsilon\to 0} \vert \Im[g(x+i\varepsilon)]\vert$ yielding
\begin{align*}
\nu(dx) = \frac{1}{\pi x} \sin\left( \frac{\arctan_2(0, q_c(x))}{2} \right) \sqrt{ \vert q_c(x) \vert } \sign \left( \frac{ \sin\left( \frac{\arctan_2(0, q_c(x))}{2} \right) }{x} \right)dx + \left( 1 - 2\min(c,1-c) \right)\delta(x),
\end{align*}
where $q_c(x) = (x^2 - 1)^2 + 4c(c-1)$. The Dirac $\delta(x)$ component in the above expression corresponds to the fact that the corresponding matrix model is of rank $2\min(m,n)$.
\subsection{Proof of Corollary \ref{cor:spike_align_matrices}}\label{proof:spike_align_matrices}
From Subsection \ref{proof:spectrum_matrices}, plugging the expression of the limiting Stieltjes transform $g(z)$ into the expressions defining $g_1(z)$ and $g_2(z)$ yields
\begin{align*}
\begin{cases}
g_1(z) = \frac{\sqrt{4c(c-1) + (z^2-1)^2}}{2z} - \frac{\sqrt{4c(c-1+z^2) + (z^2-1)^2}}{2z},\\
g_2(z) = \frac{\sqrt{4c(c-1) + (z^2-1)^2}}{2z} - \frac{\sqrt{4c(c-1-z^2) + (z^2+1)^2}}{2z}.
\end{cases}
\end{align*}
Thus, by Theorem \ref{thm:asymptotics}, we have
\begin{align*}
\begin{cases}
\alpha_1(z) = \frac{\beta}{\frac{\sqrt{4 c^{2} - 4 c + z^{4} - 2 z^{2} + 1}}{2 z} + \frac{\sqrt{4 c^{2} + 4 c z^{2} - 4 c + z^{4} - 2 z^{2} + 1}}{2 z}}\\
\alpha_2(z) = \frac{\beta}{\frac{\sqrt{4 c^{2} - 4 c + z^{4} - 2 z^{2} + 1}}{2 z} + \frac{\sqrt{4 c^{2} - 4 c z^{2} - 4 c + z^{4} + 2 z^{2} + 1}}{2 z}}\\
\alpha_3(z) = \frac{\beta z}{\sqrt{4 c^{2} - 4 c + z^{4} - 2 z^{2} + 1}}.
\end{cases}
\end{align*}
Then solving the equation $z+g(z) - \frac{\beta}{\alpha_1(z)\alpha_2(z)\alpha_3(z)}$ in $z$ provides the almost sure limit of $\lambda$ in terms of $\beta$. Specifically, 
\begin{align*}
\lambda^\infty(\beta) = \sqrt{\beta^2 + 1 + \frac{c(1-c)}{\beta^2}}.
\end{align*}
Plugging the above expression of $\lambda^\infty(\beta)$ into $\frac{1}{\sqrt{\alpha_2(z)\alpha_3(z)}}$ by replacing $z$ with the expression of $\lambda^\infty(\beta)$ provides the asymptotic alignment $a_\vx^\infty(\beta)$ as $ \frac{1}{\kappa(\beta, c)}$, with $\kappa(\beta, c) = \sqrt{\alpha_2(\lambda^\infty(\beta))\alpha_3(\lambda^\infty(\beta))}$ given by the following expression
\begin{tiny}
\begin{align*}
\kappa(\beta, c) = \sqrt{\frac{ 2\sqrt{\frac{\beta^{2} \left(\beta^{2} + 1\right) - c^{2} + c}{\beta^{2}}} \left(\beta^{2} \left(\beta^{2} + 1\right) - c \left(c - 1\right)\right)}{\sqrt{\frac{4 \beta^{4} c \left(c - 1\right) + \left(\beta^{4} - c^{2} + c\right)^{2}}{\beta^{4}}} \sqrt{\frac{\beta^{2} \left(\beta^{2} + 1\right) - c \left(c - 1\right)}{\beta^{2}}} \left(\sqrt{\frac{- 4 \beta^{2} c \left(\beta^{2} \left(\beta^{2} - c + 2\right) - c \left(c - 1\right)\right) + \left(\beta^{2} \left(\beta^{2} + 2\right) - c \left(c - 1\right)\right)^{2}}{\beta^{4}}} + \sqrt{\frac{4 \beta^{4} c \left(c - 1\right) + \left(\beta^{4} - c \left(c - 1\right)\right)^{2}}{\beta^{4}}}\right)}}.
\end{align*}
\end{tiny}
From the identities
\begin{align*}
&4\beta^4c(c-1) + (\beta^4 - c(c-1))^2 = (\beta^4 + c(c-1))^2 \\
&- 4 \beta^{2} c \left(\beta^{2} \left(\beta^{2} - c + 2\right) - c \left(c - 1\right)\right) + \left(\beta^{2} \left(\beta^{2} + 2\right) - c \left(c - 1\right)\right)^{2} = \left(\beta^{4} + \beta^{2} \left(2 - 2 c\right) - c^{2} + c\right)^{2}
\end{align*}
$\kappa(\beta, c)$ simplifies as
\begin{align*}
\kappa(\beta, c) = \beta \sqrt{\frac{   \beta^{2} \left(\beta^{2} + 1\right) - c \left(c - 1\right)  }{ (\beta^4 + c(c-1)) \left( \beta^{2} + 1 -  c \right)}}.
\end{align*}

The asymptotic alignment $a_\vy^\infty(\beta)$ is given by  $ \frac{1}{\kappa(\beta, 1-c)}$ since the dimension ratio of the component $\vy$ is $c_2=1-c$. Moreover, the critical value of $\beta$ is obtained by solving the equation $\lambda^\infty(\beta_s) = \sqrt{1 + 2\sqrt{c(1-c)}}$, i.e., when the limiting singular value gets closer to the right edge of the support of the corresponding limiting spectral measure (see Corollary \ref{cor:spectrum_matrices}). Finally, similarly as above, we can check that $\sqrt{ \alpha_1(\lambda^\infty(\beta))\alpha_2(\lambda^\infty(\beta)) }$ is equal to $1$ for $\beta\geq \beta_s$ (i.e., the alignment along the third dimension is $1$, corresponding to the component $\vz$ in \eqref{eq:spiked_tensor_model}).

\subsection{Proof of Theorem \ref{thm:semi-circle_independent_d_order}}\label{proof:circle_independent_d_order}
Denote the matrix model as
\begin{align*}
\mN \equiv \frac{1}{\sqrt N} \Phi_d(\tX, \va^{(1)}, \ldots, \va^{(d)}),
\end{align*}
with $\tX\sim \sT_{n_1,\ldots,n_d}(\mathcal{N}(0, 1))$ and $(\va^{(1)}, \ldots, \va^{(d)})\in \sS^{n_1-1}\times \cdots \times \sS^{n_d-1}$ are independent of $\tX$. We further denote the resolvent of $\mN$ as
\begin{align*}
\mQ(z) \equiv \left( \mN - z \mI_N \right)^{-1} = \begin{bmatrix}
\mQ^{11}(z) & \mQ^{12}(z) & \mQ^{13}(z) & \cdots & \mQ^{1d}(z)\\
\mQ^{12}(z)^\top & \mQ^{22}(z) & \mQ^{23}(z) & \cdots & \mQ^{2d}(z)\\
\mQ^{13}(z)^\top & \mQ^{23}(z)^\top & \mQ^{33}(z) & \ldots & \mQ^{3d}(z)\\
\vdots & \vdots & \vdots & \ddots & \vdots \\
\mQ^{1d}(z)^\top & \mQ^{2d}(z)^\top & \mQ^{3d}(z)^\top & \cdots & \mQ^{dd}(z)
\end{bmatrix}.
\end{align*}
By Borel-Cantelli lemma, we have $\frac1N\tr \mQ(z) \asto g(z)$ and for all $i\in [d]$, $\frac1N\tr \mQ^{ii}(z) \asto g_i(z)$. Applying the identity in \eqref{eq:resolvent_identity} to the symmetric matrix $\mN$ we get $\mN \mQ(z) - z\mQ(z) = \mI_N$, from which we particularly get
\begin{align*}
\frac{1}{\sqrt N} \sum_{j=2}^d \left[ \tX^{1j} \mQ^{1j}(z)^\top \right]_{i_1 i_1} - z Q_{i_1 i_1}^{11}(z) = 1,
\end{align*}
or
\begin{align}\label{eq:g1orderd}
\frac{1}{N \sqrt N}  \sum_{j=1}^2 \sum_{i_1=1}^{n_1} \left[ \tX^{1j} \mQ^{1j}(z)^\top \right]_{i_1 i_1} -\frac{z}{N}\tr \mQ^{11}(z) = \frac{n_1}{N}, 
\end{align}
where we recall that $\tX^{ij} \equiv \tX( \va^{(1)}, \ldots, \va^{(i-1)}, :, \va^{(i+1)}, \ldots, \va^{(j-1)}, :, \va^{(j+1)}, \ldots, \va^{(d)} ) \in \sM_{n_i, n_j}$.

We thus need to compute the expectation of $\frac{1}{N \sqrt N} \sum_{i_1=1}^{n_1} \left[ \tX^{1j} \mQ^{1j}(z)^\top \right]_{i_1 i_1}$ which develops as
\begin{align*}
A_j\equiv  \frac{1}{N \sqrt N} \sum_{i_1=1}^{n_1} \EE \left[ \tX^{1j} \mQ^{1j}(z)^\top \right]_{i_1 i_1} = \frac{1}{N\sqrt N} \sum_{i_1 i_j} \EE \left[ \left[ \tX^{1j} \right]_{i_1 i_j} Q^{1j}_{i_1 i_j} \right] = \frac{1}{N \sqrt{N}} \sum_{i_1,\ldots,i_d} \prod_{k\neq 1, k\neq j} a_{i_k}^{(k)} \EE \left[  \frac{\partial Q^{1j}_{i_1 i_j}}{\partial X_{i_1,\ldots,i_d} } \right]
\end{align*}
where the last equality follows from Stein's lemma (Lemma \ref{lemma:stein}). In particular, as in Appendix \ref{proof:semi-circle_independent} for the $3$-order case, it turns out that the only contributing term in the derivative $ \frac{\partial Q_{i_1 i_j}^{1j}}{\partial X_{i_1,\ldots, i_d}} $ is $- \frac{1}{\sqrt N} \prod_{k\neq 1, k\neq j} a_{i_k}^{(k)} Q_{i_1 i_1}^{11} Q_{i_j i_j}^{jj}  $ with the other terms yielding quantities of order $\mathcal{O}(N^{-1})$. Therefore, we find that
\begin{align*}
A_j &= \frac{1}{N \sqrt{N}} \sum_{i_1,\ldots,i_d} \prod_{k\neq 1, k\neq j} a_{i_k}^{(k)} \EE \left[  \frac{\partial Q^{1j}_{i_1 i_j}}{\partial X_{i_1,\ldots,i_d} } \right] = -\frac{1}{N^2} \sum_{i_1,\ldots,i_d}   \prod_{k\neq 1, k\neq j} \left( a_{i_k}^{(k)} \right)^2 \EE \left[ Q_{i_1 i_1}^{11} Q_{i_j i_j}^{jj} \right] + \mathcal{O}(N^{-1})\\
&=-\frac{1}{N^2} \sum_{i_1 i_j}  \EE \left[ Q_{i_1 i_1}^{11} Q_{i_j i_j}^{jj} \right] + \mathcal{O}(N^{-1})\\
&= -\frac1N \tr \mQ^{11}(z) \frac1N \tr \mQ^{jj}(z) + \mathcal{O}(N^{-1}) \\
&\asto - g_1(z) g_j(z) + \mathcal{O}(N^{-1}).
\end{align*}
From \eqref{eq:g1orderd}, $g_1(z)$ satisfies $-g_1(z) \sum_{j\neq 1} g_j(z) - z g_1(z) = c_1$ with $c_1=\lim \frac{n_1}{N}$. Similarly, for all $i\in [d]$, $g_i(z)$ satisfies $-g_i(z) \sum_{j\neq i} g_j(z) - z g_i(z) = c_i$ with $c_i=\lim \frac{n_i}{N}$. And since, $g(z) = \sum_{i=1}^d g_i(z)$, we have for each $i\in [d]$, $g_i(z) (g(z) - g_i(z)) + zg_i(z) + c_i = 0$, yielding 
\begin{align*}
g_i(z) = \frac{g(z) + z}{2} - \frac{\sqrt{4c_i + (g(z) + z)^2}}{2},
\end{align*}
with $g(z)$ solution to the equation $g(z) = \sum_{i=1}^d g_i(z)$ satisfying $\Im[g(z)]>0$ for $\Im[z]>0$.

\subsection{Proof of Corollary \ref{cor:semi_circle_order_d}}\label{proof:semi_circle_order_d}
Given Theorem \ref{thm:semi-circle_independent_d_order} and setting $c_i=\frac1d$ for all $i\in [d]$, we have $g_i(z) = \frac{g(z)}{d}$, thus $g(z)$ satisfies the equation
\begin{align*}
\frac{g(z)}{d} = \frac{g(z) + z}{2} - \frac{\sqrt{\frac4d + (g(z) + z)^2}}{2}.
\end{align*}
Solving in $g(z)$ yields
\begin{align*}
g(z) \in \left\{- \frac{d z}{2 \left(d - 1\right)} - \frac{\sqrt{d \left(d z^{2} - 4 d + 4\right)}}{2 \left(d - 1\right)}, - \frac{d z}{2 \left(d - 1\right)} + \frac{\sqrt{d \left(d z^{2} - 4 d + 4\right)}}{2 \left(d - 1\right)}\right\},
\end{align*}
and the limiting Stieltjes transform with $\Im[g(z)]\geq 0$ is
\begin{align*}
g(z) = - \frac{d z}{2 \left(d - 1\right)} + \frac{\sqrt{d \left(d z^{2} - 4 d + 4\right)}}{2 \left(d - 1\right)}.
\end{align*}

{\color{cblue}
\subsection{Existence of $g(z)$}\label{proof_existence_of_g}
Define the following function for $(g, z)\in \sR\times \sR_+$ and $d\geq 3$
\begin{align*}
    k(g, z) = g - \frac{g+z}{2} d + \frac12 \sum_{i=1}^d \sqrt{4c_i + (g+z)^2},
\end{align*}
with $\sum_{i=1}^d c_i = 1$. By concavity of the function $\sqrt{\cdot}$ we have $\sum_{i=1}^d \sqrt{4c_i + x^2}\leq d \sqrt{\frac{4}{d} + x^2}$ for all $x\in \sR$. Therefore, $k(g, z)$ is bounded as
\begin{align*}
    k(g, z) \leq \bar k (g, z) \equiv g - \frac{g+z}{2} d + \frac{d}{2} \sqrt{\frac{4}{d} + (g + z)^2}.
\end{align*}
From Section \ref{proof:semi_circle_order_d}, we have, for any $z> 2 \sqrt{\frac{d-1}{d}}$, there exists $g_*\in \sR$ such that $\bar k(g_*, z) = 0$. Besides, we further have $\lim_{g\to -\infty} k(g, z) = \lim_{g\to +\infty} k(g, z) = +\infty $ and hence by continuity, there exists $g(z)$ such that $k(g(z), z) = 0$ for $z$ large enough (e.g., $z > 2 \sqrt{\frac{d-1}{d}}$).
}
\subsection{Proof of Theorem \ref{thm:semi-circle_dependent_order_d}}\label{proof:semi-circle_dependent_order_d}
Given the random tensor model in \eqref{eq_asymmetric_spike_model} and its singular vectors characterized by \eqref{eq:singular_vectors_id_order_d}, we denote the associated random matrix model as
\begin{align*}
\mT \equiv \Phi_d \left( \tT, \vu^{(1)}, \ldots , \vu^{(d)}\right) = \beta \mV \mB \mV^\top + \mN,
\end{align*}
where $\mN = \frac{1}{\sqrt N} \Phi_d \left( \tX, \vu^{(1)}, \ldots , \vu^{(d)}\right)$, $\mB \in \sM_d$ with entries $B_{ij} = (1-\delta_{ij}) \prod_{k\neq i, j} \langle \vx^{(k)}, \vu^{(k)} \rangle$ 
\begin{align*}
\mV \equiv \begin{bmatrix}
\vx^{(1)} & \vzero_{n_1} &  \cdots & \vzero_{n_1}\\
\vzero_{n_2} &  \vx^{(2)} &  \cdots & \vzero_{n_2}\\
\vdots & \vdots & \ddots   & \vzero_{n_3} \\
\vzero_{n_d} & \vzero_{n_d} & \vzero_{n_d} & \vx^{(d)}
\end{bmatrix} \in \sM_{N, d}.
\end{align*}
We further denote the resolvent of $\mT$ and $\mN$ respectively as
\begin{align*}
\mR(z) \equiv \left( \mT - z \mI_N \right)^{-1} = \begin{bmatrix}
\mR^{11}(z) & \mR^{12}(z) & \mR^{13}(z) & \cdots & \mR^{1d}(z)\\
\mR^{12}(z)^\top & \mR^{22}(z) & \mR^{23}(z) & \cdots & \mR^{2d}(z)\\
\mR^{13}(z)^\top & \mR^{23}(z)^\top & \mR^{33}(z) & \ldots & \mR^{3d}(z)\\
\vdots & \vdots & \vdots & \ddots & \vdots \\
\mR^{1d}(z)^\top & \mR^{2d}(z)^\top & \mR^{3d}(z)^\top & \cdots & \mR^{dd}(z)
\end{bmatrix}.
\end{align*}
\begin{align*}
\mQ(z) \equiv \left( \mN - z \mI_N \right)^{-1} = \begin{bmatrix}
\mQ^{11}(z) & \mQ^{12}(z) & \mQ^{13}(z) & \cdots & \mQ^{1d}(z)\\
\mQ^{12}(z)^\top & \mQ^{22}(z) & \mQ^{23}(z) & \cdots & \mQ^{2d}(z)\\
\mQ^{13}(z)^\top & \mQ^{23}(z)^\top & \mQ^{33}(z) & \ldots & \mQ^{3d}(z)\\
\vdots & \vdots & \vdots & \ddots & \vdots \\
\mQ^{1d}(z)^\top & \mQ^{2d}(z)^\top & \mQ^{3d}(z)^\top & \cdots & \mQ^{dd}(z)
\end{bmatrix}.
\end{align*}
Similarly as in $3$-order case, by Woodbury matrix identity (Lemma \ref{lem:woodbury}), we have
\begin{align*}
\frac1N \tr \mR(z) &= \frac1N \tr \mQ(z) - \frac1N \tr \left[ \left( \frac{1}{\beta} \mB^{-1} + \mV^\top \mQ(z) \mV \right)^{-1} \mV^\top \mQ^2(z) \mV \right],\\
&= \frac1N \tr \mQ(z)  + \mathcal{O}(N^{-1}),
\end{align*}
since the perturbation matrix $ \left( \frac{1}{\beta} \mB^{-1} + \mV^\top \mQ(z) \mV \right)^{-1} \mV^\top \mQ^2(z) \mV$ is of bounded spectral norm (if $\Vert \mQ(z)\Vert $ is bounded, see condition in \eqref{eq:resovent_bounded}) and has finite size ($d\times d$ matrix). Therefore, the characterization of the spectrum of $\mT$ boils down to the estimation of $\frac1N \tr \mQ(z)$. Now we are left to handle the statistical dependency between the tensor noise $\mX$ and the singular vectors of $\tT$. Recalling the proof of Appendix \ref{proof:circle_independent_d_order}, we have again
\begin{align*}
\frac{1}{N \sqrt N}  \sum_{j=1}^2 \sum_{i_1=1}^{n_1} \left[ \tX^{1j} \mQ^{1j}(z)^\top \right]_{i_1 i_1} -\frac{z}{N}\tr \mQ^{11}(z) = \frac{n_1}{N}, 
\end{align*}
with $\tX^{ij} \equiv \tX( \vu^{(1)}, \ldots, \vu^{(i-1)}, :, \vu^{(i+1)}, \ldots, \vu^{(j-1)}, :, \vu^{(j+1)}, \ldots, \vu^{(d)} ) \in \sM_{n_i, n_j}$. Taking the expectation of $\frac{1}{N \sqrt N} \sum_{i_1=1}^{n_1} \left[ \tX^{1j} \mQ^{1j}(z)^\top \right]_{i_1 i_1}$ yields
\begin{align*}
A_j = \frac{1}{N \sqrt{N}} \sum_{i_1,\ldots,i_d}  \EE \left[ \frac{\partial  Q^{1j}_{i_1 i_j}  }{\partial X_{i_1,\ldots,i_d} } \prod_{k\neq 1, k\neq j} u_{i_k}^{(k)}  \right] + \frac{1}{N \sqrt{N}} \sum_{i_1,\ldots,i_d} \EE \left[ Q^{1j}_{i_1 i_j} \prod_{k\neq 1, k\neq j} \frac{\partial  u_{i_k}^{(k)}   }{\partial X_{i_1,\ldots,i_d} }  \right] = A_{j1}  +A_{j2},
\end{align*}
where we already computed the first term ($A_{j1}$) in Appendix \ref{proof:circle_independent_d_order}. Now we will show that $A_{j2}$ is asymptotically vanishing under Assumption \ref{ass:growth_order_d}. Indeed, by \eqref{eq:deriv_order_d}, the higher order terms arise from the term $ - \frac{1}{\sqrt N} \prod_{\ell \neq k} u^{(\ell)}_{i_\ell} R^{kk}_{i_k i_k}(\lambda)$, we thus only show that the contribution of this term is also vanishing. Precisely,
\begin{align*}
A_{j2} &= - \frac{1}{N^2} \sum_{i_1,\ldots,i_d} \EE \left[ Q_{i_1 i_j}^{1j} \prod_{k\neq 1, k\neq j} \prod_{\ell \neq k} u^{(\ell)}_{i_\ell} R^{kk}_{i_k i_k}(\lambda)  \right] + \mathcal{O}(N^{-1})\\
&= -\frac{1}{N^2} \sum_{i_1,\ldots, i_d} \EE \left[ \left(u_{i_1}^{(1)}\right)^{d-2} Q_{i_1 i_j}^{1j} \left(u_{i_j}^{(j)}\right)^{d-2} \prod_{k\neq 1, k\neq j}  u^{(k)}_{i_k} R^{kk}_{i_k i_k}(\lambda) \right] + \mathcal{O}(N^{-1})\\
&= -\frac{1}{N^2}\EE \left[ \left( (\vu^{(1)})^{\odot d-2} \right)^\top \mQ^{1j}(z)  \left( \vu^{(j)} \right)^{\odot d-2} \sum_{i_2, \ldots, i_{j-1}, i_{j+1}, \ldots, i_d} \prod_{k\neq 1, k\neq j}  u^{(k)}_{i_k} R^{kk}_{i_k i_k}(\lambda) \right] + \mathcal{O}(N^{-1})
\end{align*}
where $\va^{\odot q}$ denotes the vector with entries $a_i^q$. As such, $A_{j2}$ is vanishing asymptotically since the singular vectors $\vu^{(1)}, \ldots, \vu^{(d)}$ are unitary and since their entries are bounded by $1$.

We finally need to check that the derivative of $\mQ(z)$ w.r.t.\@ the entry $X_{i_1,\ldots,i_d}$ has the same expression asymptotically as the one in the independent case of Appendix \ref{proof:circle_independent_d_order}. Indeed, we have $\frac{\partial \mQ(z)}{\partial X_{i_1,\ldots, i_d}} = - \mQ(z) \frac{\partial \mN}{\partial X_{i_1,\ldots, i_d}} \mQ(z) $ with
\begin{align*}
\frac{\partial \mN}{\partial X_{i_1,\ldots, i_d}} = \begin{bmatrix}
\vzero_{n_1\times n_1} & \mC_{12} &  \cdots & \mC_{1d}\\
\mC_{12}^\top &  \vzero_{n_2\times n_2} &  \cdots & \mC_{2d}\\
\vdots & \vdots & \ddots   & \vdots \\
\mC_{1d}^\top & \mC_{2d}^\top & \cdots & \vzero_{n_d\times n_d}
\end{bmatrix} + \frac{1}{\sqrt N} \Phi_d \left( \mX, \frac{\partial \vu^{(1)}}{\partial X_{i_1,\ldots, i_d}} , \ldots, \frac{\partial \vu^{(d)}}{\partial X_{i_1,\ldots, i_d}} \right),
\end{align*}
where $\mC_{ij} = \prod_{k\neq i, k\neq j} u^{(k)}_{i_k} \ve_{i_i}^{n_i} (\ve_{i_j}^{n_j})^\top$ and $\mO = \frac{1}{\sqrt N} \Phi_d \left( \mX, \frac{\partial \vu^{(1)}}{\partial X_{i_1,\ldots, i_d}} , \ldots, \frac{\partial \vu^{(d)}}{\partial X_{i_1,\ldots, i_d}} \right) $ is of vanishing norm. Indeed, by \eqref{eq:deriv_order_d}, there exists $C>0$ independent of $N$ such that $ \Vert \frac{\partial \vu^{(1)}}{\partial X_{i_1,\ldots, i_d}} \Vert , \ldots, \Vert \frac{\partial \vu^{(d)}}{\partial X_{i_1,\ldots, i_d}} \Vert \leq \frac{C}{\sqrt N}$, therefore, the spectral norm of $\mO$ is bounded by $\frac{C'}{N^{\frac{d}{2}}}$ for some constant $C'>0$ independent of $N$. Finally, $A_j \to -g_1(z) g_j(z) + \mathcal{O}(N^{-1})$ (with $g_1(z), g_j(z)$ the almost sure limits of $\frac1N \tr \mQ^{11}(z)$ and $\frac1N\tr \mQ^{jj}(z)$ respectively), hence yielding the same limiting Stieltjes transform as the one obtained in Appendix \ref{proof:circle_independent_d_order}.
\subsection{Proof of Theorem \ref{thm:asymptotics_general}}\label{proof:asymptotics_general}
Given the identities in \eqref{eq:singular_vectors_id_order_d}, we have for all $i\in [d]$
\begin{align*}
\frac{1}{\sqrt{N}}\sum_{j_1,\ldots,j_d} x_{j_i}^{(i)} \prod_{k\neq i} u_{j_k}^{(k)} X_{j_1,\ldots,j_d} + \beta \prod_{k\neq i } \langle \vu^{(k)}, \vx^{(k)}\rangle = \lambda \langle \vu^{(i)}, \vx^{(i)}\rangle,
\end{align*}
with $\lambda$ and $\langle \vu^{(i)}, \vx^{(i)}\rangle$ concentrate almost surely around their asymptotic denoted $\lambda^\infty(\beta)$ and $a_{\vx^{(i)}}^\infty(\beta)$ respectively. Taking the expectation of the first term and applying Stein's lemma (Lemma \ref{lemma:stein}), we get
\begin{align*}
A_i = \frac{1}{\sqrt{N}}\sum_{j_1,\ldots,j_d} x_{j_i}^{(i)} \EE \left[ \frac{\partial \left( \prod_{k\neq i} u_{j_k}^{(k)} \right)}{\partial X_{j_1,\ldots,j_d} } \right]  =  \frac{1}{\sqrt{N}}\sum_{j_1,\ldots,j_d} x_{j_i}^{(i)} \sum_{k\neq i} \EE \left[ \frac{\partial u_{j_k}^{(k)}}{ \partial X_{j_1,\ldots,j_d} } \prod_{\ell\neq k, \ell \neq i} u_{j_\ell}^{(\ell)} \right],
\end{align*}
where the only contributing term in the expression of $\frac{\partial u_{j_k}^{(k)}}{ \partial X_{j_1,\ldots,j_d} } $ from \eqref{eq:deriv_order_d} is $-\frac{1}{\sqrt N} \prod_{\ell \neq k} u_{j_\ell}^{(\ell)} R_{j_k j_k}^{kk}(\lambda)$ which yields
\begin{align*}
A_i &= -\frac{1}{N} \sum_{j_1,\ldots,j_d} x^{(i)}_{j_i} \sum_{k\neq i} \EE \left[ R_{j_k j_k}^{kk}(\lambda) \prod_{\ell\neq k} u_{j_\ell}^{(\ell)} \prod_{\ell\neq k, \ell \neq i} u_{j_\ell}^{(\ell)} \right] + \mathcal{O}(N^{-1})\\
&= - \EE \left[ \langle \vu^{(i)}, \vx^{(i)}\rangle \sum_{k\neq i} \frac1N \tr \mR^{kk}(\lambda)   \right] + \mathcal{O}(N^{-1})\\
&\to -\EE \left[ \langle \vu^{(i)}, \vx^{(i)}\rangle  \right] \sum_{k\neq i} g_k(\lambda).
\end{align*}
Therefore, the almost sure limits $\lambda^\infty$ and $a_{\vx^{(i)}}^\infty$ for each $i\in [d]$ satisfy
\begin{align*}
\lambda^\infty a_{\vx^{(i)}}^\infty = \beta \prod_{k\neq i} a_{\vx^{(k)}}^\infty - a_{\vx^{(i)}}^\infty  \sum_{k\neq i} g_k(\lambda^\infty),
\end{align*}
therefore
\begin{align*}
a_{\vx^{(i)}}^\infty = \alpha_i(\lambda^\infty) \prod_{k\neq i} a_{\vx^{(k)}}^\infty, \quad \text{with}\quad \alpha_i(z) = \frac{\beta}{z + g(z) - g_i(z)},
\end{align*}
since $g(z) = \sum_{k=1}^d g_i(z)$. To solve the above equation, we simply write $x_i = a_{\vx^{(i)}}^\infty$ and $\alpha_i = \alpha_i(z)$ by omitting the dependence on $z$. We therefore have
\begin{align*}
x_i = \alpha_i \prod_{k\neq i} x_k\quad \Rightarrow \quad x_i = \alpha_i x_j \prod_{k\neq i, k\neq j} x_k = \alpha_i x_j \prod_{k\neq i, k\neq j} \left( \alpha_k \prod_{\ell\neq k} x_\ell \right),
\end{align*} 
from which we have
\begin{align*}
x_i = x_j \left(\prod_{k\neq j} \alpha_k\right) \left( \prod_{k\neq i, k\neq j} \prod_{\ell \neq k} x_\ell \right) = x_j  \left(\prod_{k\neq j} \alpha_k\right) \left( \prod_{k\neq i, k\neq j} \prod_{\ell \neq k,\ell\neq i} x_\ell \right) x_i^{d-2},
\end{align*}
thus
\begin{align*}
x_j  \left(\prod_{k\neq j} \alpha_k\right) \left( \prod_{k\neq i, k\neq j} \prod_{\ell \neq k,\ell\neq i} x_\ell \right) x_i^{d-3} = 1,
\end{align*}
and we remark that $\left( \prod_{k\neq i, k\neq j} \prod_{\ell \neq k,\ell\neq i} x_\ell \right) x_i^{d-3} = \left( \frac{x_j}{\alpha_j} \right)^{d-3} x_j^{d-2}$, hence $x_j$ is given by
\begin{align*}
x_j = \left(\frac{\alpha_j^{d-3}}{\prod_{k\neq j} \alpha_k}\right)^{\frac{1}{2d-4}},
\end{align*}
which ends the proof.
{\color{cblue}
\paragraph*{Alternative expression of $q_i(z)$} From the above, we have $\lambda^\infty + g(\lambda^\infty) = \beta \prod_{i=1}^d x_i$ and
\begin{align*}
    (\lambda^\infty + g(\lambda^\infty) - g_i(\lambda^\infty)) x_i = \beta \prod_{j\neq i}^d x_j \quad \Rightarrow \quad (\lambda^\infty + g(\lambda^\infty) - g_i(\lambda^\infty)) x_i^2 = \beta \prod_{ i=1}^d x_i.
\end{align*}
Therefore, 
\begin{align*}
    x_i = \sqrt{ \frac{ \lambda^\infty + g(\lambda^\infty)  }{ \lambda^\infty + g(\lambda^\infty) - g_i(\lambda^\infty) } } = \sqrt{ 1 + \frac{g_i(\lambda^\infty) }{ \lambda^\infty + g(\lambda^\infty) - g_i(\lambda^\infty) } },
\end{align*}
with $z + g(z) - g_i(z) = \frac{-c_i}{g_i(z)}$ since $g_i(z)$ satisfies $g_i^2(z) - (g(z) + z) g_i(z) - c_i = 0$. Hence, we find
\begin{align*}
    x_i = \sqrt{ 1 - \frac{g_i^2(\lambda^\infty)}{c_i} }.
\end{align*}
}
\subsection{Additional lemmas}
{\color{cblue}
\begin{lemma}[Spectral norm of random Gaussian tensors]\label{lemma_spectral_norm_of_X} Let $\tX\sim \sT_{n_1, \ldots, n_d}(\mathcal{N}(0, 1))$, then the spectral norm of $\Vert \tX \Vert$ can be bounded, with probability at least $1-\delta$ for $\delta>0$, as
\begin{align*}
    \Vert \tX \Vert \leq \sqrt{ 2 \left[ \left( \sum_{i=1}^d n_i \right) \log\left( \frac{2d}{\log(3/2)}\right) + \log\left( \frac{2}{\delta} \right) \right] }.
\end{align*}
\end{lemma}
\begin{proof}
By definition, the spectral norm of $\tX$ is given as
\begin{align}\label{eq_spectral_norm}
    \Vert \tX \Vert = \sup_{\vu^{(i)} \in \sS^{n_i - 1},\, i\in [d]} \tX (\vu^{(1)}, \ldots, \vu^{(d)})
\end{align}
For some $\varepsilon>0$, let $\gE_1, \ldots, \gE_d$ be $\varepsilon$-nets of $\sS^{n_1-1}, \ldots, \sS^{n_d-1}$ respectively. Since $\sS^{n_1-1}\times \cdots \times \sS^{n_d-1}$ is compact, there exists a maximizer of \eqref{eq_spectral_norm}. And with the $\varepsilon$-net argument, there exists $\ve^{(i)}\in \gE_i$ for each $i\in [d]$ such that
\begin{align*}
    \Vert \tX \Vert = \tX(\ve^{(1)} + \vdelta^{(1)}, \ldots, \ve^{(d)} + \vdelta^{(d)}),
\end{align*}
such that $\Vert \vdelta^{(i)} \Vert \leq \varepsilon$ for $i\in [d]$. Therefore, one has
\begin{align*}
    \Vert \tX \Vert \leq \tX(\ve^{(1)} , \ldots, \ve^{(d)}) + \left( \sum_{i=1}^d \varepsilon^i \binom{d}{i} \right) \Vert \tX \Vert.
\end{align*}
For $\varepsilon = \frac{\log(3/2)}{d}$, one has $\sum_{i=1}^d \varepsilon^i \binom{d}{i} \leq \sum_{i=1}^d \frac{(\varepsilon d)^i}{i!} \leq e^{\varepsilon d} -1 = \frac12 $. As such, the spectral norm of $\tX$ can be bounded as
\begin{align}\label{eq_bound_spectral_norm_e_net}
    \Vert \tX \Vert \leq 2 \max_{\ve^{(i)} \in \gE_i,\, i\in [d]} \tX(\ve^{(1)}, \ldots, \ve^{(d)}).
\end{align}
Since the entries of $\tX$ are i.i.d. standard Gaussian random variables, we have $\EE\left[ e^{tX_{i_1\ldots i_d}} \right] \leq e^{t^2/2}$, hence using Hoeffding's inequality we have for any $\vu^{(i)} \in \sS^{n_i - 1}$ with $i\in [d]$
\begin{align*}
    &\sP\left\lbrace \tX(\vu^{(1)},\ldots, \vu^{(d)}) \geq t \right\rbrace = \sP \left\lbrace e^{s \tX(\vu^{(1)},\ldots, \vu^{(d)})} \geq e^{st} \right\rbrace\leq e^{-st} \EE \left[  e^{s \tX(\vu^{(1)},\ldots, \vu^{(d)})} \right]\\
    &\leq \exp \left\lbrace -st + \frac{s^2}{2} \sum_{i_1, \ldots, i_d} (u_{i_1}^{(1)})^2 \cdots (u_{i_d}^{(d)})^2 \right\rbrace = \exp\left\lbrace -st + \frac{s^2}{2}\right\rbrace.
\end{align*}
Minimizing the right-hand side w.r.t. $s$ yields $\sP\left\lbrace \tX(\vu^{(1)},\ldots, \vu^{(d)}) \geq t \right\rbrace \leq e^{-t^2/2}$. With the same arguments, we further have $\sP\left\lbrace \tX(\vu^{(1)},\ldots, \vu^{(d)}) \leq -t \right\rbrace \leq e^{-t^2/2}$. And taking the union of the two cases yields
\begin{align*}
    \sP\left\lbrace \vert\tX(\vu^{(1)},\ldots, \vu^{(d)}) \vert \geq t \right\rbrace \leq 2 e^{-t^2/2}.
\end{align*}
Back to \eqref{eq_bound_spectral_norm_e_net}, since $\vert \gE_i \vert \leq \left( \frac{2}{\varepsilon}\right)^{n_i}$, involving the union bound gives us
\begin{align*}
    \sP \left\lbrace \Vert \tX \Vert \geq t \right\rbrace &\leq \sum_{\ve^{(1)} \in \gE_1, \ldots, \ve^{(d)} \in \gE_d} \sP \left\lbrace \tX(\ve^{(1)}, \ldots, \ve^{(d)}) \geq \frac{t}{2} \right\rbrace\\
    &\leq \left( \frac{2d}{\log(3/2)}\right)^{\sum_{i=1}^d n_i} 2 \exp\left( - \frac{t^2}{8} \right),
\end{align*}
which yields the final bound for an appropriate choice of $t$.
\end{proof}
}
\begin{lemma}[Woodbury matrix identity]\label{lem:woodbury} Let $\mA\in \sM_n$, $\mB\in \sM_k$, $\mU\in \sM_{n,k}$ and $\mV\in \sM_{k,n}$, we have
\begin{align*}
\left( \mA + \mU \mB \mV \right)^{-1} = \mA^{-1} - \mA^{-1} \mU \left( \mB^{-1} + \mV \mA^{-1} \mU \right)^{-1} \mV \mA^{-1}.
\end{align*}
\end{lemma}

\bibliographystyle{imsart-number} 
\bibliography{references}       


\end{document}

%% file: math.tex
\usepackage{amsmath,amsfonts,bm}

\def\eqref#1{equation~\ref{#1}}

\def\1{\bm{1}}

\def\vzero{{\bm{0}}}

\def\va{{\bm{a}}}
\def\vb{{\bm{b}}}
\def\vc{{\bm{c}}}
\def\vd{{\bm{d}}}
\def\ve{{\bm{e}}}

\def\vg{{\bm{g}}}
\def\vh{{\bm{h}}}

\def\vk{{\bm{k}}}

\def\vu{{\bm{u}}}
\def\vv{{\bm{v}}}
\def\vw{{\bm{w}}}
\def\vx{{\bm{x}}}
\def\vy{{\bm{y}}}
\def\vz{{\bm{z}}}

\def\vdelta{{\bm{\delta}}}

\def\mA{{\bm{A}}}
\def\mB{{\bm{B}}}
\def\mC{{\bm{C}}}

\def\mI{{\bm{I}}}

\def\mM{{\bm{M}}}
\def\mN{{\bm{N}}}
\def\mO{{\bm{O}}}

\def\mQ{{\bm{Q}}}
\def\mR{{\bm{R}}}
\def\mS{{\bm{S}}}
\def\mT{{\bm{T}}}
\def\mU{{\bm{U}}}
\def\mV{{\bm{V}}}

\def\mX{{\bm{X}}}

\def\mZ{{\bm{Z}}}

\DeclareMathAlphabet{\mathsfit}{\encodingdefault}{\sfdefault}{m}{sl}
\SetMathAlphabet{\mathsfit}{bold}{\encodingdefault}{\sfdefault}{bx}{n}
\newcommand{\tens}[1]{\bm{\mathsfit{#1}}}
\def\tA{{\tens{A}}}
\def\tB{{\tens{B}}}
\def\tC{{\tens{C}}}

\def\tT{{\tens{T}}}

\def\tW{{\tens{W}}}
\def\tX{{\tens{X}}}
\def\tY{{\tens{Y}}}

\def\gE{{\mathcal{E}}}

\def\sC{{\mathbb{C}}}

\def\sM{{\mathbb{M}}}

\def\sP{{\mathbb{P}}}

\def\sR{{\mathbb{R}}}
\def\sS{{\mathbb{S}}}
\def\sT{{\mathbb{T}}}

\newcommand{\Var}{\mathrm{Var}}

\DeclareMathOperator*{\argmin}{arg\,min}

\DeclareMathOperator{\sign}{sign}

%% file: aap-sample.bbl
\begin{thebibliography}{25}

\bibitem{anandkumar2014tensor}
\begin{barticle}[author]
\bauthor{\bsnm{Anandkumar},~\bfnm{Animashree}\binits{A.}},
  \bauthor{\bsnm{Ge},~\bfnm{Rong}\binits{R.}},
  \bauthor{\bsnm{Hsu},~\bfnm{Daniel}\binits{D.}},
  \bauthor{\bsnm{Kakade},~\bfnm{Sham~M}\binits{S.~M.}} \AND
  \bauthor{\bsnm{Telgarsky},~\bfnm{Matus}\binits{M.}}
(\byear{2014}).
\btitle{Tensor decompositions for learning latent variable models}.
\bjournal{Journal of machine learning research}
\bvolume{15}
\bpages{2773--2832}.
\end{barticle}
\endbibitem

\bibitem{baik2005phase}
\begin{barticle}[author]
\bauthor{\bsnm{Baik},~\bfnm{Jinho}\binits{J.}},
  \bauthor{\bsnm{Arous},~\bfnm{G{\'e}rard~Ben}\binits{G.~B.}} \AND
  \bauthor{\bsnm{P{\'e}ch{\'e}},~\bfnm{Sandrine}\binits{S.}}
(\byear{2005}).
\btitle{Phase transition of the largest eigenvalue for nonnull complex sample
  covariance matrices}.
\bjournal{The Annals of Probability}
\bvolume{33}
\bpages{1643--1697}.
\end{barticle}
\endbibitem

\bibitem{arous2021long}
\begin{barticle}[author]
\bauthor{\bsnm{Ben~Arous},~\bfnm{G{\'e}rard}\binits{G.}},
  \bauthor{\bsnm{Huang},~\bfnm{Daniel~Zhengyu}\binits{D.~Z.}} \AND
  \bauthor{\bsnm{Huang},~\bfnm{Jiaoyang}\binits{J.}}
(\byear{2021}).
\btitle{Long Random Matrices and Tensor Unfolding}.
\bjournal{arXiv preprint arXiv:2110.10210}.
\end{barticle}
\endbibitem

\bibitem{arous2019landscape}
\begin{barticle}[author]
\bauthor{\bsnm{Ben~Arous},~\bfnm{Gerard}\binits{G.}},
  \bauthor{\bsnm{Mei},~\bfnm{Song}\binits{S.}},
  \bauthor{\bsnm{Montanari},~\bfnm{Andrea}\binits{A.}} \AND
  \bauthor{\bsnm{Nica},~\bfnm{Mihai}\binits{M.}}
(\byear{2019}).
\btitle{The landscape of the spiked tensor model}.
\bjournal{Communications on Pure and Applied Mathematics}
\bvolume{72}
\bpages{2282--2330}.
\end{barticle}
\endbibitem

\bibitem{benaych2011eigenvalues}
\begin{barticle}[author]
\bauthor{\bsnm{Benaych-Georges},~\bfnm{Florent}\binits{F.}} \AND
  \bauthor{\bsnm{Nadakuditi},~\bfnm{Raj~Rao}\binits{R.~R.}}
(\byear{2011}).
\btitle{The eigenvalues and eigenvectors of finite, low rank perturbations of
  large random matrices}.
\bjournal{Advances in Mathematics}
\bvolume{227}
\bpages{494--521}.
\end{barticle}
\endbibitem

\bibitem{biroli2020iron}
\begin{barticle}[author]
\bauthor{\bsnm{Biroli},~\bfnm{Giulio}\binits{G.}},
  \bauthor{\bsnm{Cammarota},~\bfnm{Chiara}\binits{C.}} \AND
  \bauthor{\bsnm{Ricci-Tersenghi},~\bfnm{Federico}\binits{F.}}
(\byear{2020}).
\btitle{How to iron out rough landscapes and get optimal performances: averaged
  gradient descent and its application to tensor PCA}.
\bjournal{Journal of Physics A: Mathematical and Theoretical}
\bvolume{53}
\bpages{174003}.
\end{barticle}
\endbibitem

\bibitem{capitaine2009largest}
\begin{barticle}[author]
\bauthor{\bsnm{Capitaine},~\bfnm{Mireille}\binits{M.}},
  \bauthor{\bsnm{Donati-Martin},~\bfnm{Catherine}\binits{C.}} \AND
  \bauthor{\bsnm{F{\'e}ral},~\bfnm{Delphine}\binits{D.}}
(\byear{2009}).
\btitle{The largest eigenvalues of finite rank deformation of large Wigner
  matrices: convergence and nonuniversality of the fluctuations}.
\bjournal{The Annals of Probability}
\bvolume{37}
\bpages{1--47}.
\end{barticle}
\endbibitem

\bibitem{de2021random}
\begin{barticle}[author]
\bauthor{\bsnm{Goulart},~\bfnm{Jos{\'e} Henrique de~Morais}\binits{J.~H.
  d.~M.}}, \bauthor{\bsnm{Couillet},~\bfnm{Romain}\binits{R.}} \AND
  \bauthor{\bsnm{Comon},~\bfnm{Pierre}\binits{P.}}
(\byear{2021}).
\btitle{A Random Matrix Perspective on Random Tensors}.
\bjournal{arXiv preprint arXiv:2108.00774}.
\end{barticle}
\endbibitem

\bibitem{handschy2019phase}
\begin{bphdthesis}[author]
\bauthor{\bsnm{Handschy},~\bfnm{Madeline~Curtis}\binits{M.~C.}}
(\byear{2019}).
\btitle{Phase Transition in Random Tensors with Multiple Spikes},
\btype{PhD thesis},
\bpublisher{University of Minnesota}.
\end{bphdthesis}
\endbibitem

\bibitem{hitchcock1927expression}
\begin{barticle}[author]
\bauthor{\bsnm{Hitchcock},~\bfnm{Frank~L}\binits{F.~L.}}
(\byear{1927}).
\btitle{The expression of a tensor or a polyadic as a sum of products}.
\bjournal{Journal of Mathematics and Physics}
\bvolume{6}
\bpages{164--189}.
\end{barticle}
\endbibitem

\bibitem{huang2020power}
\begin{barticle}[author]
\bauthor{\bsnm{Huang},~\bfnm{Jiaoyang}\binits{J.}},
  \bauthor{\bsnm{Huang},~\bfnm{Daniel~Z}\binits{D.~Z.}},
  \bauthor{\bsnm{Yang},~\bfnm{Qing}\binits{Q.}} \AND
  \bauthor{\bsnm{Cheng},~\bfnm{Guang}\binits{G.}}
(\byear{2020}).
\btitle{Power iteration for tensor pca}.
\bjournal{arXiv preprint arXiv:2012.13669}.
\end{barticle}
\endbibitem

\bibitem{jagannath2020statistical}
\begin{barticle}[author]
\bauthor{\bsnm{Jagannath},~\bfnm{Aukosh}\binits{A.}},
  \bauthor{\bsnm{Lopatto},~\bfnm{Patrick}\binits{P.}} \AND
  \bauthor{\bsnm{Miolane},~\bfnm{Leo}\binits{L.}}
(\byear{2020}).
\btitle{Statistical thresholds for tensor {PCA}}.
\bjournal{The Annals of Applied Probability}
\bvolume{30}
\bpages{1910--1933}.
\end{barticle}
\endbibitem

\bibitem{landsberg2012tensors}
\begin{barticle}[author]
\bauthor{\bsnm{Landsberg},~\bfnm{Joseph~M}\binits{J.~M.}}
(\byear{2012}).
\btitle{Tensors: geometry and applications}.
\bjournal{Representation theory}
\bvolume{381}
\bpages{3}.
\end{barticle}
\endbibitem

\bibitem{lesieur2017statistical}
\begin{binproceedings}[author]
\bauthor{\bsnm{Lesieur},~\bfnm{Thibault}\binits{T.}},
  \bauthor{\bsnm{Miolane},~\bfnm{L{\'e}o}\binits{L.}},
  \bauthor{\bsnm{Lelarge},~\bfnm{Marc}\binits{M.}},
  \bauthor{\bsnm{Krzakala},~\bfnm{Florent}\binits{F.}} \AND
  \bauthor{\bsnm{Zdeborov{\'a}},~\bfnm{Lenka}\binits{L.}}
(\byear{2017}).
\btitle{Statistical and computational phase transitions in spiked tensor
  estimation}.
In \bbooktitle{Proc. IEEE International Symposium on Information Theory (ISIT)}
\bpages{511--515}.
\end{binproceedings}
\endbibitem

\bibitem{lim2005singular}
\begin{binproceedings}[author]
\bauthor{\bsnm{Lim},~\bfnm{Lek-Heng}\binits{L.-H.}}
(\byear{2005}).
\btitle{Singular values and eigenvalues of tensors: a variational approach}.
In \bbooktitle{Proc. IEEE International Workshop on Computational Advances in
  Multi-Sensor Adaptive Processing}
\bpages{129--132}.
\end{binproceedings}
\endbibitem

\bibitem{marvcenko1967distribution}
\begin{barticle}[author]
\bauthor{\bsnm{Mar{\v{c}}enko},~\bfnm{Vladimir~A}\binits{V.~A.}} \AND
  \bauthor{\bsnm{Pastur},~\bfnm{Leonid~Andreevich}\binits{L.~A.}}
(\byear{1967}).
\btitle{Distribution of eigenvalues for some sets of random matrices}.
\bjournal{Mathematics of the USSR-Sbornik}
\bvolume{1}
\bpages{457}.
\end{barticle}
\endbibitem

\bibitem{montanari2014statistical}
\begin{barticle}[author]
\bauthor{\bsnm{Montanari},~\bfnm{Andrea}\binits{A.}} \AND
  \bauthor{\bsnm{Richard},~\bfnm{Emile}\binits{E.}}
(\byear{2014}).
\btitle{A statistical model for tensor {PCA}}.
\bjournal{arXiv preprint arXiv:1411.1076}.
\end{barticle}
\endbibitem

\bibitem{nocedal2006numerical}
\begin{bbook}[author]
\bauthor{\bsnm{Nocedal},~\bfnm{Jorge}\binits{J.}} \AND
  \bauthor{\bsnm{Wright},~\bfnm{Stephen}\binits{S.}}
(\byear{2006}).
\btitle{Numerical optimization}.
\bpublisher{Springer Science \& Business Media}.
\end{bbook}
\endbibitem

\bibitem{peche2006largest}
\begin{barticle}[author]
\bauthor{\bsnm{P{\'e}ch{\'e}},~\bfnm{Sandrine}\binits{S.}}
(\byear{2006}).
\btitle{The largest eigenvalue of small rank perturbations of Hermitian random
  matrices}.
\bjournal{Probability Theory and Related Fields}
\bvolume{134}
\bpages{127--173}.
\end{barticle}
\endbibitem

\bibitem{perry2020statistical}
\begin{binproceedings}[author]
\bauthor{\bsnm{Perry},~\bfnm{Amelia}\binits{A.}},
  \bauthor{\bsnm{Wein},~\bfnm{Alexander~S}\binits{A.~S.}} \AND
  \bauthor{\bsnm{Bandeira},~\bfnm{Afonso~S}\binits{A.~S.}}
(\byear{2020}).
\btitle{Statistical limits of spiked tensor models}.
In \bbooktitle{Annales de l'Institut Henri Poincar{\'e}, Probabilit{\'e}s et
  Statistiques}
\bvolume{56}
\bpages{230--264}.
\bpublisher{Institut Henri Poincar{\'e}}.
\end{binproceedings}
\endbibitem

\bibitem{rabanser2017introduction}
\begin{barticle}[author]
\bauthor{\bsnm{Rabanser},~\bfnm{Stephan}\binits{S.}},
  \bauthor{\bsnm{Shchur},~\bfnm{Oleksandr}\binits{O.}} \AND
  \bauthor{\bsnm{G{\"u}nnemann},~\bfnm{Stephan}\binits{S.}}
(\byear{2017}).
\btitle{Introduction to tensor decompositions and their applications in machine
  learning}.
\bjournal{arXiv preprint arXiv:1711.10781}.
\end{barticle}
\endbibitem

\bibitem{stein1981estimation}
\begin{barticle}[author]
\bauthor{\bsnm{Stein},~\bfnm{Charles~M}\binits{C.~M.}}
(\byear{1981}).
\btitle{Estimation of the mean of a multivariate normal distribution}.
\bjournal{The annals of Statistics}
\bpages{1135--1151}.
\end{barticle}
\endbibitem

\bibitem{sun2014tensors}
\begin{barticle}[author]
\bauthor{\bsnm{Sun},~\bfnm{Will~Wei}\binits{W.~W.}},
  \bauthor{\bsnm{Hao},~\bfnm{Botao}\binits{B.}} \AND
  \bauthor{\bsnm{Li},~\bfnm{Lexin}\binits{L.}}
(\byear{2014}).
\btitle{Tensors in Modern Statistical Learning}.
\bjournal{Wiley StatsRef: Statistics Reference Online}
\bpages{1--25}.
\end{barticle}
\endbibitem

\bibitem{tao2012topics}
\begin{bbook}[author]
\bauthor{\bsnm{Tao},~\bfnm{Terence}\binits{T.}}
(\byear{2012}).
\btitle{Topics in random matrix theory}
\bvolume{132}.
\bpublisher{American Mathematical Soc.}
\end{bbook}
\endbibitem

\bibitem{zare2018extension}
\begin{barticle}[author]
\bauthor{\bsnm{Zare},~\bfnm{Ali}\binits{A.}},
  \bauthor{\bsnm{Ozdemir},~\bfnm{Alp}\binits{A.}},
  \bauthor{\bsnm{Iwen},~\bfnm{Mark~A}\binits{M.~A.}} \AND
  \bauthor{\bsnm{Aviyente},~\bfnm{Selin}\binits{S.}}
(\byear{2018}).
\btitle{Extension of PCA to higher order data structures: An introduction to
  tensors, tensor decompositions, and tensor PCA}.
\bjournal{Proceedings of the IEEE}
\bvolume{106}
\bpages{1341--1358}.
\end{barticle}
\endbibitem

\end{thebibliography}
